\documentclass[reqno, 11pt, a4paper]{amsart}
\usepackage{amsmath, amssymb, amsthm}

\usepackage{appendix}
\usepackage{graphicx, color}
\usepackage{enumerate}
\numberwithin{equation}{section}
\usepackage{a4wide}
\usepackage[active]{srcltx}
\usepackage[table]{xcolor}
\definecolor{cyan(process)}{rgb}{0.0, 0.72, 0.92}
\usepackage{xcolor}
\definecolor{columbiablue}{rgb}{0.61, 0.87, 1.0}
\usepackage{wrapfig}
\usepackage[]{pstricks}
\usepackage{fancyhdr}
\definecolor{sandstone}{HTML}{786D5F}
\definecolor{beaublue}{rgb}{0.74, 0.83, 0.9}
\definecolor{cherryblossompink}{rgb}{1.0, 0.72, 0.77}
\definecolor{light-gray}{gray}{0.95}

\newcommand{\kh}[1]{{\color{red} #1}}
\definecolor{apricot}{rgb}{0.98, 0.81, 0.69}

\usepackage[pagewise, displaymath, mathlines]{lineno}
\newcommand*\patchAmsMathEnvironmentForLineno[1]{
  \expandafter\let\csname old#1\expandafter\endcsname\csname #1\endcsname
  \expandafter\let\csname oldend#1\expandafter\endcsname\csname end#1\endcsname
  \renewenvironment{#1}
     {\linenomath\csname old#1\endcsname}
     {\csname oldend#1\endcsname\endlinenomath}}
\newcommand*\patchBothAmsMathEnvironmentsForLineno[1]{
  \patchAmsMathEnvironmentForLineno{#1}
  \patchAmsMathEnvironmentForLineno{#1*}}
\AtBeginDocument{
\patchBothAmsMathEnvironmentsForLineno{equation}
\patchBothAmsMathEnvironmentsForLineno{align}
\patchBothAmsMathEnvironmentsForLineno{flalign}
\patchBothAmsMathEnvironmentsForLineno{alignat}
\patchBothAmsMathEnvironmentsForLineno{gather}
\patchBothAmsMathEnvironmentsForLineno{multline}
}

\newcommand{\Add}[1]{\textcolor{black}{#1}}	
\newcommand{\Erase}[1]{\if0{#1}\fi}

\usepackage{mathtools}
\usepackage{mathrsfs, enumerate}
\usepackage{hyperref}
\usepackage[truedimen,top=3truecm, bottom=3truecm, left=2.5truecm, right=2.5truecm, includefoot]{geometry} 

\usepackage{tikz}
\usetikzlibrary{shapes,arrows,positioning,fit,backgrounds,calc}
\usepackage{pgfplots}
\pgfplotsset{compat=1.16}

\usepackage[percent]{overpic}

\begin{document}

\title[Anomalous behavior from oscillators in the high-temperature regime]{Derivation of anomalous behavior from Interacting Oscillators in the High-Temperature Regime}

\author[P. Gon\c{c}alves]{Patr\'{i}cia Gon\c{c}alves}
\address{Center for Mathematical Analysis, Geometry and Dynamical Systems, Instituto Superior Técnico, Universidade de Lisboa, Av. Rovisco Pais, 1049-001 Lisboa, Portugal.}
\email{pgoncalves@tecnico.ulisboa.pt}

\author[K. Hayashi]{Kohei Hayashi}
\address{RIKEN Interdisciplinary Theoretical and Mathematical Science, 2-1 Hirosawa, Wako, Saitama 351-0198, Japan.}
\email{koheihayashi0826@gmail.com}


\maketitle

\theoremstyle{plain}
\newtheorem{theorem}{Theorem}[section] 
\newtheorem{lemma}[theorem]{Lemma}
\newtheorem{corollary}[theorem]{Corollary}
\newtheorem{proposition}[theorem]{Proposition}

\theoremstyle{definition}
\newtheorem{definition}[theorem]{Definition}
\newtheorem{remark}[theorem]{Remark}
\newtheorem{assumption}[theorem]{Assumption}
\newtheorem{example}[theorem]{Example}

\makeatletter
\renewcommand{\theequation}{%
\thesection.\arabic{equation}}
\@addtoreset{equation}{section}
\makeatother

\makeatletter
\renewcommand{\p@enumi}{A}
\makeatother

\begin{abstract}
A microscopic model of interacting oscillators, which admits two conserved quantities, volume, and energy, is investigated.
We begin with a system driven by a general nonlinear potential under high-temperature regime by taking the inverse temperature of the system asymptotically small. 
As a consequence, one can extract a principal part (by a simple Taylor expansion argument), which is driven by the harmonic potential,  and we show that previous results for the harmonic chain are covered with generality. 
We consider two fluctuation fields, which are defined as a linear combination of {the fluctuation fields of} the two conserved quantities, volume, and energy, and we show that the fluctuations of one field converge to a solution of an additive stochastic heat equation, which corresponds to the Ornstein-Uhlenbeck process, in a weak asymmetric regime, or to a solution of the stochastic Burgers equation, in a stronger asymmetric regime. On the other hand, the fluctuations of the other field cross from an additive stochastic heat equation to a fractional diffusion equation given by a skewed L\'evy process. 
\end{abstract}

\section{Introduction}
A long standing problem in the field of mathematical physics has to do with the rigorous mathematical derivation of the macroscopic evolution equations of the conserved quantities in Newtonian particle systems. 
By assuming a  stochastic dynamics instead of a deterministic one, the problem becomes mathematically tractable and  some answers can be given with success. 
Over the last four decades both mathematicians and theoretical physicists have been giving many fruitful and relevant contributions related to the aforementioned problem and one of the challenges that has been in the spotlight is the derivation of the well-known hydrodynamic limit from stochastic interacting particle systems, as well as, the characterization of the fluctuations of locally conserved quantities around that limit. 
In this framework, many types of partial differential equations (PDE) and stochastic PDEs~(SPDEs) have been studied and derived from several underlying random dynamics. 
These equation can be of different nature depending on the type of the underlying dynamics. 
The hydrodynamic limit consists in showing that the empirical measure associated to each conserved quantity converges to a deterministic measure which is absolutely continuous with respect to the Lebesgue measure and whose density is a solution to a PDE, the hydrodynamic equation. 
Since the limit is deterministic, probabilistically speaking, the hydrodynamic limit is a law of large numbers for the conserved quantity(ies) of the system, whereas the fluctuations is a central limit theorem, since the limit is random and  described by a solution to an SPDE.  
These two derivations are done through a scaling limit procedure where the scaling parameter $n$ connects the macroscopic space, a continuous space where the solutions of the macroscopic equations will be defined; to the microscopic space, a  discrete space where the random system evolves according to a prescribed dynamics.  
Since two scales for space are considered, two scales for time naturally emerge, a macroscopic time $t$ and a microscopic time $tn^a$, where the value of the parameter $a$ highly depends on the underlying dynamics. 
For microscopic systems with only one conservation law, there is no ambiguity on the choice of the fluctuation fields that one should look at. 
Nevertheless, for multi-component systems one has many ways of considering the fluctuation fields associated to the conserved quantities and moreover a special feature of these systems is the fact that different time scales might coexist, and this never occurs for systems with only one conservation law.

In \cite{spohn2015nonlinear}, with a focus on anharmonic chains of oscillators, it was developed the nonlinear fluctuating hydrodynamics theory (NLFH) for the equilibrium time-correlations of the conserved quantities of that model (which has several conservation laws) and analytical predictions were done based on a mode-coupling approximation. 
Roughly speaking, the approach of \cite{spohn2015nonlinear} starts at the macroscopic level, i.e. assumes that a hyperbolic system of conservation laws is governing the macroscopic evolution of the empirical conserved quantities. 
Then one adds a diffusion term and a dissipation term to the system of coupled PDEs and then linearizes the system at second order, with respect to the equilibrium averages of the conserved quantities. 
A fundamental role is played by the normal modes, i.e. the eigenvectors of the linearized equation and these modes evolve with different velocities and in different time scales. 
These modes might be described by different forms of anomalous super-diffusion or standard diffusions and this description depends on the value of certain coupling constants. 
These coupling constants fix the value of the quadratic terms in the equation, i.e. terms in the evolution equation of each quantity which are written in terms of products of the conserved quantities and these constants fix the limiting processes that one should obtain. 
For systems with two conservation laws all the possible limits are summarized in the tables in the appendix, see Tables \ref{tab:mcm_1}, \ref{tab:mcm_2} and \ref{tab:mcm_3}.  
Surprisingly, besides the usual diffusive behavior, several forms of anomalous behavior can be obtained either by means of  fractional behavior (but for very particular values of fractional {power})  or the Kardar-Parisi-Zhang (KPZ) behavior.

The KPZ behavior that has been extensively derived from microscopic systems and it is the one that we refer to in this article is characterized through what is called the KPZ equation or its companion the stochastic Burgers (SB) equation that we now briefly describe.  
The KPZ equation, which was introduced in ~\cite{kardar1986dynamic} is the following SPDE
\begin{equation*}
\partial_t h = \nu \partial_x^2 h
+ \tau (\partial_x h)^2 + \sqrt{D} \dot{W}.
\end{equation*}
Above $t\in[0,\infty)$ and $x\in\mathbb{R}$, denote temporal and spatial variables, respectively. 
In addition, $\nu,D>0$ and $\tau\in \mathbb{R}$ are some constants and $\dot{W}=\dot{W}(t,x)$ denotes the one-dimensional space-time white-noise. 
The KPZ equation is conjectured to be a universal SPDE describing the fluctuations of randomly growing interfaces of 1-d stochastic dynamics close to a stationary state.
Throughout this paper, we focus on the one-dimensional case for space.  
As a similar object, the tilt {$u=\partial_xh$} formally satisfies the following stochastic Burgers (SB) equation.  
\begin{equation*}
\partial_t u = \nu \partial_x^2 u
+ \tau \partial_x u^2 + \sqrt{D} \partial_x \dot{W}.
\end{equation*}
What is notable for the above KPZ/SB equations is their universality, since they have been derived from various types of microscopic systems~\cite{bertini1997stochastic, gonccalves2014nonlinear, gonccalves2015stochastic, hayashi2023derivation, diehl2017kardar, jara2019scaling, jara2020stationary}.
Whereas these results are concerned with the scalar-valued case, a vector-valued version, which macroscopically corresponds to a coupled system of KPZ/SB equations, has also been studied~\cite{bernardin2021derivation, butelmann2022scaling, hayashi2022derivation} for a system of interacting diffusion processes.

Above we have mentioned that fractional behavior has been predicted from NLFH theory for the normal modes of systems with several conservation laws. 
According to the predictions of \cite{spohn2015nonlinear}, the fractional behavior in systems with two conservation is given by a $\alpha$-L\'evy process with $\alpha\in\{3/2,5/3, \textrm{Gold}\}$ where $ \textrm{Gold}=(1+\sqrt 5)/2$ is the Golden mean number. 
The fact that only these exponents appear in the limit is still not well understood.

In this paper, we consider a model of interacting oscillators which was introduced in~\cite{bernardin2012anomalous} that we refer to as the BS model. 
The dynamics results from a superposition of an Hamiltonian dynamics  (depending on a given potential, say $V:\mathbb{R}\to \mathbb{R}$, and whose strength is regulated by a sequence $\alpha_n$) and a noise that exchanges variables at nearest-neighbor positions. 
This model admits two conserved quantities, volume and energy, and anomalous diffusion phenomena have been observed from these quantities. 
The paper~\cite{bernardin2014anomalous} pointed out that the model exhibits anomalous behavior  when the potential is given by the so-called Toda lattice potential~(a.k.a.~the Kac-van Moerbecke potential) $V(\eta)=e^{-\eta}-1+\eta$. 
Additionally, \cite{bernardin20163, bernardin2018weakly} identified the limiting process of the energy fluctuation field, for the case of the harmonic potential $V(\eta)=\eta^2/2$, by a $3/2$-L\'{e}vy, while the volume fluctuation field has diffusive behavior. 
In fact, the limit behavior of the energy fluctuation field is diffusive in the regime where the noise is stronger than the Hamiltonian dynamics but when the Hamiltonian dynamics is stronger, then the behavior is described by a skewed $3/2$-fractional L\'{e}vy process, see {Figure \ref{fig:energy}} where it is assumed that $\alpha_n=O(n^{-\kappa})$ and the system is speeded up in the time scale $tn^a$: 

\begin{figure}[htb!]
\begin{center}
\begin{tikzpicture}[scale=0.17]
\draw (0,25) node[left]{$a$};
\draw (25,0) node[below]{$\kappa$};
\draw (8,0) node[below]{$ \frac 13$};
\draw (0,0) node[left]{$0$};
\draw (0,11) node[left]{$\frac 32$};
\draw (0,20) node[left]{$2$};
\fill[light-gray] (0,11) -- (8,20) -- (25,20) -- (25,25) -- (0,25) -- cycle;
\fill[fill=pink, fill opacity=0.5] (0,0) -- (25,0) -- (25,20) -- (8,20) -- (0,11)--cycle;
\draw[-,=latex,blue!50!green,ultra thick] (8,20) -- (25, 20) node[midway,above,sloped] {\textbf{ \tiny{Diffusion} }};
\draw[-,=latex,gray,red, ultra thick] (0,11) -- (8,20) node[midway,above,sloped] {\textbf{\tiny{ L\'evy $\tfrac 32$}}};
\node[circle,fill=black,inner sep=0.8mm] at (8,20) {};
\node[] at (8,21) [above] {\textbf{\tiny{L\'evy $\tfrac 32$+Diffusion}}};
\node[] at (17.5,10) {\textbf{\tiny{no evolution}}};
\draw[-,=latex, dashed] (8,-0.1) -- (8,20);
\draw[->,>=latex] (0,0) -- (26,0);
\draw[->,>=latex] (0,0) -- (0,26);
\end{tikzpicture}
\end{center}
\caption{ $V(\eta)$ (Energy) fluctuations.}
\label{fig:energy}
\end{figure}
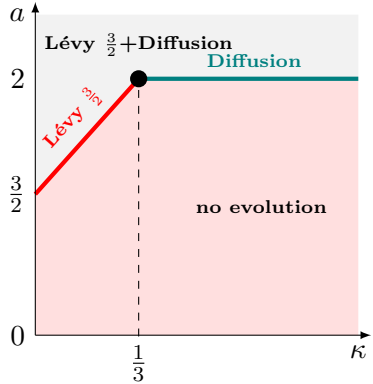
In \cite{bernardin2018nonlinear} it was proved that the last result is also valid for an anharmonic setting when the harmonic potential is perturbed by a weak quartic function, i.e. for $V_n(\eta)=\eta^2/2+\gamma_n \eta^4/4$ and for $\gamma_n=O(n^{-1/4})$.

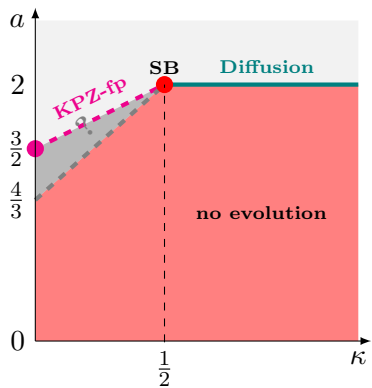
\begin{figure}[htb!]
\begin{center}
\begin{tikzpicture}[scale=0.17]
\draw (0,25) node[left]{$a$};
\draw (25,0) node[below]{$\kappa$};
\draw (10,0) node[below]{$\frac 12$};
\draw (0,0) node[left]{$0$};
\draw (0,15) node[left]{$\frac 32$};
\draw (0,11) node[left]{$\frac 43$};
\draw (0,20) node[left]{$2$};
\fill[light-gray] (0,11) -- (10,20) -- (25,20) -- (25,25) -- (0,25) -- cycle;
\fill[fill=red, fill opacity=0.5] (0,0) -- (25,0) -- (25,20) -- (10,20) -- (0,11)--cycle;
\fill[fill=gray, fill opacity=0.5] (0,11) -- (10,20) -- (0,15)--cycle;
\draw[-,=latex,blue!50!green,ultra thick] (10,20) -- (25, 20) node[midway,above,sloped] {\textbf{ \tiny{Diffusion}}};
\draw[-,=latex,magenta,dashed, ultra thick] (0,15) -- (10,20) node[midway,above,sloped] {\textbf {\tiny{KPZ-fp }}};
\draw[-,=latex,gray,dashed, ultra thick] (0,11) -- (10,20) node[midway,above,sloped] {\textbf ?};
\node[circle,fill=red,inner sep=0.8mm] at (10,20) {};
\node[circle,fill=magenta,inner sep=0.8mm] at (0,15) {};
\node[] at (10,20) [above] {\textbf{\tiny{SB}}};
\node[] at (17.5,10) {\textbf{\tiny{no evolution}}};
\draw[-,=latex, dashed] (10,-0.1) -- (10,20);
\draw[->,>=latex] (0,0) -- (26,0);

\draw[->,>=latex] (0,0) -- (0,26);
\end{tikzpicture}
\end{center}
\caption{ $V'(\eta)$ fluctuations.}
\label{fig:KvM}
\end{figure}

More recently,  \cite{ahmed2022microscopic} analyzed, in the case of the Toda lattice potential, the limiting behavior of the  fluctuation field for another conserved quantity, which corresponds to the derivative of the energy, i.e. for $V'(\eta)$, and they derived the SB equation as the limiting object when the Hamiltonian dynamics is as strong as the exchange noise dynamics, or a diffusive behavior  when the noise is stronger. 
The results for this quantity can be summarized in {Figure \ref{fig:KvM}}. 
In the regime where the Hamiltonian dynamics is stronger it is believed that the behavior of the quantity $V'(\eta)$ should be given in terms of the so-called KPZ fixed point (KPZ-fp), constructed in \cite{matetski2021kpz}.  
Nevertheless, for the other normal mode, only diffusive behavior can be obtained, but in the regime where the intensity of the noise is stronger than the Hamiltonian dynamics, and all the other regimes are still open. 

As we have seen above, depending on the chosen potential one can either get the limit behavior for the conserved quantities as diffusive, KPZ or a $3/2$-L\'{e}vy. 
With this in mind, our motivation in this article is to see how universal is the aforementioned behavior. 
Therefore, in this paper, we study the above model of interacting oscillators which is driven by a general form of the potential $V$ under the high-temperature regime. 
To clarify our idea, let $\beta>0$ be the inverse temperature of the system and consider the infinite temperature limit $\beta\to0$.
Then, a rescaled potential $V_\beta(\eta)=\beta^{-2}V(\beta\eta)$ is expanded as 
\begin{equation}\label{eq:taylor_expansion_introduction}
V_\beta(\eta) = \frac{1}{2} V^{\prime \prime}(0) \eta^2 
+ \frac{1}{6} \beta V^{(3)}(0) \Add{\eta^3}
+ \cdots,
\end{equation}
provided that the potential satisfies a normalizing condition $V(0)=V^\prime(0)=0$.
From this simple Taylor expansion argument, we can extract the harmonic potential as a principal object with a small error term when $\beta\to0$. 
From this we expect that the previous results can be covered with more generality since the above argument holds for an almost arbitrary form of the potential, as long as, it satisfies some regularity condition, for instance. 
Moreover, note that the potential possibly has a nonlinear perturbation which is proportional to $\beta$, or in other words we can call it \textit{asymmetry} since it comes from a cubic function. 
From this observation, it is also expected that results that come from a nonlinearity are also re-derived from our model. 
Our main results state that these two pieces of anticipation are indeed justified.

In this paper, we consider the limiting behavior of volume and energy fluctuation fields.
A natural object to consider is the pair of volume and energy fluctuation fields themselves.  
To obtain a nonlinear term in the limit, however, we are obliged to take moving frames with different values of velocities for the two fields. 
This makes it difficult to obtain a closed system of limiting equations and only a result for linear fluctuations is obtained~(see Theorem \ref{thm:pair_fluctuations}). 
Instead, we study a linear combination of the volume and energy fluctuation fields, i.e. the field associated to the quantity $\eta+\mathfrak{u}V(\eta)$ ($\mathfrak{u}$ is a constant to be fixed) which is a scalar-valued field, with a common value of velocity. We are then in the context of NLFH theory and we can consider the fields of the normal modes associated to the system. 
Then, it turns out that we have two possible choices of linear combinations, i.e. two choices for the value of $\mathfrak{u}$  and the corresponding values of velocity, whose corresponding fluctuation fields converge either to  solutions of the Ornstein-Uhlenbeck equation, in the weak asymmetry regime; or to the SB equation or a $3/2$-L\'{e}vy fractional diffusion equation, in the stronger asymmetry regime.

We could be tempted to use the NLFH theory developed in \cite{spohn2015nonlinear}, to set up the correct linear combination and respective velocity, and then to obtain the  corresponding predictions on the limiting equations for each field. 
Nevertheless, the starting point is to consider the column matrix given with the average flux of each conserved quantity. Unfortunately, we are not able to perform these computations with great generality because we are not able to rewrite these average currents in terms of the averages of the conserved quantities, and as a consequence, we cannot pose the problem correctly. 
In \cite{spohn2015nonlinear} NLFH theory was applied to the BS model for a fixed potential and despite the difficulties we just mentioned, some numerical simulations have been done and predictions for the normal modes were obtained, for more on this see the appendix. 
To stand on exact results, alternatively, what we do is to look at Dynkin's formula applied to a generic field which is written as a linear combination of the fields of energy and volume, i.e. $\eta+ \mathfrak{u} V(\eta)$ and both evolving in the same reference frame with a given value of velocity $v$. 
The choice $\mathfrak{u}$ and the velocity of reference frames $v$ is such that in Dynkin's formula, the lower order terms with respect to the conserved quantities are null. 
From this observation, we can infer what are the correct fields to look at and the respective velocities. 
Since we are concerned with two conservation laws we have two unknowns $(\mathfrak{u},v)$ to obtain from a system with two equations. 
By eliminating these diverging terms in Dynkin's formula, what remains are higher-order terms that are, in principle, easier to control. 
Nevertheless, our infinitesimal generator when acting at the conserved quantities, energy, and volume, adds more terms to the equations, and, as a consequence, in the evolution of each quantity, we get a hierarchy of equations that have to be carefully estimated and then properly truncated. 
To that end, we use two different methods.
The first one is based on the second-order Boltzmann-Gibbs principle, which is used to derive the crossover from the Ornstein-Uhlenbeck equation to the SB equation~(Theorem \ref{thm:SBE_derivation}).
The second-order Boltzmann-Gibbs principle was introduced in \cite{gonccalves2014nonlinear} to treat quadratic fields of the exclusion process, and was extended to many other models in \cite{gonccalves2015stochastic, gonccalves2017second, blondel2016convergence, gonccalves2017stochastic, gonccalves2020derivation}, for instance.  
This is the same argument that has been used to obtain the diagram in \Erase{figure}\Add{Figure} \ref{fig:KvM} for the fluctuations of $V'(\eta)$ and we obtain exactly the same behavior from a generic potential satisfying mild assumptions.

The second method we employ for the derivation of fractional behavior~(Theorem \ref{thm:3/2L_derivation}) is to consider higher-order fields of the quantities that appear in the evolution equations and to properly control the new terms \Erase{that appear }in the expansion. 
This strategy was developed in \cite{bernardin20163} to analyze anomalous behavior of energy for a harmonic potential and in \cite{bernardin2018nonlinear} for showing persistence on the universality behavior observed for the harmonic potential but in the last case, the potential is perturbed with a quartic term. 
In our specific model we consider a quadratic field associated to the quantity $V'(\eta)$, i.e. a two-dimensional field given in terms of the quantities $V'(\eta)$ and by looking at the evolution of this quantity we get new quantities which are written in terms of both volume and energy  
and that we have to truncate to see the leading terms. 
In the case of the harmonic potential the evolution of the quadratic field is quite simple since it only involves the energy and the quadratic field itself. In this sense the hierarchy is finite since we start from the energy field and in its evolution we see the quadratic field for the other quantity, namely the volume, and in the evolution of the quadratic field for the volume we get back to the energy field. Then the only difficulty is to properly link the equations of the two fields, but this can be made by means of choosing properly the test functions that solve a Poisson equation. 
In our model this can also be done but one has to be careful since other terms appear in the evolution equations. Since we can tune the value of the inverse temperature $\beta$ we are able to control the new terms that appear in the evolution equations and truncate the hierarchy in one step as for the harmonic potential.

The \Add{two }fluctuation fields that we consider, namely the fields  associated to the quantities $\eta+\mathfrak{u} V(\eta)$ for the choices $\mathfrak{u} =V^{(3)}(0)\beta$ and $\mathfrak{u} =-1/\lambda$ (where $\lambda$ is the average of $V'(\eta)$ with respect to the stationary measure) asymptotically behave as the fields for $V'(\eta)$ and $V(\eta)$ and for this reason our results then extend those of the fluctuations of $V'(\eta)$ for the Toda Lattice potential and of the fluctuations of the energy $V(\eta)$ for the harmonic potential. 
We therefore show that the behavior of the SB equation and the $3/2$-L\'evy behavior are universal.

\begin{center}
\begin{figure}[htpb]
\includegraphics[keepaspectratio, width=0.7\textwidth]{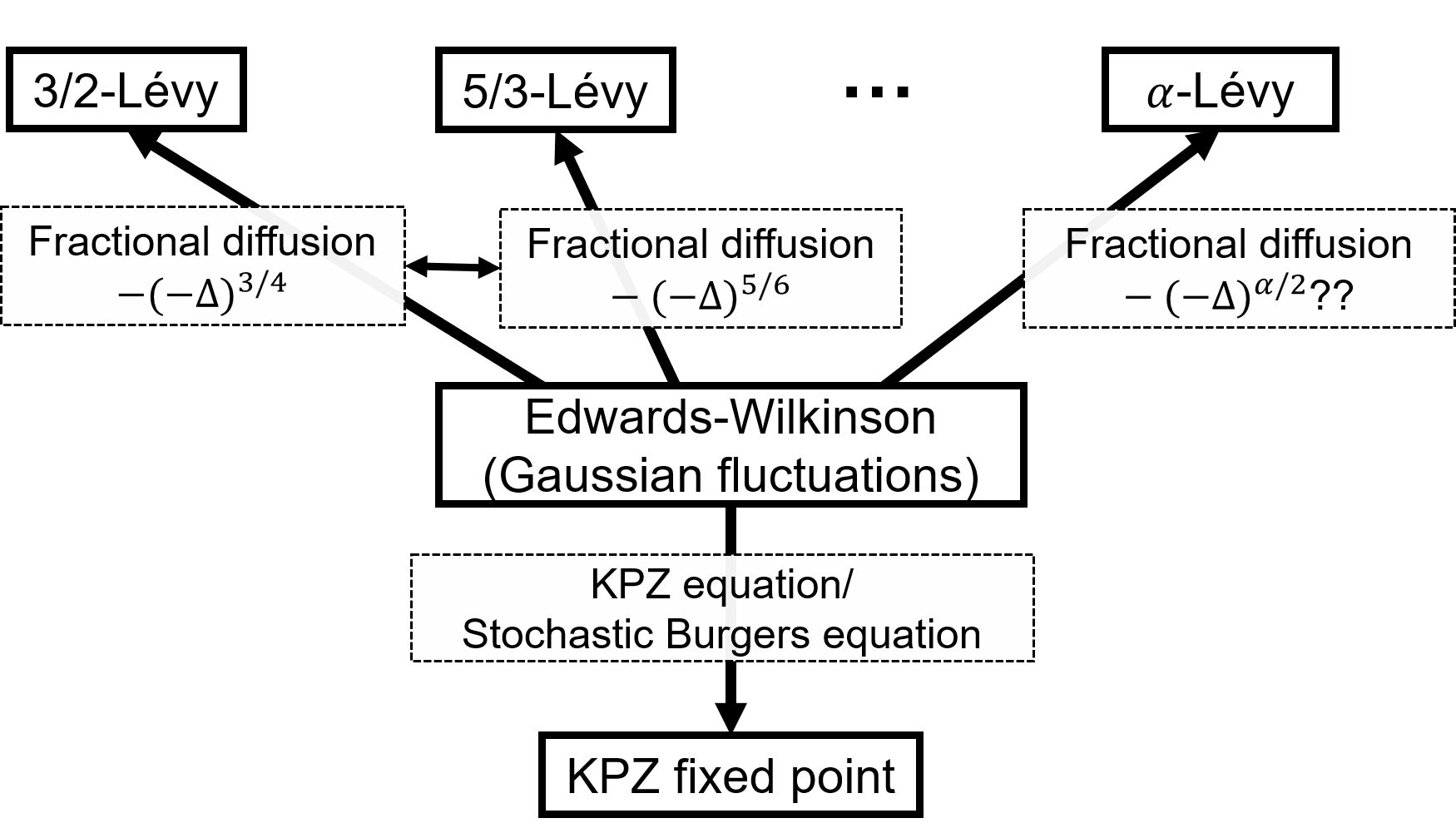} 
\caption{Classification of some universality classes. }
\label{fig:universality_classes}
\end{figure}
\end{center}

Finally, we note the appearance of fractional diffusions from many other contexts, for example, from a 1-d infinite chain of coupled\Erase{, charged} harmonic oscillators \Erase{with a magnetic field}
\Add{~\cite{cane2022superdiffusion, saito20195, komorowski2020fractional}}.
In \cite{saito20195} it is proved that the density of the energy distribution (by means of the Wigner distribution) has a space-time evolution given by a linear phonon Boltzmann equation, whose solution when properly scaled, in the limit becomes a solution of the fractional diffusion equation with exponent $5/3$. 
In \cite{cane2022superdiffusion} it is considered the same model as in \cite{saito20195} but particles are submitted to the action of a magnetic field of intensity $B$. 
The case $B=0$  was studied in \cite{jara2015superdiffusion}, where it was proved that the energy super-diffuses as a $3/2$-fractional diffusion, while if $B\neq0$ it is described by a $5/3$-fractional diffusion, see \cite{saito20195}. 
In \cite{cane2022superdiffusion} it is quantified the intensity of the magnetic field to pass from one regime to the other and the description of the transition mechanism is given in terms of a L\'evy process that interpolates between the two fractional universality classes: $3/2$-L\'evy and $5/3$-L\'evy. 
As we can see in the previous examples, by either changing the dynamics or by changing the potential, the limiting laws can have a variety of forms and they are universal in the sense that they can be obtained from a general collection of microscopic dynamics and they do not depend on the special features of the underlying microscopic dynamics but on their phenomenology.
The crossover that has been established for those universal laws can be obtained by tuning certain parameters at the microscopic dynamics, which permit a comparison of one dynamics with respect to another. 
When one dominates, the system falls in one university class, and then it can cross to another universality class when the other dynamics dominate instead. 
With the previous examples, we have just seen a possible way to cross from the Edwards-Wilkinson (EW) \cite{edwards1982surface} class, i.e. Gaussian fluctuations, to the KPZ universality class, by considering the fluctuation field for $V'(\eta)$ in the case of the Toda lattice potential and by changing the parameter regulating the Hamiltonian dynamics. 
When both dynamics are tuned with the same strength we see the SB equation connecting both universality classes. We highlight that the crossover from the SB equation to the (conjectured) KPZ fixed point is still missing. 
Moreover, we have seen how to cross from the EW class to the $3/2$-L\'evy class and these classes are connected by a process that is obtained as the superposition of the two dynamics. 
There are other mechanisms to cross between EW and $3/2$-L\'evy class, see for example \cite{bernardin2015normal, bernardin2018nonlinear}. 
In \cite{bernardin2018nonlinear} a perturbation of the quadratic potential with a small anharmonicity of the quartic form (a version of the Fermi-Pasta-Ulam (FPU) chain) is considered and it is realized that for a certain strength of the quartic potential the limit fluctuations continue to behave as $3/2$-L\'evy, i.e. the same behavior as the purely quadratic potential. 
In Figure~\ref{fig:universality_classes}, we can see a cartoon in $(1+1)$-dimensions representing some universality classes that have already been obtained from microscopic dynamics. 
According to \cite{popkov2015fibonacci, popkov2015universality} for multi-component systems with several conservation laws, there can be a family of universality classes given by $\alpha$-L\'evy processes with an exponent $\alpha$ written as the quotient of consecutive Fibonacci numbers. 
However, we have less understanding for the other universality classes, where the existence of them is predicted from NLFH, and there are still open problems. 
As future works, we have strong interests on characterizing the limit behavior of fluctuations in the possible universality classes in  a mathematically rigorous way. 

\subsection*{Organization of the Paper}
In Section~\ref{sec:result}, we give a precise definition of our model and state our main results.
Our first result is about a pair of fluctuations corresponding to the two conserved quantities: volume and energy. 
After that, we consider two possible choices of the linear combination of these original fields, from which the SB equation and $3/2$-L\'{e}vy fractional diffusion equation are derived. 
Section~\ref{sec:trivial} is devoted to the computation of Dynkin's martingale decomposition for the original fluctuation fields and the proof of the first result.  
Then, in Section~\ref{sec:cancellation}, we present a way, based on the expression of Dynkin's martingale that allows choosing the linear combination of the original fields, in order to cancel some diverging terms, provided the asymmetric part of the generator is larger than the value considered in the first result for the original fluctuation fields.  
The proof for the convergence to the SB equation and to the $3/2$-L\'{e}vy fractional diffusion will be given in Section~\ref{sec:KPZ} and Section~\ref{sec:3/2Levy}, respectively. 
In the appendix we discuss the predictions from NLFH.

\subsection*{Notation}
Given two real-valued functions $f$ and $g$ depending on a variable $u\in\mathbb{R}^d$ we will write $f(u)\lesssim g(u)$ if there exists a constant $C>0$ such that $f(u)\le C g(u)$ for any $u$. 
Moreover, we write $f=O(g)$ (resp. $f=o(g)$) in the neighborhood of $u_0$ if $|f|\le |g|$ in the neighborhood of $u_0$ (resp. $\lim_{u\to u_0}f(u)/g(u)=0$). 
Sometimes it will be convenient to make precise the dependence of the constant $C$ on some extra parameters and this will be done by the standard notation $C(\lambda)$ if $\lambda$ is the extra parameter. 
Finally, we denote by $\langle \cdot, \cdot \rangle_{L^2(\mathbb{R})}$  the inner product in $L^2(\mathbb{R})$, i.e. for any $f,g\in L^2(\mathbb{R})$ 
\begin{equation*}
\langle f, g\rangle_{L^2(\mathbb{R})}
\coloneqq \int_{\mathbb{R}} f(x) g(x) dx, 
\end{equation*}
 and by $\| \cdot\|_{L^2(\mathbb{R})}$  the $L^2(\mathbb{R})$-norm, i.e. 
$
\|f \|_{L^2(\mathbb{R})}
\coloneqq ( \langle f, f\rangle_{L^2(\mathbb{R})} )^{1/2}$.

\section{Statement of Results}
\label{sec:result}
\subsection{Model Description}
\label{subsec:model_description}
Here we define our model and state our  main results. 
Let $\mathscr{X}=\mathbb{R}^{\mathbb{Z}}$ be a space of oscillators whose elements are denoted by $\eta = (\eta_j)_{j \in \mathbb{Z}}$.  
For each $\eta\in\mathscr{X}$ and $\alpha>0$, define 
\begin{equation*}
||| \eta |||_\alpha 
= \sum_{j\in\mathbb{Z}} | \eta_j| e^{-\alpha|j|} ,
\end{equation*}
and let $\Omega_\alpha$ be the set of configuration $\eta\in\mathcal{X}$ such that $|||\eta|||_\alpha<+\infty$. 
The normed space $(\Omega_\alpha,|||\cdot |||_\alpha)$ turns out to be a Banach space. 
Set $\Omega=\cap_{\alpha>0} \Omega_\alpha$. 
Note that the space $\Omega$ is a complete metric space with respect to the distance 
\begin{equation*}
d(\eta_1,\eta_2)
= \sum_{\ell\in\mathbb{N}} 2^{-\ell} \min \{ 1, ||| \eta_1 -\eta_2 |||_{1/\ell}\} . 
\end{equation*}
In this paper, we consider the dynamics of oscillators on the configuration space $\Omega$ driven by a general potential $V$, which satisfies the following assumptions.

\begin{assumption}
\label{ass:oscillator_potential}
Assume that a smooth convex \Add{non-negative} function $V:\mathbb{R}\to\mathbb{R}$ satisfies $V(0)=V'(0)=0$ and $V^{\prime\prime}(0)=1$.
Moreover, assume for each $k\in\{0,\ldots,5\}$ that, the derivative $V^{(k)}(\eta)$ has at most exponential growth, that is, there exists a constant $\gamma_V>0$ such that 
\begin{equation*}
\max_{0\le k\le 5} \sup_{\eta\in\mathbb{R}} \big|  e^{-\gamma_V |\eta|} V^{(k)}(\eta) \big| <+\infty.
\end{equation*}
Here we used the convention $V^{(0)}=V$. 
\end{assumption}

In what follows, we fix a nonlinear function $V$ satisfying Assumption \ref{ass:oscillator_potential}\Erase{, which makes the function $V$ to be non-negative,} and for each $\beta>0$, we set $V_{\beta}(\eta)=\beta^{-2}V(\beta\eta)$.  
Throughout this paper, we are interested in the case when $\beta$ is small, which means that the temperature of the system goes to infinity. 
In this regime, by a Taylor expansion \eqref{eq:taylor_expansion_introduction}, the scaled potential $V_{\beta}$ converges to the harmonic potential $V(\eta)=\eta^2/2$, when $\beta\to0$. 
Now, we define the dynamics.  
Let $n>0$ be a scaling parameter and let $L = S + \alpha_n A$ with $\alpha_n=O(n^{-\kappa})$, where the operators $A$ and $S$, act on a smooth local function $f:\Omega\to\mathbb R$, as
\begin{equation*}
A f(\eta) = \sum_{j \in \mathbb{Z}}
\big( V^\prime_{\beta_n}(\eta_{j-1}) - V^\prime_{\beta_n}(\eta_{j+1}) \big)
\partial_{\eta_j} f(\eta)
= \sum_{j \in \mathbb{Z}} (\xi_{j-1} - \xi_{j+1})\partial_{\eta_j} f(\eta),
\end{equation*}
and 
\begin{equation*}
S f(\eta) =\frac{1}{2} \sum_{j \in \mathbb{Z}}
\big( f(\eta^{j,j+1}) - f(\eta) \big),
\end{equation*}
respectively. 
Above we introduced the notation 
\begin{equation}
\label{eq:def_xi}
\xi_j=V_\beta^\prime(\eta_j).
\end{equation}
We also introduced the notation  $\eta^{j,j+1}$ for the configuration obtained from $\eta$ after we exchange the occupation variables at sites $j$ and $j+1$: $(\eta^{j,j+1})_k=\eta_{j+1}$ when $k=j$, $(\eta^{j,j+1})_k=\eta_{j}$ when $k=j+1$ and otherwise $(\eta^{j,j+1})_k=\eta_{k}$.
Moreover, let $\theta(n) = n^a$ be a \Add{time} scaling factor for some $a>0$ and set $L_n = \theta(n)L$. 
To be concerned in the high temperature regime, we will take the inverse temperature in such a way that \Add{$$\beta=\beta_n\to0$$} as $n\to\infty$. 
\Erase{Throughout the paper, we consider the Markov process which is generated by the operator $L_n$.}
In \cite{bernardin2012anomalous} it was assumed that the potential $V:\mathbb R\to \mathbb R$ is a non-negative smooth function, such that  $Z_{\beta_n, b,\lambda}$ is well defined for $b>0$ and $\lambda \in\mathbb R$ and moreover, the potential satisfies the following condition: $0\le V^{\prime\prime}\le C$ for some $C>0$.
We observe that these conditions are sufficient to have well-defined dynamics, but they are not necessary. 
For example, it is shown in \cite{bernardin2014anomalous} that the Toda lattice potential $V(\eta)=e^{-\eta}-1+\eta$ defines a well-defined dynamics on $\Omega$ with some restrictions on parameters of stationary states, which we will describe below. 
In the whole article, we assume that our choice of potential is such that the dynamics in infinite volume is well-defined.
Let $\{ \eta_j(t):t\ge 0, j\in\mathbb{Z}\}$ be \Erase{a}\Add{the} Markov process on $\Omega$ generated by $L_n$ where we omit the dependency of $n$. 

\subsection{Invariant Measure and Static Estimates}
\label{sec:static_estimate}
Similarly to the interacting diffusion case \cite{diehl2017kardar, jara2020stationary, hayashi2022derivation}, the following product Gibbs measure associated to the potential $V_\beta$, whose common marginal is given by 
\begin{equation}\label{eq:inv_mea}
\nu_{\beta, \mathsf b,\lambda}(d\eta_j)
= \frac{1}{Z_{\beta, \mathsf b,\lambda}} \exp (-\mathsf bV_\beta(\eta_j) + \lambda \eta_j) d\eta_j,
\end{equation}
is invariant for this process.
Above, $\mathsf b>0$, $\lambda\in\mathbb{R}$ are constants  and \begin{equation}\label{def:partition_function}
Z_{\beta, \mathsf b,\lambda}\coloneqq\int_{-\infty}^{+\infty}\exp(-\mathsf bV_\beta (\eta)+\lambda \eta) d\eta,
\end{equation} 
is the normalizing constant, which is called the partition function, that makes of $\nu_{\beta, \mathsf b,\lambda}$ a probability measure, and $\beta>0$ is the inverse temperature whose value is taken in such a way that $Z_{\beta,\mathsf b,\lambda}$ is finite and thus the measure $\nu_{\beta,\mathsf b,\lambda}$ is well-defined. 
The next lemma assures that when we take the inverse temperature $\beta$ sufficiently small, an assumption that we will impose throughout the paper, the partition function $Z_{\beta,\mathsf b,\lambda}$ is finite and a uniform exponential moment bound holds. 
In what follows, we assume that the inverse temperature $\beta$ satisfies $\beta<1$ for simplicity.

\begin{lemma}
\label{lem:static_estimate}
Fix $b>0$ and $\lambda\in\mathbb{R}$. 
For any $\gamma>0$, there exists $C_\gamma>0$ and $\beta_c=\beta_c(\gamma)>0$ such that 
\begin{equation}
\label{eq:uniform_moment_bound}
\sup_{\beta<\beta_c} E_{\nu_{\beta,\mathsf b,\lambda}} \big[e^{\gamma |\eta_j|} \big] < C_\gamma,
\end{equation}
\Add{where $E_{\mu}$ denotes the expectation with respect to a probability measure $\mu$.}
\end{lemma}
\begin{proof}
First we show that the partition function has a uniform bound. 
Let we choose $K_0>0$ in such a way that $K_0\ge 4\lambda/b$.  Moreover, we assume that $\beta_n>0$ is sufficiently small, so that 
\begin{equation*}
\beta_n M_* K_0^2 \le 
 K_0 ,
\end{equation*}
where we set $M_*= \sup_{\delta\in(0,1)} 
\big| V^{(3)}(\delta K_0) \big|$.  
Recall that the Taylor expansion for $V_\beta$ yields
\begin{equation*}
 V'_{\beta}(K_0)
= K_0 + (\beta/2)
V^{(3)}(\delta \beta K_0) K_0^2  
\end{equation*}
for some $\delta=\delta(K_0)\in(0,1)$. 
Recalling that we assumed $\beta<1$, \Erase{this has the bound}\Add{we get the bound} 
\begin{equation*}
V^\prime_\beta(K_0)
\ge K_0 - (\beta/2) M_*  K_0^2
\ge K_0/2. 
\end{equation*}
Similarly, we have that 
\begin{equation*}
- V'_\beta(-K_0) 
\ge -K_0/2. 
\end{equation*}
On the other hand, recalling the convexity and the non-negativity of $V$, we have that  
\begin{equation*}
V_\beta(K_1) \ge 
V_\beta(K_2) + V_\beta^\prime (K_2)(K_1-K_2)
\ge V_\beta^\prime(K_2) (K_1-K_2), 
\end{equation*}
for any $K_1,K_2\in\mathbb{R}$. 
As a consequence, for each $\eta\in\mathbb{R}$, taking $K_1=\eta$ and $ K_2= \mathrm{sgn}(\eta) K_0$ in the above bound, we have 
\begin{equation*}
V_\beta(\eta)
\ge
2\mathsf b^{-1}|\lambda|(|\eta|-K_0).
\end{equation*}
Hence, the partition function has the bound   
\begin{equation*}
Z_{\beta,\mathsf b,\lambda}
= \int_{\mathbb{R}} e^{-\mathsf bV_\beta(\eta) + \lambda\eta}
d\eta 
\lesssim \int_{\mathbb{R}} e^{-2|\lambda \eta| + \lambda \eta} d\eta. 
\end{equation*}
In particular, the density function  $\varphi_{\beta,\mathsf b,\lambda}(\eta)\coloneqq \exp( -\mathsf bV_\beta(\eta)+\lambda \eta )$ inside the integral of $Z_{\beta,\mathsf b,\lambda}$ is bounded by an $L^1(\mathbb{R})$-function, which is independent of \Erase{$\beta_n$}\Add{$\beta$}. 
Therefore, noting that $\varphi_{\beta,\mathsf b,\lambda}(\eta)\to \exp (-\mathsf b|\eta|^2/2 + \lambda\eta)$ as $\beta\to0$ for each $\eta\in\mathbb{R}$, we have 
\begin{equation*}
\lim_{\beta\to0} 
\Erase{Z_{\beta_n,b,\lambda}}
\Add{Z_{\beta,\mathsf b,\lambda}}
= \int_{\mathbb{R}} \lim_{\beta\to0}\varphi_{\beta,\mathsf b,\lambda}(\eta)d\eta 
= e^{\lambda^2/2}, 
\end{equation*}
by the dominated convergence theorem. 

Now, we estimate the exponential moment
\begin{equation*}
E_{\nu_{\beta,\mathsf b,\lambda}}\big[ e^{\gamma |\eta_j| }\big]
= \frac{1}{Z_{\beta,\mathsf b,\lambda}}
\int_{\mathbb{R}}
e^{ - \mathsf bV_\beta(\eta_j) + \gamma|\eta_j|+ \lambda\eta_j} d\eta_j .
\end{equation*}
By repeating the same argument for the bound of the partition function, we see that the density $\exp \big(-\mathsf bV_\beta(\eta_j)+ \gamma|\eta_j|+\lambda\eta_j \big) $ is bounded by an $L^1(\mathbb{R})$ function which is independent of $\beta$. 
Hence, combining with the boundedness of the partition function, we complete the proof. 
\end{proof}

\begin{example}
\label{ex:fpu}
Here we consider the case where the nonlinear function $V$ is given as the so-called FPU-$\alpha$ potential which has the form $V(\eta)=\eta^2/2+ \alpha \eta^3 + \eta^4/4$ with $\alpha\in \mathbb{R}$. 
This potential is convex if, and only if, $3\alpha^2 \le 1$ and it also satisfies all the conditions in the Assumption \ref{ass:oscillator_potential}. 
Moreover, the partition function $Z_{\beta,\mathsf b,\lambda}$ is finite if $\beta^{-1}\ge 24\mathsf b^{-1}|\lambda\alpha|$ and the uniform moment bound \eqref{eq:uniform_moment_bound} holds with $\beta_c=\mathsf b(24(|\lambda|+\gamma)|\alpha|)^{-1}$. 
\end{example}

\begin{example}
\label{ex:kvm}
As another example, the Toda lattice potential $V(\eta)=e^{-\eta}-1+\eta$ satisfies all the conditions in the Assumption \ref{ass:oscillator_potential} and its invariant measure is a log-gamma distribution, which is well-defined when $\beta^{-1}> \mathsf b^{-1}|\lambda|$.
Moreover, we can see that $E_{\nu_{\beta,\mathsf b,\lambda}}[e^{\gamma |\eta|}]$ is bounded by a constant which is independent of $\beta$, provided $\beta^{-1}>\mathsf b^{-1}(\gamma+|\lambda|)$.
In other words, the assertion \eqref{eq:uniform_moment_bound} holds with $\beta_c=\mathsf b(\gamma+|\lambda|)^{-1}$. 
\end{example}

In the sequel, we consider fixed values of $\mathsf b,\lambda$ and take sufficiently small $\beta$, which forces the measure $\nu_{\beta,\mathsf b,\lambda}$ to be well-defined and the uniform exponential moment bound \eqref{eq:uniform_moment_bound} holds with $\gamma=2\gamma_V$. 
The above uniform exponential moment bound will be used to bound error terms of the Taylor expansion, combined with the  exponential growth of derivatives of the potential. 
For instance, by Taylor's theorem, we have 
\begin{equation}
\label{eq:taylor_expansion_upto_beta}
\xi_j = \eta_j + \frac{1}{2}V^{(3)}(0) \beta \eta_j^2
+ \varepsilon_j , \quad 
\varepsilon_j=\frac{1}{3!}\beta^2 V^{(4)}(\delta \beta\eta_j) \eta_j^3,
\end{equation}
for some $\delta=\delta(\eta_j)\in(0,1)$. 
Recall that the derivative $V^{(4)}(\cdot)$ has at most exponential growth, which yields the bound $\varepsilon_j \lesssim \beta^2 e^{2\gamma_V |\eta_j|}$ where $\gamma_V$ is the constant in Assumption \ref{ass:oscillator_potential}.  
Thus, we have the bound 
\begin{equation*}
E_{\nu_{\beta,\mathsf b,\lambda}}[\varepsilon_j] 
\le \beta^2 E_{\nu_{\beta,\mathsf b,\lambda}}\big[ e^{2\gamma_V|\eta_j|} \big]
\lesssim \beta^2, 
\end{equation*}
where we used \eqref{eq:uniform_moment_bound} in the second inequality.  
Such an argument will be repeatedly used to show that error terms of the Taylor expansion are negligible, by taking sufficiently large $\beta$. 
(See for example the statement after~\eqref{eq:taylor_expansion_upto_beta_squared} below.)

\subsection{Main Results}
In the sequel, we fix $\lambda\in \mathbb{R}$ and set $\mathsf b=1$ for simplicity. 
Then we take $\beta=\beta_n$ depending on the scaling parameter $n$ and we simply write \Add{$$\nu_n = \nu_{\beta_n, \mathsf b, \lambda}$$} for this setting.
We fix a time horizon $T$.  
We consider our interacting oscillator model $\{ \eta_j(t);j \in \mathbb{Z} \}$ with generator $L_n$, starting from the invariant measure $\nu_n$.   
We denote by ${D} ([0,T],\Omega) $ 
the space of \Erase{the }c\`adl\`ag (right-continuous and with left limits) trajectories taking values in $\Omega$. 
Let $\mathbb P_n$ be the probability measure \Erase{in}\Add{on} $ {D} ([0,T],\Omega) $ which is induced by $\nu_n$ and let $\mathbb E_n$ denote the expectation with respect to  $\mathbb P_n$.
The above interacting oscillator model admits two conserved quantities\Erase{; volume and energy:}
\begin{equation*}
\sum_{j \in \mathbb{Z}} \eta_j  \quad \textrm{and}\quad 
\sum_{j \in \mathbb{Z}} \zeta_j,  
\end{equation*}
\Add{which are refered to as volume and energy, respectively,} where 
\begin{equation}
\label{eq:def_zeta}
\zeta_j = V_\beta (\eta_j).
\end{equation} 
Let $\mathcal S(\mathbb R)$ be the space of Schwartz functions and \Erase{$S'(\mathbb R)$}\Add{$\mathcal{S}'(\mathbb R)$} its dual, i.e. the set of linear continuous functionals defined on $\mathcal S(\mathbb R)$ and taking real values. Let \Erase{$ {D} ([0,T],S'(\mathbb R)) $}\Add{$D([0,T],\mathcal S'(\mathbb R))$} be
the space of \Erase{the }c\`adl\`ag (right-continuous and with left limits) trajectories in $\mathcal S'(\mathbb R)$.

Then, as natural objects to investigate, we introduce the volume and the energy fluctuation fields as elements of $ D([0,T],\mathcal{S}^\prime(\mathbb{R}))$ that are defined for each  $\varphi \in \mathcal{S}(\mathbb{R})$ in the following manner:
\begin{equation}
\label{eq:volume_fluctuation_definition}
\mathcal{V}^n_t (\varphi)
= \frac{1}{\sqrt{n}} \sum_{j \in \mathbb{Z}}
\overline{\eta}_j(t) T^-_{f_1t} \varphi^n_j,
\end{equation}
and
\begin{equation}
\label{eq:energy_fluctuation_definition}
\mathcal{E}^n_t (\varphi)
= \frac{1}{\sqrt{n}} \sum_{j \in \mathbb{Z}}
\overline{\zeta}_j(t) T^-_{f_2t} \varphi^n_j. 
\end{equation}
Here we used the bar notation over a random variable to mean the deviation from its expectation with respect to the invariant measure $\nu_n$: $\overline{\eta}_j=\eta_j - E_{\nu_n}[\eta_j]$, for instance.
In addition, $T^-_\cdot$ denotes a shift operator $T^-_{v}\varphi^n_j = \varphi^n_{j-v}$ for each $v\in \mathbb{R}$, and  $\varphi^n_j = \varphi(j/n)$.
In the above definition, $f_1=f_1(n)$ and $f_2=f_2(n)$ are constants depending on $n$, which may take different values for volume and energy.  

In what follows, the velocities $f_1$ and $f_2$ are carefully calibrated, depending on the time scale $n^a$, the weak asymmetry $\alpha_n$ and the inverse temperature $\beta_n$, which also depends on $n$. 
Our first result is concerned with the limiting behavior of the pair of fluctuation fields $\mathcal{Z}^n=(\mathcal{V}^n,\mathcal{E}^n)$.

\begin{theorem}[Linear fluctuations]
\label{thm:pair_fluctuations}
Let $\mathcal{V}^n$ and $\mathcal{E}^n$ be the volume and the energy fluctuation fields which are defined by \eqref{eq:volume_fluctuation_definition} and \eqref{eq:energy_fluctuation_definition}, respectively.
We take $\theta(n)=n^2$, $\alpha_n=O(n^{-\kappa})$, $\beta_n=O(n^{-\delta})$ and $(f_1,f_2)=(2\theta(n)\alpha_n,0)$ and assume $\kappa > 1/2$ and $\kappa + \delta > 1$.
Moreover, assume $\lambda=0$. 
Then, the pair of fluctuation fields $(\mathcal{V}^n, \mathcal{E}^n)$ converges in distribution in $D([0,T]),\mathcal{S}^\prime(\mathbb{R})^2)$ to some $(u^1,u^2)$, which satisfies the following system of uncorrelated stochastic heat equations with additive noise:
\begin{equation*}
\begin{aligned}
 \partial_t u^i
= \frac{1}{2}\partial_x^2 u^i 
+ \sigma^i \partial_x \dot{W}^i.
\end{aligned}
\end{equation*}
Here $\dot{W}^1=\dot{W}^1(t,x)$ and $\dot{W}^2=\dot{W}^2(t,x)$ denote  independent space-time white-noises and we set $\sigma^1=1$ and $\sigma^2=1/\sqrt{2}$.
\end{theorem}

\Add{
\begin{remark}
When $\lambda\neq0$, in the time evolution of the pair of fluctuation fields a \Erase{blowing}\Add{diverging} term survives, see the terms of degree one  in the second line of \eqref{eq:Z_anti_symmetric}. 
Currently, we are not aware of how to treat these terms and therefore, we imposed the condition $\lambda=0$.
\end{remark}
}

\begin{remark}
The analysis in the sub-diffusive time scale, i.e., $\theta(n)=n^a$ with $a<2$, and derivation of some linear drift term, as explored in \cite[Remark 2]{ahmed2022microscopic}, are possible extensions.
However, this is not necessary for our goal to derive anomalous behavior, so that we will not stick to these cases here. 
\end{remark}

\Add{On the contrary, when $\kappa \le 1/2$, we expect to derive a non-trivial SPDE.} 
Indeed, in this sub-critical regime, we can choose distinct moving frames $f_1$ and $f_2$, in order that nonlinear terms survive in the limit. 
These nonlinear terms, however, cannot be written in terms of the fluctuation fields, since we might have different \Add{frame} velocities for volume and energy. 
To avoid such a difficulty, instead, we consider  linear combinations of volume and energy fluctuations with a common velocity $v_n$ as follows: 
\begin{equation}
\label{eq:xi_fluctuation_def}
\mathcal{X}^n_t(\mathfrak{u}_n;\varphi)
= \frac{1}{\sqrt{n}} \sum_{j\in\mathbb Z}
\big( 
\overline{\eta}_j(t) + \mathfrak{u}_n \overline{\zeta}_j(t) \big)
T^-_{v_nt} \varphi^n_j , 
\end{equation}
We will choose the constant $\mathfrak{u}_n$ and the velocity $v_n$ properly, in such a way that we can characterize the limiting behavior of the new field $\mathcal{X}^n$. 
We can show that this is established for the following two cases:
\begin{equation}
\label{eq:two_regimes}
\begin{aligned}
&\text{(i) } 
\mathfrak{u}_n=\mathfrak{u}_n^1\coloneqq c_3\beta_n, \quad 
v_n=v^1_n \coloneqq \theta(n)\alpha_n (2+2\lambda c_3 \beta_n). 
\qquad \qquad \qquad \qquad \qquad \qquad  \\ 
& \text{(ii) } 
\mathfrak{u}_n=\mathfrak{u}_n^2\coloneqq-1/\lambda~(
\lambda\neq0) , \quad 
v_n= v^2_n\coloneqq 0.
\end{aligned}
\end{equation}
\Add{
For simplicity, we defined 
\begin{equation}
\label{eq:derivative_coefficient}
c_k=V^{(k)}(0).     
\end{equation}
}

\begin{remark}
Note that we can consider the fluctuation fields associated with $\zeta_j + \mathfrak{u}_n\eta_j$, instead of \eqref{eq:xi_fluctuation_def}. 
Then, the role of $\mathfrak{u}_n$ becomes that of $\mathfrak{u}_n^{-1}$, so  that the condition $\lambda\neq0$ in the case (ii) can be removed. 
For this reason, we will consider \Erase{a scaled fluctuation fields}\Add{the fluctuation fields multiplied by $\lambda$} in the case (ii). 
\end{remark}

Recall the Taylor expansion \eqref{eq:taylor_expansion_upto_beta}. 
Then, we notice that the fluctuation field of the conserved quantity $\eta_j + c_3\beta_n \zeta_j$ in the case (i) roughly matches that of \Erase{$\xi_j$}\Add{$\xi_j=V'_\beta(\eta_j)$}. 
On the other hand, for the case~(ii), recall that $E_{\nu_n}[\xi_j]=\lambda$ and that we took along the article $b=1$ for simplicity. 
Then, we have 
\begin{equation*}
\begin{aligned}
(\overline{\xi}_j)^2
= \lambda^2 -2\lambda (\eta_j - \lambda^{-1} \zeta_j) + \varepsilon_j,
\end{aligned}
\end{equation*}
where $\varepsilon_j$ is an error term which satisfies the bound $\varepsilon_j \lesssim \beta_n e^{2\gamma_V | \eta_j|}$. 
Hence the fluctuations of $\overline{\eta}_j -\lambda^{-1}\overline{\zeta}_j$ asymptotically match those of $(\overline{\xi}_j)^2$ when $\beta_n$ is small. 
Our main theorems are concerned with the derivation of the SBE and the $3/2$-L\'{e}vy anomalous diffusion for the  fluctuations of $\xi_j$ and $(\overline{\xi}_j)^2$, respectively.
Let us recall the notion of the stationary energy solution of the stochastic Burgers equation which was introduced in~\cite{gonccalves2014nonlinear}. 
Let $\nu,D > 0$ and $\Lambda \in \mathbb{R}$ be fixed constants and consider the $(1+1)$-dimensional stochastic Burgers equation
\begin{equation}
\label{eq:SBE_general}
\partial_t u = \nu \partial_x^2 u 
+ \Lambda \partial_x u^2 + \sqrt{D} \partial_x \dot{W}.
\end{equation}
We begin with the definition of stationarity. 

\begin{definition}
We say that an $\mathcal{S}^\prime (\mathbb{R})$-valued process $u = \{ u_t  : t \in [0,T] \} $ satisfies condition \textbf{(S)} if for all $t \in [0,T]$, the random variable $u_t$ has the same distribution as the space white-noise with variance $D/(2\nu)$. 
\end{definition}

For a process $u = \{ u_t: t \in [0,T]\}$ satisfying the condition \textbf{(S)}, we define 
\begin{equation}
\label{eq:def_quadratic_function_approximation}
\mathcal{A}^\varepsilon_{ s, t } (\varphi ) = \int_s^t \int_{\mathbb{R} } u_r (\iota_\varepsilon (x; \cdot) )^2 \partial_x \varphi (x ) dx dr .  
\end{equation}
for every $0 \le s < t \le T $, $\varphi \in \mathcal{S} (\mathbb{R} ) $ and $\varepsilon > 0 $. 
Here we defined the function $\iota_\varepsilon (x ; \cdot ) : \mathbb{R} \to \mathbb{R}  $ by $\iota_{ \varepsilon } (x ; y) =  \varepsilon^{ - 1 } \mathbf{1}_{ [ x , x  + \varepsilon ) } (y) $ for each $x \in \mathbb{R} $ and $\epsilon>0$. 
Although the function $\iota_\varepsilon(x,\cdot)$ does not belong to the Schwartz space, the quantity \eqref{eq:def_quadratic_function_approximation} is well-defined when $u$ satisfies the condition \textbf{(S)}.

\begin{definition}
Let $u = \{ u_t :t \in [0,T]\}$ be a process satisfying the condition \textbf{(S)}. 
We say that the process $u$ satisfies the energy estimate \textbf{(EC)} if there exists a constant $\kappa > 0$ such that for any $\varphi \in \mathcal{S} (\mathbb{R} )$, any $0 \le s < t \le T$ and any $0 < \delta < \varepsilon < 1 $,  
\begin{equation*}
\Erase{\mathbb{E}_n}\Add{\mathbb E} \big[ \big| \mathcal{A}^\varepsilon_{ s, t } (\varphi ) - \mathcal{A}^\delta_{ s, t } (\varphi ) \big|^2 \big] 
\le \kappa \varepsilon (t-s) \| \partial_x \varphi \|^2_{ L^2(\mathbb{R})} .
\end{equation*}
Here $\mathbb{E}$ denotes the expectation with respect to the measure of a probability space where the process $u$ lives. 
\end{definition}

Then the following result is proved in \cite{gonccalves2014nonlinear}. 

\begin{proposition}
\label{prop:nonlinear}
Assume $\{ u_t:t\in [0,T]\} $ satisfies the conditions \textbf{(S)} and \textbf{(EC)}. Then there exists an $\mathcal{S}^\prime (\mathbb{R} )$-valued process $\{ \mathcal{A}_t : t \in [0, T ] \} $ with continuous trajectories such that  
\begin{equation*}
\mathcal{A}_t (\varphi ) = \lim_{ \varepsilon \to 0 } \mathcal{A}^\varepsilon_{ 0, t } (\varphi) ,
\end{equation*}
in $L^2 $ for every $t \in [0,T]$ and $\varphi \in \mathcal{S}(\mathbb{R})$.  
\end{proposition}

From the last proposition, thinking that the singular term $\partial_x u^2 $ is given by the last quantity, we can define a solution of \eqref{eq:SBE_general} as follows.

\begin{definition}
\label{def:energysol}
We say that an $\mathcal{S}^\prime(\mathbb{R})$-valued process $u=\{u (t, \cdot) : t\in [0,T] \}$ is a stationary energy solution of the stochastic Burgers equation \eqref{eq:SBE_general} if 
\begin{enumerate}
\item The process $u$ satisfies the conditions \textbf{(S)} and \textbf{(EC)}. 
\item For all $\varphi \in \mathcal{S} (\mathbb{R} )$, the process 
\begin{equation*}
u_t(\varphi) - u_0 (\varphi) - \nu \int_0^t u_s (\partial_x^2 \varphi ) ds 
+ \Lambda \mathcal{A}_t (\varphi) ,
\end{equation*}
is a martingale with quadratic variation $D \| \partial_x \varphi \|^2_{ L^2 (\mathbb{R} ) } t $ where $\mathcal{A}_\cdot$ is the process obtained in Proposition \ref{prop:nonlinear}. 
\item For all $\varphi \in \mathcal{S} (\mathbb{R} )$, writing $\hat{u}_t = u_{T-t}$ and $\hat{ \mathcal{A} }_t = - (\mathcal{A}_T - \mathcal{A}_{ T- t })$, the process
\begin{equation*}
\hat{u}_t (\varphi) - \hat{u}_0 (\varphi ) - \nu \int_0^t \hat{u}_s (\partial_x^2 \varphi) ds 
+ \Lambda \hat{\mathcal{A}}_t (\varphi) ,
\end{equation*}
is a martingale with quadratic variation $D \| \partial_x \varphi \|^2_{L^2(\mathbb{R})}t$. 
\end{enumerate}
\end{definition}

Then there exists a unique-in-law stationary energy solution of \eqref{eq:SBE_general}. 
Existence was shown in \cite{gonccalves2014nonlinear} and then uniqueness was proved in \cite{gubinelli2018energy}.

\begin{theorem}[The SBE regime]
\label{thm:SBE_derivation}
Let $\mathfrak{u}_n=\mathfrak{u}_n^1=c_3 \beta_n$ and $v_n =v_n^1= \theta(n)\alpha_n (2+2\lambda c_3\beta_n)$.
We consider the diffusive scaling $\theta(n)=n^2$, and assume $\lim_{n\to\infty}\sqrt{n}\alpha_n\beta_n =1$ and $\lim_{n\to\infty} n\beta_n^4=0$. 
Moreover, we assume $c_4-c_3^2=0$.
Let $\mathcal{X}^n$ be the fluctuation field defined by \eqref{eq:xi_fluctuation_def}.
Then, $\mathcal{X}^n$ converges in distribution in $D([0,T],\mathcal{S}^\prime(\mathbb{R}))$ to the stationary energy solution of the stochastic Burgers equation
\begin{equation}
\label{eq:sbe_main_theorem}
\begin{aligned}
 \partial_t u = \frac{1}{2} \partial_x^2 u 
- c_3 \partial_x u^2
+ \partial_x \dot{W} ,
\end{aligned}
\end{equation}
where $\dot{W}=\dot{W}(t,x)$ denotes the one-dimensional space-time white-noise.
\end{theorem}

\begin{remark}
If $c_3=0$, then above $\mathfrak u_n^1=0$ and  the quantity we are looking at is just the volume, but in that case the limit is given by an Ornstein-Uhlenbeck process. We also observe that, independently of the value of $c_3$, if  $\lim_{n\to\infty}\sqrt{n}\alpha_n\beta_n =0$, then the limit is again an  Ornstein-Uhlenbeck process.
\end{remark}

\begin{remark}
Note that in the Taylor expansion of $\xi_j$\Add{, which is defined in \eqref{eq:def_xi},} the quantity $c_4-c_3^2$ is the coefficient of the cubic term\Add{, which is proportional to $\eta_j^3$}, see \eqref{eq:taylor_expansion_upto_beta_squared} below. 
If the condition $c_4-c_3^2=0$ does not hold, then it is possible that another term appears in the limit. 
However, it is not clear how to write this cubic term with respect to  volume and energy fluctuation fields, therefore, we decided to avoid the issue by choosing that coefficient to be zero, but it is certainly an interesting question that we leave for a future work. 
\end{remark}

Our second result is concerned with the case (ii) of \eqref{eq:two_regimes}.
For $j\in\mathbb{Z}$, we define the correlation function 
\begin{equation}
\label{eq:correlation}
S_j(t) = 
\frac{1}{2} \mathbb{E}_n \big[
\big(\overline{\zeta}_j(0) -\lambda \overline{\eta}_j(0) \big)
\big(\overline{\zeta}_j(t) -\lambda\overline{\eta}_j(t) \big) \big]. 
\end{equation}

\Add{We denote by $k_* \in \{ 3,4,\ldots \}\cup\{\infty\}$ the smallest integer such that $V^{(k_*)}(0)\neq 0$. 
As examples of other cases of the potential, let us suppose that the function $V$ is the FPU-$\alpha$ potential which is given in Example \ref{ex:fpu}.
Then, $k_*=3$ when $\alpha\neq0$, whereas $k_*=4$ when $\alpha=0$.  
Note that when $\kappa^*=\infty$ the nonlinear function $V(\eta)$ is proportional to the purely harmonic potential $\eta^2/2$ and moreover for any $\beta\in\mathbb R$ we have that $\beta^{-2}V(\beta \cdot)=V(\cdot)$.  The next result has been shown in \cite{bernardin2018nonlinear} when the nonlinear potential $V$ is the purely harmonic one given above and we  extend it to generic potentials.}

\begin{theorem}[The 3/2-L\'{e}vy regime]
\label{thm:3/2L_derivation}
Let $\theta(n)=n^a$ and $\alpha_n=\gamma n^{-\kappa}$ with $\gamma>0$, $\kappa\ge 0$ and $a=\min\{ 3/2+3\kappa/2, 2\}$.  
 Assume $\beta_n = O(n^{-b})$ with $b\ge 1/(2k_*-4)$ when $k_*\ge 4$ whereas $b\ge 1/4$ when $k_*=3$.  
Let $\mathcal{X}^n$ be \Add{the} fluctuation field defined by
\begin{equation*}
\mathcal{X}^n_t(\varphi) 
= \frac{1}{\sqrt{n}}\sum_{j\in\mathbb{Z}}
(\overline{\zeta}_j - \lambda \overline{\eta}_j) \varphi(j/n), 
\end{equation*}
for each $\varphi\in \mathcal{S}(\mathbb{R})$. 
Let $f,g$ be smooth functions on $\mathbb{R}$ with compact support. 
Then, 
\Add{ 
\begin{equation*}
\lim_{n\to\infty}
\frac{1}{2} \mathbb E_n \big[\mathcal{X}^n_0(g) \mathcal{X}^n_t(f) \big]
= \iint_{\mathbb{R}^2}
f(x)g(y) P^{\gamma,\kappa}_t(x-y) dx dy,
\end{equation*}
}
and $\{P^{\gamma,\kappa}_t(x):t\ge 0, x\in\mathbb{R}\}$ is the fundamental solution of the equation 
$\partial_t u= \mathbb{L}_{\gamma,\kappa}u,$
where 
\begin{equation}
\label{eq:levy_operator}
\mathbb{L}_{\gamma,\kappa}
= \frac{1}{2}\mathbf{1}_{\kappa\ge1/3}\Delta
- \gamma^{3/2} \mathbf{1}_{\kappa\le1/3} \mathscr{L},
\end{equation}
with $\mathscr{L}=-\frac{1}{\sqrt{2}}[(-\Delta)^{3/4}-\nabla(-\Delta)^{1/4}]$. 
\end{theorem}

\begin{remark}
We observe that the result of \Add{the} last theorem gives (as long as the value of $b$ satisfies the assumptions of the theorem) that the second mode has exactly the same behavior as the energy in the case of an harmonic potential $V(\eta)=\eta^2/2$, i.e. the same diagram as in Figure \ref{fig:energy}.
\end{remark}

\begin{remark}
In~\cite{bernardin2018nonlinear}, a perturbation given by a quartic function, namely, $V(\eta)=\eta^2/2+ \overline{\gamma} \eta^4$, is studied and the authors proved that the same equation driven by the operator \eqref{eq:levy_operator} is derived provided $\overline{\gamma}$ decays faster than $n^{-1/4}$ for large $n$. 
Observing that $\overline{\gamma}=\beta_n^2$, this bound is better than the one we obtained in Theorem~\ref{thm:3/2L_derivation}, whereas their proof is based on the exact form of the potential. 
Additionally, we expect that Theorem~\ref{thm:3/2L_derivation} is valid for any order of perturbation which is given by the quartic function. 
\end{remark}

\Add{
\begin{remark}
The interested reader can find an heuristic argument on the value of the critical exponent $\kappa=1/3$ that appears  in the statement of Theorem \ref{thm:3/2L_derivation} in the introduction of \cite{bernardin2018weakly}. 
The critical scale $\alpha_n=n^{-1/3}$ is obtained by solving a macroscopic differential equation which  well approximates the time evolution of the correlation function. 
\end{remark} 
}

\section{Proof of Theorem \ref{thm:pair_fluctuations}}
\label{sec:trivial}
\subsection{The martingale Decomposition}
In the sequel, we define discrete derivative operators as follows. 
\begin{equation}
\label{eq:definition_discrete_derivatives}
\begin{aligned}
&\nabla^{1,n} \varphi^n_j = 
n(\varphi^n_{j+1} - \varphi^n_{j}), \quad  
\nabla^{2,n} \varphi^n_j =
\frac{n}{2}(\varphi^n_{j+1} - \varphi^n_{j-1}),   \\
&\Delta^n \varphi^n_j = n^2 (\varphi^n_{j+1} + \varphi^n_{j-1} - 2 \varphi^n_j) . 
\end{aligned}
\end{equation}
Here, note that both $\nabla^{n,1}$ and $\nabla^{2,n}$ approach to the continuous derivative $\partial$, though, the rates of convergence are different. 
Indeed, we can show that $\nabla^{1,n}\varphi_j-\partial_x\varphi_j=O(n^{-1})$ while $\nabla^{2,n}\varphi_j- \partial_x\varphi_j = O(n^{-2})$ with the help of the mean-value Theorem.
In this section, we give a proof of Theorem \ref{thm:pair_fluctuations}, assuming $\kappa>1/2$, namely, the strength of asymmetry is sufficiently weak.

Our starting point is a martingale decomposition for the pair of fluctuation fields $(\mathcal{V}^n,\mathcal{E}^n)$. 
Hereafter we set $\mathcal{Z}^n_t = (\mathcal{V}^n_t, \mathcal{E}^n_t)$ and recall \eqref{eq:volume_fluctuation_definition} and \eqref{eq:energy_fluctuation_definition}.
To compute the correlation between two martingales associated with $\mathcal{V}^n$ and $\mathcal{E}^n$, similarly to \cite{ahmed2022microscopic}, we apply Dynkin's formula, see, for example, Lemma A.1.5.1 of \cite{kipnis1998scaling}, to $\mathcal{Z}^n_t (\overrightarrow{\varphi}) =  \mathcal{V}^n_t(\varphi_1) + \mathcal{E}^n_t(\varphi_2)$ for each $\overrightarrow{\varphi}=(\varphi_1,\varphi_2)\in \mathcal{S}(\mathbb{R})^2$. 
Then, we have that 
\begin{equation*}
\mathcal{N}^n_t (\overrightarrow{\varphi})
= \mathcal{Z}^n_t(\overrightarrow{\varphi}) 
- \mathcal{Z}^n_0 (\overrightarrow{\varphi})
- \int_0^t (\partial_s + L_n)\mathcal{Z}^n_s (\overrightarrow{\varphi})ds ,
\end{equation*}
and $\mathcal{N}^n_t(\varphi)^2 -\langle \mathcal{N}^n(\varphi)\rangle_t$ where 
\begin{equation}
\label{eq:qv_vector}
\begin{aligned}
\langle \mathcal{N}^n (\overrightarrow{\varphi})\rangle_t
&= \int_0^t \big( L_n \mathcal{Z}^n_s (\overrightarrow{\varphi})^2
- 2 \mathcal{Z}^n_s(\overrightarrow{\varphi}) L_n 
\mathcal{Z}^n_s (\overrightarrow{\varphi}) \big) ds \\
&= \frac{\theta(n)}{2n^3} \int_0^t \sum_{j \in \mathbb{Z}}
(\eta_{j}(s)- \eta_{j+1}(s))^2 
\big( \nabla^{1,n} T^-_{f_1t} 
\varphi_1(j/n) \big)^2 ds\\
&\quad+ \frac{\theta(n)}{2n^3} \int_0^t \sum_{j \in \mathbb{Z}}
(\zeta_{j}(s) - \zeta_{j+1}(s))^2 
\big( \nabla^{2,n} T^-_{f_2t} 
\varphi_2(j/n) \big)^2 ds\\
&\quad+ \frac{\theta(n)}{n^3} \int_0^t \sum_{j \in \mathbb{Z}}
(\eta_{j}(s) - \eta_{j+1}(s)) (\zeta_{j}(s)- \zeta_{j+1}(s)) 
\nabla^{1,n} T^-_{f_1t} 
\varphi_1(j/n)  
\nabla^{1,n} T^-_{f_2t} 
\varphi_2(j/n) ds 
\end{aligned}
\end{equation}
are martingales with respect to the natural filtration of the process.
In what follows, we give a generic computation for any $\lambda$ and at the final step we set $\lambda=0$ to show Theorem \ref{thm:pair_fluctuations}. 
Since the limiting measure is product and homogeneous,  
\begin{equation*}
\begin{aligned}
\lim_{n\to \infty} E_{\nu_n}[(\eta_j-\eta_{j+1})(\zeta_j - \zeta_{j+1})]
&= E_{\eta_j \sim \mathcal{N}(\lambda,1)}
[(\eta_j-\eta_{j+1})(\eta_j^2/2-\eta_{j+1}^2/2)]\\
&= E_{\eta_j \sim \mathcal{N}(\lambda,1)}
[\eta_j^3]-E_{\eta_j \sim \mathcal{N}(\lambda,1)}
[\eta_{j+1}\eta_j^2]\\
&=\lambda^3+3\lambda-\lambda(\lambda^2+1)=2\lambda .
\end{aligned}
\end{equation*}
Here the limiting procedure is justified as we mentioned in Section~\ref{sec:static_estimate}. 
Similarly, we have that 
\begin{equation*}
\begin{aligned}
\lim_{n\to \infty} E_{\nu_n}[(\eta_j-\eta_{j+1})^2]
&=2 \mathrm{Var}_{\eta_j \sim \mathcal{N}(\lambda,1)}[\eta_j]
=2, 
\end{aligned}
\end{equation*}
and 
\begin{equation*}
\begin{aligned}
\lim_{n\to \infty} E_{\nu_n}[(\zeta_j-\zeta_{j+1})^2]
&=\frac{1}{2} E_{\eta_j \sim \mathcal{N}(\lambda,1)}[\eta_j^4 - \eta_j^2 \eta_{j+1}^2] \\
&= \frac{1}{2}[(\lambda^4+6\lambda^2+3) - (\lambda^2+1)^2]
= 2\lambda^2+1 ,
\end{aligned}
\end{equation*}
where we used the fact that $E[X^{2n}]=(2n-1)!!$ when $X$ is drawn from the standard normal distribution. 
In particular, when $\lambda=0$ the expectation of the third term in the utmost right-hand side of \eqref{eq:qv_vector} vanishes as $n\to +\infty$ and the noise terms are diagonalized.
In summary, under the diffusive scaling $\theta(n)=n^2$, we have that 
\begin{equation}
\label{eq:qv_limit}
\lim_{n\to \infty}\mathbb{E}_n [\langle \mathcal{N}^n(\overrightarrow{\varphi}) \rangle_t]
=  t \| \partial_x \varphi_1\|^2_{L^2(\mathbb{R})}
+ \frac{2\lambda^2+1}{2} t \|\partial_x \varphi_2\|^2_{L^2(\mathbb{R})} 
+ 2\lambda t \langle \partial_x \varphi_1, 
\partial_x \varphi_2 \rangle_{L^2(\mathbb{R})}.
\end{equation}

Now, we compute the action of the generator $L_n$ on $\mathcal{Z}^n$. 
First, for the symmetric part, we get
\begin{equation*}
\begin{aligned}
& \int_0^t \theta(n)S \mathcal{V}^n_s (\varphi)ds 
= \frac{\theta(n)}{2n^{5/2}} \int_0^t \sum_{j \in \mathbb{Z}}
\overline \eta_j(s) \Delta^n T^-_{f_1s}\varphi^n_j ds , \\
& \int_0^t \theta(n)S \mathcal{E}^n_s (\varphi)ds 
= \frac{\theta(n)}{2n^{5/2}} \int_0^t \sum_{j \in \mathbb{Z}}
\overline{\zeta}_j(s) \Delta^n T^-_{f_2s}\varphi^n_j ds ,
\end{aligned}
\end{equation*}
for each $\varphi\in \mathcal{S}(\mathbb{R})$ where we used the short-hand notation $\varphi^n_j=\varphi(j/n)$.  
In particular, the symmetric part $S$  should always be accelerated by the diffusive scaling $\theta(n)=n^2$ in order to obtain a non-trivial limit.  

Next, we consider the anti-symmetric part of the generator. 
The action of the anti-symmetric part, after some rearrangement, is calculated as follows.

\begin{lemma}
\label{lem:action_antisymmetric_base}
We have that 
\begin{equation*}
\begin{aligned}
\int_0^t (\partial_s + \theta(n) \alpha_n A)
\mathcal{V}^n_s(\varphi) ds 
= \frac{\theta(n)\alpha_n}{n^{3/2}} 
\int_0^t \sum_{j \in \mathbb{Z}}
\bigg(2 \overline{\xi}_j(s) - \frac{f_1(n)}{\theta(n)\alpha_n} \overline{\eta}_j(s) \bigg) 
\partial_x T^-_{f_1s} \varphi^n_jds 
+ E^{1,n}_t ,
\end{aligned}
\end{equation*}
and 
\begin{equation*}
\begin{aligned}
&\int_0^t (\partial_s + \theta(n) \alpha_n A)
\mathcal{E}^n_s(\varphi) ds \\
&\quad = \frac{\theta(n)\alpha_n}{n^{3/2}} 
\int_0^t \sum_{j \in \mathbb{Z}}
\bigg(\overline{\xi}_{j+1}(s) \overline{\xi}_j(s)
+ 2 \lambda \overline{\xi}_j 
- \frac{f_2(n)}{\theta(n)\alpha_n} \overline{\zeta}_j(s) \bigg) 
\partial_x T^-_{f_2s} \varphi^n_jds 
+ E^{2,n}_t,
\end{aligned}
\end{equation*}
where,  
\begin{equation}
\label{eq:action_generator_vector_error}
\mathbb{E}_n \bigg[\sup_{0\le t\le T} \big| E^{1,n}_t 
\big|^2 \bigg] 
\lesssim T^2 \frac{\theta(n)^2 \alpha_n^2}{n^6} \quad \quad\textrm{and}\quad  \quad
\mathbb{E}_n \bigg[\sup_{0\le t\le T} \big| E^{2,n}_t 
\big|^2 \bigg] 
\lesssim T^2 \frac{\theta(n)^2 \alpha_n^2}{n^4} .
\end{equation}
\end{lemma}
\begin{proof}
We begin with the computation for $\mathcal{V}^n$. 
For each $\varphi\in \mathcal{S}(\mathbb{R})$, note that 
\begin{equation*}
\begin{aligned}
& \int_0^t (\partial_s + \theta(n) \alpha_n A)
\mathcal{V}^n_s(\varphi)ds \\
&\quad= \frac{\theta(n)\alpha_n}{n^{3/2}} \int_0^t \sum_{j \in \mathbb{Z}} 
2\overline{\xi}_j(s) \nabla^{2,n} T^-_{f_1s}\varphi^n_j ds 
- \frac{1}{\sqrt{n}} \int_0^t \sum_{j \in \mathbb{Z}} 
\frac{f_1(n)}{n} \overline{\eta}_j(s) \partial_x T^-_{f_1s}\varphi^n_j ds \\
&\quad= \frac{\theta(n)\alpha_n}{n^{3/2}} \int_0^t 
\sum_{j \in \mathbb{Z}}
\bigg(2\overline{\xi}_j(s) - \frac{f_1(n)}{\theta(n)\alpha_n} \overline{\eta}_j(s) \bigg) 
\partial_x T^-_{f_1s} \varphi^n_j ds 
+ E^{1,n}_t,
\end{aligned}
\end{equation*}
where $E^{n,1}_t$ satisfies the desired estimate \eqref{eq:action_generator_vector_error}. 
In the last line we replaced the discrete derivative by the continuous one with the variance of order $O(\theta(n)^2 \alpha_n^2 n^{-6})$, which is estimated as follows. 
Recall from \eqref{eq:definition_discrete_derivatives} the definition of the discrete derivative. 
Then, \Add{recalling the definition of $\nabla^{2,n}$ in\eqref{eq:definition_discrete_derivatives}, due to Taylor's theorem,} \Erase{$|\nabla^{2,n}\varphi^n_j - \partial_x \varphi^n_j |\lesssim  n^{-2} \partial_x^2\varphi^n_j$}\Add{$|\nabla^{2,n}\varphi^n_j - \partial_x \varphi^n_j |\lesssim  n^{-3} \partial_x^3\varphi^n_j$}  and thus by the Schwarz's inequality  
\begin{equation*}
\begin{aligned}
&\mathbb{E}_n \bigg[\sup_{0\le t\le T} \bigg| \frac{\theta(n)\alpha_n}{n^{3/2}} 
\int_0^t \sum_{j\in\mathbb{Z}} \overline{\xi}_j(s)
T^-_{f_1s} \big(\nabla^{2,n} \varphi^n_j - \partial_x \varphi^n_j \big) ds 
\bigg|^2 \bigg] \\
&\quad \lesssim \frac{\theta(n)^2\alpha_n^2T}{n^3} \int_0^T \sum_{j \in \mathbb{Z}}
\mathbb{E}_n \big[ \overline{\xi}_j(s)^2 \big]  
T^-_{f_1s} \big(\nabla^{2,n} \varphi^n_j - \partial_x \varphi^n_j \big)^2 ds 
\lesssim T^2 \frac{\theta(n)^2\alpha_n^2}{n^{6}}.
\end{aligned}
\end{equation*}
Here we used the fact that \Erase{$E_n [(\overline{\xi}_j)^2]$}\Add{$E_{\nu_n} [(\overline{\xi}_j)^2]$} is bounded by a constant which is independent of $n$ by Lemma \ref{lem:static_estimate}.
On the other hand, for the energy fluctuation field we have that 
\begin{equation*}
\begin{aligned}
& \int_0^t (\partial_s + \theta(n)\alpha_n A)\mathcal{E}^n_s(\varphi)ds \\
& \quad= \frac{\theta(n)\alpha_n}{n^{3/2}} \int_0^t \sum_{j \in \mathbb{Z}} 
\xi_j(s) \xi_{j+1}(s) \nabla^{1,n} T^-_{f_2s} \varphi^n_j ds 
- \frac{1}{\sqrt{n}} \int_0^t \sum_{j \in \mathbb{Z}}
\frac{f_2(n)}{n}\overline{\zeta}_j(s)\partial_x T^-_{f_2s} \varphi^n_j ds\\ 
& \quad= \frac{\theta(n)\alpha_n}{n^{3/2}} \int_0^t \sum_{j \in \mathbb{Z}} 
\bigg( \overline{\xi}_{j}(s) \overline{\xi}_{j+1}(s) 
+ \lambda (\overline{\xi}_j(s) + \overline{\xi}_{j+1}(s))
\bigg) \nabla^{1,n} T^-_{f_2s} \varphi^n_j ds \\ 
&\qquad- \frac{1}{\sqrt{n}} \int_0^t \sum_{j \in \mathbb{Z}}
\frac{f_2(n)}{n}\overline{\zeta}_j(s)\partial_x T^-_{f_2s} \varphi^n_j ds.
\end{aligned}
\end{equation*}
Similarly to the volume fluctuation, we replace the discrete derivative by the continuous one, and since variables are all centered, the error term of this replacement satisfies the bound in the assertion.
Note that the cost of this replacement is worse than the one above since now $|\nabla^{1,n}\varphi^n_j - \partial_x \varphi^n_j|\lesssim n^{-1}$. 
Finally, by a summation-by-parts the last display equals to 
\begin{equation*}
\begin{aligned}
& \quad \frac{\theta(n)\alpha_n}{n^{3/2}} \int_0^t \sum_{j \in \mathbb{Z}} 
\bigg( \overline{\xi}_{j}(s) \overline{\xi}_{j+1}(s) 
+ 2\lambda \overline{\xi}_j(s) 
- \frac{f_2(n)}{\theta(n)\alpha_n} \overline{\zeta}_j(s) 
\bigg) \partial_x T^-_{f_2s} \varphi^n_j ds 
+ E^{2,n}_t ,
\end{aligned}
\end{equation*}
where $E^{2,n}_t$ satisfies the bound \eqref{eq:action_generator_vector_error}. 
This is justified as follows: 
\begin{equation*}
\begin{aligned}
&\mathbb{E}_n \bigg[\sup_{0\le t \le T} \bigg| 
\int_0^t  \frac{\theta(n)\alpha_n}{n^{3/2}}\sum_{j\in \mathbb{Z}}
(\overline{\xi}_{j+1}(s)- \overline{\xi}_{j}(s))
\partial_x T^-_{f_2s} \varphi^n_j ds 
\bigg|^2 \bigg]\\
&\quad = \mathbb{E}_n \bigg[\sup_{0\le t \le T} \bigg|
\int_0^t  \frac{\theta(n)\alpha_n}{n^{5/2}}
\sum_{j\in \mathbb{Z}}
\overline{\xi}_{j}(s)
T^-_{f_2s} \nabla^{1,n} \partial_x \varphi^n_j ds 
\bigg|^2 \bigg]\\
&\quad \le \frac{T\theta(n)^2\alpha_n^2}{n^5} 
\int_0^T \mathbb{E}_n \bigg[ \bigg( \sum_{j \in\mathbb{Z}} \overline{\xi}_j(s) \nabla^{1,n} \partial_x T^-_{f_2s} \varphi^n_j \bigg)^2 \bigg] ds  
\lesssim \frac{T^2\theta(n)^2\alpha_n^2}{n^4} 
\| \partial_x^2 \varphi\|^2_{L^2(\mathbb{R})}.
\end{aligned}
\end{equation*}
\end{proof}

With these representations for the action of the generator at hand, we have that
\begin{equation}
\label{eq:Z_anti_symmetric}
\begin{aligned}
& \int_0^t (\partial_s + \theta(n) \alpha_n A)
\mathcal{Z}^n_s (\overrightarrow{\varphi}) ds \\
&\quad= 
\frac{\theta(n)\alpha_n}{n^{3/2}} \int_0^t 
\sum_{j \in \mathbb{Z}}
\bigg(2\overline{\xi}_j(s) - \frac{f_1(n)}{\theta(n)\alpha_n} \overline{\eta}_j(s) \bigg) 
\partial_x T^-_{f_1s} 
\varphi^{n,1}_j ds \\
&\qquad + \frac{\theta(n)\alpha_n}{n^{3/2}} \int_0^t \sum_{j \in \mathbb{Z}} 
\bigg( \overline{\xi}_{j}(s) \overline{\xi}_{j+1}(s) 
+ 2\lambda \overline{\xi}_j (s) 
- \frac{f_2(n)}{\theta(n)\alpha_n} \overline{\zeta}_j(s) 
\bigg) 
\partial_x T^-_{f_2s} 
\varphi^{n,1}_j ds 
+ E^n_t ,
\end{aligned}
\end{equation}
where $E^n_t$ satisfies the bound 
\begin{equation*}
\mathbb{E}_n\bigg[ \sup_{0\le t\le T} \big| E^n_t \big|^2 \bigg]
\lesssim T^2 \frac{\theta(n)^2\alpha_n^2}{n^4}, 
\end{equation*}
which vanishes as $n\to\infty$ provided $\alpha_n=o_n(1)$ and $\theta(n)=n^2$.
Now, recall the assumptions in Theorem \ref{thm:pair_fluctuations} and the Taylor expansion \eqref{eq:taylor_expansion_upto_beta}.  
We note that the first term on the utmost right-hand side of \eqref{eq:Z_anti_symmetric} is estimated for any $\varphi$ as 
\begin{equation*}
\mathbb{E}_n \bigg[\sup_{0\le t \le T} \bigg|
\frac{2\theta(n)\alpha_n}{n^{3/2}} \int_0^t \sum_{j\in\mathbb{Z}} 
(\overline{\xi}_j -\overline{\eta}_j) (s) 
\partial_x T^-_{f_1s} \varphi^n_j ds \bigg|^2 \bigg]
\lesssim T^2 \frac{\beta_n^2\theta(n)^2\alpha_n^2}{n^2} \| \partial_x \varphi \|^2_{L^2(\mathbb{R})} ,
\end{equation*}
which vanishes under the assumptions of Theorem \ref{thm:pair_fluctuations}: $\theta(n)=n^2$ and $n\beta_n =o_n(1)$. 
On the other hand, the quantity  
\begin{equation*}
\int_0^t \sum_{j\in\mathbb{Z}}
\overline{\xi}_j(s) \overline{\xi}_{j+1}(s) 
\partial_x T^-_{f_2s} 
\varphi^{n,2}_j ds ,
\end{equation*}
remains non-trivial but not divergent as $n\to\infty$.
This will be verified in the SBE regime with the help of the second-order Boltzmann-Gibbs principle, see Proposition \ref{prop:2BG}.
In particular, the above quadratic term which is multiplied by $\theta(n)\alpha_nn^{-3/2}$ vanishes in $L^2$ when $\theta(n)=n^2$ and $\kappa>1/2$. 
Therefore, combined with the additional assumption $\lambda=0$ and $f_2=0$, the quadratic term in \eqref{eq:Z_anti_symmetric} vanishes in $L^2$ as $n\to\infty$.

Hence, the anti-symmetric part \eqref{eq:Z_anti_symmetric} is negligible in the limit $n\to\infty$, and we have a martingale decomposition 
\begin{equation*}
\mathcal{Z}^n_t (\overrightarrow{\varphi}) 
= \mathcal{Z}^n_0 (\overrightarrow{\varphi}) 
+ \int_0^t \mathcal{Z}^n_s (\partial_x^2 \overrightarrow{\varphi}) ds
+ \mathcal{N}^n_t (\overrightarrow{\varphi}) 
+ \Add{\tilde E^n_t}, 
\end{equation*}
\Add{
where $\tilde E^n_t$ is an error term which is negligible in the sense that 
\begin{equation*}
\lim_{n\to\infty}
\mathbb E_n \Big[\sup_{0\le t\le T} \big|\tilde E^n_t\big|^2 \Big]
= 0. 
\end{equation*}
}
From the above decomposition, we can deduce the assertions of Theorem \ref{thm:pair_fluctuations}. 
We do not present more steps on this because it is very similar to the approach of \cite{ahmed2022microscopic} and we refer the interested readers to that article for details.

\section{Choice of Fluctuation Fields}
\label{sec:cancellation}
In this section, we find the linear combination of the volume and energy fluctuations, from which we derive the stochastic Burgers equation and the $3/2$-L\'{e}vy anomalous diffusion.    
Throughout this section, we consider $\theta(n)=n^2$. 
We apply Dynkin's formula, see, for example, \cite[Lemma A.1.5.1]{kipnis1998scaling}, 
for $\mathcal{X}^n(\cdot) = \mathcal{X}^n(\mathfrak{u}_n;\cdot)$ defined by \eqref{eq:xi_fluctuation_def} with a common velocity $v_n$. 
Then, for each $\varphi\in \mathcal{S}(\mathbb{R})$, 
\begin{equation*}
\mathcal{M}^n_t(\varphi)
= \mathcal{X}^n_t (\varphi)
- \mathcal{X}^n_0 (\varphi)
-\int_0^t (\partial_s + L_n) \mathcal{X}^n_s(\varphi) ds, 
\end{equation*}
and $\mathcal{M}^n_t(\varphi)^2 -\langle \mathcal{M}^n (\varphi) \rangle_t$ where 
\begin{equation*}
\langle \mathcal{M}^n(\varphi)\rangle_t
= \int_0^t \big( L_n \mathcal{X}^n_s(\varphi)^2 - 2 \mathcal{X}^n_s(\varphi) L_n \mathcal{X}^n_s(\varphi) \big) ds,
\end{equation*}
are martingales with respect to the natural filtration of the process. 
From simple, but long, computations, 
\begin{equation}
\label{eq:qv_scalar}
\begin{aligned}
\langle \mathcal{M}^n (\varphi) \rangle_t 
&= \frac{\theta(n)}{2n^3} \int_0^t \sum_{j \in \mathbb{Z}} 
(\eta_j(s) -\eta_{j+1}(s))^2 \big( \nabla^{n,1} T^-_{v_ns} \varphi^n_j \big)^2 ds \\
&\quad+ \frac{(\mathfrak{u}_n)^2 \theta(n)}{2n^3} \int_0^t \sum_{j \in \mathbb{Z}} 
(\zeta_j(s) -\zeta_{j+1}(s))^2 \big( \nabla^{n,1} T^-_{v_ns} \varphi^n_j \big)^2 ds \\
&\quad+ \frac{\mathfrak{u}_n\theta(n)}{n^3}
\int_0^t \sum_{j\in\mathbb{Z}} 
(\eta_j(s) - \eta_{j+1}(s)) (\zeta_j(s) -\zeta_{j+1}(s))
(\nabla^{n,1} T^-_{v_ns}\varphi^n_j)^2 ds. 
\end{aligned}
\end{equation}
We proceed by computing the action of the generator on $\mathcal{X}^n$. 
First, for the symmetric part, we have that 
\begin{equation*}
\int_0^t \theta(n) S \mathcal{X}^n_s(\varphi) ds 
= \frac{\theta(n)}{2n^{5/2}}
\int_0^t \sum_{j\in\mathbb{Z}} 
(\overline{\eta}_j + \mathfrak{u}_n 
\overline{\zeta}_j ) 
\Delta^n T^-_{v_ns} \varphi^n_j ds ,
\end{equation*}
which is replaced by $\frac{\theta(n)}{2n^2}\int_0^t \mathcal{X}^n_s(\partial_x^2 \varphi) ds$ whose error is of order $O(\theta(n)^2 n^{-6})$, i.e.:
\begin{equation*}
\begin{aligned}
\mathbb{E}_n \bigg[ \sup_{0\le t \le T} \bigg|
\int_0^t \theta(n)S\mathcal{X}^n_s(\varphi)ds  -\frac{\theta(n)}{2n^2} \int_0^t \mathcal{X}^n_s(\partial_x^2 \varphi) ds \bigg|^2 \bigg]
\lesssim 
T^2 \frac{\theta(n)^2}{n^6} .
\end{aligned}
\end{equation*}

On the other hand, the action of the anti-symmetric part can be computed as follows. 

\begin{lemma}
We have that 
\begin{equation}
\label{eq:action_antisymmetric_X}
\begin{aligned}
&\int_0^t (\partial_s + \theta(n) \alpha_n A) \mathcal{X}^n_s(\varphi) ds\\
&\quad= \frac{\mathfrak{u}_n \theta(n) \alpha_n}{n^{3/2}}
\int_0^t \sum_{j \in \mathbb{Z}} 
\overline{\xi}_j(s) \overline{\xi}_{j+1}(s) 
\partial_x T^-_{v_ns} \varphi^n_j ds \\
&\qquad+ \bigg( (2+2\lambda\mathfrak{u}_n) \frac{\theta(n)\alpha_n}{n^{3/2}} 
- \frac{v_n}{n^{3/2}} \bigg) 
\int_0^t \sum_{j \in \mathbb{Z}} 
\overline{\eta}_j 
\partial_x T^-_{v_ns} \varphi^n_j ds \\
&\qquad+ \bigg( (2+2\lambda\mathfrak{u}_n) \frac{\theta(n)\alpha_n c_3\beta_n}{n^{3/2}} - \frac{\mathfrak{u}_nv_n}{n^{3/2}} \bigg) 
\int_0^t \sum_{j \in \mathbb{Z}} 
\overline{\zeta}_j 
\partial_x T^-_{v_ns} \varphi^n_j ds \\
&\qquad+ (2+2\lambda\mathfrak{u}_n) \frac{\theta(n)\alpha_n (c_4-c_3^2)\beta_n^2}{3n^{3/2}}
\int_0^t \sum_{j \in \mathbb{Z}} \zeta_j^n(s) \eta^n_j(s)
\partial_x T^-_{v_ns} \varphi^n_j ds 
+ E^n_t ,
\end{aligned}
\end{equation}
where $E^n_t$ satisfies the bound 
\begin{equation*}
\mathbb{E}_n \bigg[\sup_{0\le t\le T} \big| E^n_t 
\big|^2 \bigg] 
\lesssim T^2 
\bigg( (2+2\lambda\mathfrak{u}_n)^2 \frac{\theta(n)^2\alpha_n^2 \beta_n^6}{n^2} 
\vee \frac{\theta(n)^2\alpha_n^2}{n^6} 
\vee \frac{(\mathfrak{u}_n)^2\theta(n)^2\alpha_n^2}{n^4} \bigg) . 
\end{equation*}
\end{lemma}
\begin{proof}
\Add{
Recalling that from Lemma \ref{lem:action_antisymmetric_base}, we have 
\begin{equation*}
\begin{aligned}
&\int_0^t \theta(n)\alpha_n A \mathcal{X}^n_s(\varphi) ds\\
&\quad= \frac{\theta(n)\alpha_n}{n^{3/2}}
\int_0^t \sum_{j \in \mathbb{T}_n}
\big( (2+ 2\lambda \mathfrak{u}_n)\overline{\xi}_j(s) 
+ \mathfrak{u}_n \overline{\xi}_j(s)\overline{\xi}_{j+1}(s) \big)
\partial_x T^-_{v_ns} \varphi^n_j ds
+ E^n_t ,
\end{aligned}
\end{equation*}
where the error term $E^n_\cdot$, which, by an abuse of notation, is denoted by the same notation as in  the assertion, satisfies the bound
\begin{equation*}
\mathbb{E}_n\bigg[ \sup_{0\le t\le T} \big| E^n_t 
\big|^2 \bigg] 
\lesssim T^2 \bigg( \frac{\theta(n)^2\alpha_n^2}{n^6} 
\vee \frac{(\mathfrak{u}_n)^2\theta(n)^2\alpha_n^2}{n^4} \bigg) .
\end{equation*}
}
Hence, \Erase{the anti-symmetric part of the action of the generator is represented as}
\begin{equation}
\label{eq:anti-symmetric_action_in_lemma}
\begin{aligned}
&\int_0^t (\partial_s + \theta(n) \alpha_n A) \mathcal{X}^n_s(\varphi) ds\\
&\quad= \frac{\theta(n)\alpha_n}{n^{3/2}}
\int_0^t \sum_{j \in \mathbb{Z}}
\big( (2+2\lambda \mathfrak{u}_n) \overline{\xi}_j(s) 
+ \mathfrak{u}_n \overline{\xi}_j(s)\overline{\xi}_{j+1}(s)  \big)
\partial_x T^-_{v_ns} \varphi^n_j ds \\
&\qquad - \frac{v_n}{n^{3/2}} 
\int_0^t \sum_{j \in \mathbb{Z}}
\big( \overline{\eta}_j(s) 
+ \mathfrak{u}_n \overline{\zeta}_j(s) \big) 
\partial_x T^-_{v_ns} \varphi^n_j ds + E^n_t .
\end{aligned}
\end{equation}
To proceed further, we make use of the following Taylor expansion up to order $\beta_n^3$. 
Recall that $c_k=V^{(k)}(0)$. 
By Taylor's theorem, we have that 
\begin{equation}\label{eq:taylor_expansion_upto_beta_squared}
\begin{aligned}
\xi_j - \eta_j 
= c_3 \beta_n \zeta_j 
+ \frac{c_4-c_3^2}{3} \beta_n^2 \zeta_j \eta_j 
+ \varepsilon_j,
\end{aligned}
\end{equation}
where $\varepsilon_j$ depends on the fifth derivative of the potential. 
As a consequence of Assumption \ref{ass:oscillator_potential}, we note that 
\begin{equation*}
\mathbb{E}_n\bigg[ \sup_{0\le t \le T} \bigg| 
\int_0^t \sum_{j\in\mathbb{Z}} \overline{\varepsilon}_j \varphi^n_j
ds \bigg|^2 \bigg] 
\lesssim T^2n\beta_n^6 \|\varphi\|^2_{L^2(\mathbb{R})},  
\end{equation*}
with the help of \eqref{eq:uniform_moment_bound} on the uniform moment bound. 
Now we substitute the expansion \eqref{eq:taylor_expansion_upto_beta_squared} in the first linear term on the right-hand side of \eqref{eq:anti-symmetric_action_in_lemma}. 
Then, we obtain the desired expression with an additional cost whose variance is proportional to $\beta_n^6$. 
\end{proof}

Note that the assumption $c_4-c_3^2=0$ forces the last term of \eqref{eq:action_antisymmetric_X} to be equal to zero. 
In addition, to cancel linear fluctuations, which are divergent in our regime, we choose $\mathfrak{u}_n$ and $v_n$ in such a way that 
\begin{equation}
\begin{cases}
\label{eq:linear_terms_cancellation}
\begin{aligned}
& \theta(n)\alpha_n(2+2\lambda\mathfrak{u}_n) - v_n = 0 ,\\
& \theta(n)\alpha_n c_3 \beta_n (2+2\lambda\mathfrak{u}_n) - \Erase{\mathfrak{u}}\Add{\mathfrak u_n}v_n = 0 .
\end{aligned}
\end{cases}
\end{equation}
Namely, the constant \Erase{$\mathfrak{u}$}\Add{$\mathfrak u_n$} should satisfy 
\begin{equation*}
\mathfrak{u}_n (2+2\lambda \mathfrak{u}_n)
= c_3 \beta_n (2+2\lambda \mathfrak{u}_n).
\end{equation*}

This quadratic equation has two solutions $\mathfrak{u}_n^1=c_3 \beta_n$ and $\mathfrak{u}_n^2 = - 1/\lambda$, which give $v^1_n = \theta(n)\alpha_n (2+2\lambda c_3 \beta_n)$ and $v^2_n= 0$, respectively. 
Therefore, when the relationship \eqref{eq:linear_terms_cancellation} is satisfied, the second and the third terms on the utmost right-hand side of \eqref{eq:action_antisymmetric_X} cancel. 
In the next section we are going to analyse the convergence of the fluctuation fields for each one of the previous quantities.

\section{Proof of Theorem \ref{thm:SBE_derivation}: The SB Equation}
\label{sec:KPZ}
\subsection{The Martingale Decomposition}
We consider the case $\mathfrak{u}_n=\mathfrak{u}_n^1= c_3 \beta_n$ and  $v_n=v_n^1=\theta(n)\alpha_n (2+2\lambda c_3 \beta_n)$ under the diffusive scaling $\theta(n)=n^2$. 
In this section, we are concerned with the following fluctuation field given in \eqref{eq:xi_fluctuation_def} with $\mathfrak{u}_n=\mathfrak{u}^1_n$ and $v_n=v^1_n$, i.e.,
\begin{equation*}
\mathcal{X}^n_t(c_3\beta_n; \varphi) = \frac{1}{\sqrt{n}} \sum_{j \in \mathbb{Z}}
(\overline{\eta}_j(t) + c_3 \beta_n \overline{\zeta}_j(t))
T^-_{v_n^1 t} \varphi^n_j.
\end{equation*}
We hereafter write $\mathcal{X}^n_t(c_3\beta_n;\varphi)=\mathcal{X}^n_t(\varphi)$ for simplicity. 
In this case, by \eqref{eq:action_antisymmetric_X}, the action of the anti-symmetric part of the generator is computed as follows.  
\begin{equation*}
\begin{aligned}
\int_0^t (\partial_s + \theta(n)\alpha_n A)\mathcal{X}^n_s(\varphi) ds 
&= \frac{\theta(n)\alpha_n c_3 \beta_n}{n^{3/2}} 
\int_0^t \sum_{j \in \mathbb{Z}}
\overline{\xi}_j(s) \overline{\xi}_{j+1}(s)
\nabla^{1,n} T^-_{v_n^1 s} \varphi^n_j ds 
+ E^n_t ,
\end{aligned}
\end{equation*}
where the remainder term $E^n_t$ satisfies the bound 
\begin{equation}
\label{eq:restriction}
\mathbb{E}_n \bigg[\sup_{0\le t\le T} \big| E^n_t\big|^2 \bigg]
\lesssim 
\frac{\theta(n)^2\alpha_n^2\beta_n^6}{n^2} 
\vee \frac{\theta(n)^2\alpha_n^2}{n^6}
\vee \frac{\beta_n^2\theta(n)^2\alpha_n^2}{n^4} .
\end{equation}
The action of the anti-symmetric part of the generator on the fluctuation field gives rise to the nonlinear term in the limiting equation, the SBE, provided $\lim_{n\to\infty}\sqrt{n}\alpha_n \beta_n = 1$, which will be justified by the second-order Boltzmann-Gibbs principle. 
Moreover, note that the error term vanishes as $n\to\infty$ when we additionally assume 
\begin{equation}
\label{eq:second_restriction}
n\beta_n^4=o_n(1).
\end{equation}

In summary, we have a martingale decomposition
\begin{equation}
\label{eq:KPZ_martingale_decomposition}
\mathcal{X}^n_t(\varphi)
= \mathcal{X}^n_0(\varphi)
+ \mathcal{S}^n_t(\varphi)
+ \mathcal{B}^n_t(\varphi)
+ \mathcal{M}^n_t(\varphi)
+ \mathcal{R}^n_t(\varphi), 
\end{equation}
where $\mathcal{M}^n_t$ is the Dynkin's martingale whose quadratic variation is given by \eqref{eq:qv_scalar} with $\mathfrak{u}_n=\mathfrak{u}^1_n=c_3 \beta_n$, $v_n=v^1_n$ and $\theta(n)=n^2$, and 
\begin{equation*}
\mathcal{S}^n_t (\varphi)
=\frac{1}{2\sqrt{n}} \int_0^t \sum_{j \in\mathbb{Z}} 
\big( \overline{\eta}_j+ c_3\beta_n \overline{\zeta}_j \big)(s) \Delta^n T^-_{v_n^1s} \varphi^n_j ds, 
\end{equation*}
\begin{equation*}
\mathcal{B}^n_t (\varphi)
= c_3 \int_0^t \sum_{j \in\mathbb{Z}} 
\overline{\xi}_j(s) \overline{\xi}_{j+1}(s) 
\nabla^{1,n} T^-_{v_n^1s} \varphi^n_j ds, 
\end{equation*}
and $\mathcal{R}^n_\cdot$ is a reminder term that vanishes in the following sense.
\begin{equation*}
\begin{aligned}
\lim_{n\to\infty} \mathbb{E}_n \bigg[ \sup_{0\le t\le T} \big| \mathcal{R}^n_t (\varphi) \big|^2 \bigg] =0 .
\end{aligned}
\end{equation*}
Following arguments below, it turns out that the reminder term $\mathcal{R}^n$ does not affect the limit. 
With the decomposition \eqref{eq:KPZ_martingale_decomposition} at hand, we give a proof of Theorem \ref{thm:SBE_derivation}. 
We will show that each term in the martingale decomposition is tight and then identify their limit points.

\begin{figure}[htb!]
\begin{center}
\begin{tikzpicture}[scale=0.17]
\draw (0,25) node[left]{$a$};
\draw (25,0) node[below]{$\kappa$};
\draw (8,0) node[below]{$ \frac 12-b$};
\draw (0,0) node[left]{$0$};
\draw (0,11) node[left]{$\frac 32+b$};
\draw (0,7) node[left] {$ \frac 43 (b+1)$};
\draw (0,20) node[left]{$2$};
\fill[light-gray] (0,11) -- (8,20) -- (25,20) -- (25,25) -- (0,25) -- cycle;
\fill[fill=apricot, fill opacity=0.5] (0,0) -- (25,0) -- (25,20) -- (8,20) -- (0,7)--cycle;
\fill[fill=gray, fill opacity=0.5] (0,7) -- (0,11) -- (8,20)--cycle;
\draw[-,=latex,gray, dashed, ultra thick] (0,7) -- (8,20) node[midway,below,sloped] {\textbf{\tiny{ ?}}};
\draw[-,=latex,blue!50!green,ultra thick] (8,20) -- (25, 20) node[midway,above,sloped] {\textbf{ \tiny{Diffusion} }};
\draw[-,=latex, magenta, dashed, ultra thick] (0,11) -- (8,20) node[midway,above,sloped] {\textbf{\tiny{KPZ-fp}}};
\node[circle,fill=red,inner sep=0.8mm] at (8,20) {};
\node[] at (8.5,20.5) [above] {\quad\textbf{\tiny{\textcolor{red}{SB}}}};
\node[] at (17.5,10) {\textbf{\tiny{no evolution}}};
\draw[-,=latex, dashed] (8,-0.1) -- (8,25);
\draw[->,>=latex] (0,0) -- (26,0);
\draw[->,>=latex] (0,0) -- (0,26);
\end{tikzpicture}
\end{center}
\caption{ $V'(\eta)$ fluctuations.}
\label{fig:xi_fluctuation}
\end{figure}
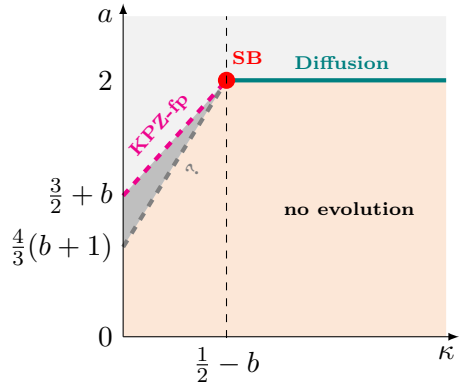

\begin{remark}
We observe that  from \cite[Theorem 4]{gonccalves2017second} (see in fact equation (5.30) in \cite{ahmed2022microscopic} for the precise bound) we know that the second moment of the second term in \eqref{eq:action_antisymmetric_X} is of order $\beta_n^2(\theta_n)^{3/2} \alpha_n^2n^{-2}$. 
Since  $\theta_n=n^a$, $\beta_n=n^{-b}$ and $\alpha_n=\alpha n^{-\kappa}$  we see that the quadratic term has no contribution to the limit if $a<\frac 43(\kappa+b+1)$. 
When  $a=\frac 43(\kappa+b+1)$ (see the line in gray color in the figure below) the same result should be true and in analogy to Figure \ref{fig:KvM} we believe that the same should be true up to the line $a=3/2+\kappa+b$, i.e. the line in magenta color in \Erase{the figure below}\Add{Figure \ref{fig:xi_fluctuation}}.
Moreover, if we require that the nonlinear term survives in the limit then we need to impose that the factor in front of the second term in \eqref{eq:action_antisymmetric_X} is of order $O(1)$. This means that the relationship $a-\Erase{k}\Add{\kappa}=3/2+b$ has to hold. Since we also need that the error terms in \eqref{eq:restriction} vanish as $n\to+\infty$ we need the following conditions to be true:
\begin{equation*}
a- \kappa 
< \min\{(1+3b) , 3 , (2 + b)\}. 
\end{equation*}
In fact, last conditions simultaneously hold when $1/4<b<3/2$.
Moreover, when $b<1/2$, the symmetric and martingale parts also survive.
Therefore,  for  $1/4<b<1/2$ we will derive the SB equation for $\kappa=1/2-b$ and $a=3/2+\kappa+b=2$. 
\end{remark}

\subsection{Tightness}
In this part, we recall some basic notions of the Skorohod space and the tightness of a sequence in the c\`{a}dl\`{a}g space for readers' convenience.
To begin with a general setting, let $E$ be a complete, separable metric space, endowed with a distance $d_E$.
Let $D([0,T],E)$ be the space of all right continuous functions with left limits taking values on $E$. 
Let $\lambda $ be the set of all strictly increasing continuous functions $\lambda$ from $[0,T]$ into itself. 
Then, we define 
\begin{equation*}
\| \lambda \| 
= \sup_{s\neq t} \bigg| \log \frac{\lambda(t)-\lambda(s)}{t-s} \bigg|
\end{equation*}
and define for each $X,Y\in D([0,T],E)$ 
\begin{equation*}
d(X,Y)
= \inf_{\lambda \in\Lambda} \max \big\{ \| \lambda \|, \sup_{0\le t\le T} d_E(X_t, Y_{\lambda(t)}) \big\}. 
\end{equation*}
Then it is known that the Skorohod space $D([0,T],E)$ endowed with the metric $d$ is a complete separable metric space, see \cite[Chapter 3]{billingsley1968convergence}. 
Next, in order to characterize the convergence of a sequence of paths in the Skorohod space, we make use of the following modified modulus of continuity: for each $X=\{X_t:t\in[0,T]\} \in D([0,T],E)$, set 
\begin{equation*}
w^\prime_X(\gamma)
= \inf_{\{t_i\}_{0\le t \le N}}
\max_{0\le i <N} \sup_{t_i\le s<t <t_{i+1}}
d_E(X_s, X_t),
\end{equation*}
where the first infimum is taken over all partitions $\{ t_i\}_{0\le i\le N}$ of the interval $[0,T]$ such that $0=t_0<t_1<\cdots <t_N=T$ and $t_i-t_{i-1}>\gamma$ for each $i=1,\ldots, N$. 
Then, the relative compactness of a sequence in the Skorohod space is characterized by the following Prohorov's theorem \cite[Theorem 4.1.3]{kipnis1998scaling}. 

\begin{proposition}
\label{prop:prohorov}
Let $\{ \mathbb{P}_n\}_n$ be a sequence of probability measures on $D([0,T],E)$. 
The sequence is relatively compact if, and only if,
\begin{itemize}
\item[(1)] For each $t\in[0,T]$ and each $\varepsilon>0$ there exists a compact set $K(t,\varepsilon)$ in $E$ such that $\mathbb{P}_n(X_t \notin K(t,\varepsilon))\le \varepsilon$. 
\item[(2)] For each $\varepsilon>0$, we have $
\lim_{\gamma\to0} \limsup_{n\to \infty}
\mathbb{P}_n (w^\prime_X(\gamma) > \varepsilon) = 0$. 
\end{itemize}
\end{proposition}

Here note that the modulus of continuity has the bound $w^\prime_X(\gamma) \le w_X(2\gamma)$ where 
\begin{equation*}
w_X(\gamma) = \sup_{|t-s|\le \gamma} d_E(X_s,X_t)
\end{equation*}
for each $X\in D([0,T],E)$. 
Therefore, to show that a sequence in the Skorohod space is relatively compact, which is equivalent to the sequence being tight since the space is complete and separable, it suffices to show the following condition (2') instead of the condition (2) in Proposition \ref{prop:prohorov}. 

\begin{itemize}
    \item[(2')]  For each $\varepsilon>0$, we have $
\lim_{\gamma\to0} \limsup_{n\to \infty}
\mathbb{P}_n (w_X(\gamma) > \varepsilon) = 0$. 
\end{itemize}

Note that once the condition (2') is verified, combined with the condition (1) of Proposition \ref{prop:prohorov}, then all limit points of a sequence $\{\mathbb{P}_n\}_n$ are concentrated on continuous paths.  

Now we return to our current situation and recall the martingale decomposition \eqref{eq:KPZ_martingale_decomposition}. 
Our central aim is to show the tightness of each sequence. 

\begin{lemma}
\label{lem:tightness}
The sequences $\{\mathcal{X}^n_t : t \in [0, T ] \}_{ n \in \mathbb{N} } $, $\{ \mathcal{M}^n_t : t \in [0, T ] \}_{ n \in \mathbb{N} } $, $\{ \mathcal{S}^n_t : t \in [0, T ] \}_{ n \in \mathbb{N} } $ and $\{ \mathcal{B}^n_t : t \in [0, T ] \}_{ n \in \mathbb{N} } $, when the processes start from the invariant measure $\nu_n$, are tight with respect to the Skorohod topology on $D([0,T],\mathcal{S}^\prime (\mathbb{R})) $.
\end{lemma}

Here, note that the space of Schwartz distributions $\mathcal{S}^\prime(\mathbb{R})$ is metrizable, which turns out to be separable and complete with respect to the strong topology.
(See \cite[Section 2.3]{gonccalves2014nonlinear} for a precise description of the topology.) 
To prove the tightness of a sequence of processes, the following Mitoma's criterion \cite[Theorem 4.1]{mitoma1983tightness} is helpful.

\begin{proposition}[Mitoma's criterion]
\label{Mitoma}
A sequence of $\mathcal{S}^\prime (\mathbb{R} ) $-valued processes $\{ \mathcal{Y}^n_t : t \in [0, T ] \}_{ n \in \mathbb{N} } $ with trajectories in $D ([0, T ] , \mathcal{S}^\prime (\mathbb{R} ) ) $ is tight with respect to the uniform topology if, and only if, the sequence $\{ \mathcal{Y}^n_t (\varphi ) : t \in [0, T ] \}_{ n \in \mathbb{N}}$ of real-valued processes is tight with respect to the Skorohod topology on $D([0,T],\mathbb{R}) $ for any $\varphi \in \mathcal{S}(\mathbb{R}) $.   
\end{proposition}

In addition, we will make use of the following continuity criterion.
(See for example \cite[Theorem 1.2.1]{revuz2013continuous}.)

\begin{proposition}[The Kolmogorov-Chentsov criterion]
\label{prop:kolmogorov_chentsov}
Let $\{ X_t : t\in [0,T]\}$ be a Banach-valued process \Add{such that} which there exists constants $\kappa,\gamma_1\Add{,} \gamma_2>0$ \Add{satisfying} 
\begin{equation*}
\mathbb{E} \big[ \big\|X_t-X_s\big\|^{\gamma_1} \big]
\le \kappa | t-s|^{1+\gamma_2},
\end{equation*}
for any $s,t\in[0,T]$. 
Here $\| \cdot\|$ denotes the norm of the Banach space on which the process takes values. 
Then, there is a modification $\tilde{X}$ of $X$ such that 
\begin{equation*}
\mathbb{E} \bigg[ \bigg( \sup_{s\neq t} \frac{\|\tilde{X}_t -\tilde{X}_s\|}{|t-s|^{\alpha}} \bigg)^{\gamma_1} \bigg] < +\infty
\end{equation*}
for any $\alpha \in [0, \gamma_2/\gamma_1)$. 
In particular, the paths of $\tilde{X}$ are almost-surely $\alpha$-H\"{o}lder continuous. 
\end{proposition}

In what follows, we prove Lemma \ref{lem:tightness}. 
With the help of Mitoma's criterion, it suffices to show the tightness of sequences $\{ \mathcal{X}^n_t (\varphi) : t \in [0, T ] \}_{ n \in \mathbb{N} } $, $\{ \mathcal{S}^n_t (\varphi) : t \in [0, T ] \}_{ n \in \mathbb{N} } $, $\{ \mathcal{B}^n_t (\varphi) : t \in [0, T ] \}_{ n \in \mathbb{N} } $ and $\{ \mathcal{M}^n_t (\varphi) : t \in [0, T ] \}_{ n \in \mathbb{N}}$ in $D ([0,T],\mathbb{R})$ for any given test function $\varphi \in \mathcal{S} (\mathbb{R})$. 
In order to prove the tightness of a real-valued sequence $\{X^n_t:t\ge0 \}_n$ in $D([0,T],\mathbb{R})$, according to Prohorov's theorem, recall that it suffices to show two conditions, namely the condition (1) in Proposition \ref{prop:prohorov} and the condition (2'): 
\begin{equation}
\label{eq:tighness_condition_prohorov}
\lim_{\delta\to0} \limsup_{n\to \infty}
\mathbb{P}_n \bigg( \sup_{\substack{|t-s|\le \delta\\ 0\le s,t \le T}} | X^n_t-X^n_s| > \varepsilon \bigg)
= 0,
\end{equation}
for any $\varepsilon>0$. 
The condition (1) of Proposition \ref{prop:prohorov} on fixed times easily follows for our sequences. 
Hence, in what follows, our task is to verify the condition \eqref{eq:tighness_condition_prohorov} for the sequences $\{ \mathcal{X}^n_t (\varphi) : t \in [0, T ] \}_{ n \in \mathbb{N} } $, $\{ \mathcal{S}^n_t (\varphi) : t \in [0, T ] \}_{ n \in \mathbb{N} } $, $\{ \mathcal{B}^n_t (\varphi) : t \in [0, T ] \}_{ n \in \mathbb{N} } $ and $\{ \mathcal{M}^n_t (\varphi) : t \in [0, T ] \}_{ n \in \mathbb{N}}$. 
For the initial field $\{ \mathcal{X}^n_0\}_n$, it is enough to observe that by characteristic functions we can show that it converges to a Gaussian field, and in particular it is tight. 
Thus, in what follows, we focus on the tightness of the martingale, symmetric and anti-symmetric parts, from which the tightness of the fields $\{ \mathcal{X}^n_\cdot\}_n$ is deduced.

\subsubsection{Martingale Part}
Next, we deal with the martingale part. 
We will make use of the following fourth-moment estimate.  
\begin{lemma}
\label{lem:martingale_fourth_moment_estimate}
For each smooth compactly supported function $\varphi$, and for each $s,t\in[0,T]$ such that $s<t$ we have that 
\begin{equation*}
\mathbb{E}_n \big[ \big( \mathcal{M}^n_t(\varphi) - \mathcal{M}^n_s(\varphi) \big)^4 \big]
\lesssim  (t-s)^2 + n^{-3} (t-s) . 
\end{equation*}
\end{lemma}
In Lemma \ref{lem:martingale_fourth_moment_estimate} we assumed that each test function has a compact support, in order to assure the integrability of an exponential martingale in the proof, whereas we need the fourth moment estimate for functions in the Schwartz space $\mathcal{S}(\mathbb{R})$.
On the other hand, note that we can approximate any element of $\mathcal{S}(\mathbb{R})$ by smooth compactly supported functions with respect to the Sobolev \Add{$H^1(\mathbb R)$-}norm. 
Then, the tightness of the martingales, whose test function is taken from this restricted space, turns out to be sufficient for the proof by conducting the approximation in the martingale decomposition. 
\begin{proof}[Proof of Lemma \ref{lem:martingale_fourth_moment_estimate}]
The proof is based on an expansion of an exponential martingale, similarly to \cite{gonccalves2015stochastic}.
We set $\Xi_j=\eta_j + c_3\beta_n V_{\beta_n}(\eta_j)$ in the sequel. 
First, note that the process 
\begin{equation*}
\exp\big( \rho \mathcal{X}^n_t(\varphi) \big) 
= \prod_{j\in\mathbb{Z}} \exp \bigg( \frac{\rho}{\sqrt{n}} \overline{\Xi}_j T^-_{v_n^1t}\varphi^n_j \bigg) ,
\end{equation*}
is contained in the space $L^2(\nu_n)$ for sufficiently small $\rho$.  
Indeed, recalling the static estimate from Section \ref{sec:static_estimate} and noting that $b=1$, we obtain: 
\begin{equation*}
\begin{aligned}
E_{\nu_n}\big[\exp(\rho \Xi_j) \big]
= \frac{1}{Z_n} \int_{\mathbb{R}}
\exp\{ (\rho c_3 \beta_n -1) V_{\beta_n}(\eta_j)
+ (\rho+\lambda) \eta_j \} d\eta_j 
\end{aligned}
\end{equation*}
is finite when $\rho$ is sufficiently small.
In addition, note that we calculate 
\begin{equation}
\label{eq:oscillator_exponential_martingale_calculation}
\begin{aligned}
\exp\big(-\rho \mathcal{X}^n_t(\varphi)\big) (\partial_t + L_n ) \exp\big(\rho \mathcal{X}^n_t(\varphi)\big)
&= \sum_{j\in\mathbb{Z}} \big(\exp (\rho(\Xi_j-\Xi_{j+1})T^-_{v_n^1t}(\varphi^n_{j+1} -\varphi^n_j ) ) - 1\big)\\
&\quad+ \rho c_3\sqrt{n}\alpha_n\beta_n 
\sum_{j\in\mathbb{Z}} \overline{\xi}_j \overline{\xi}_{j+1}
\nabla^{1,n} T^-_{v_n^1t} \varphi^n_j ,
\end{aligned}
\end{equation}
which turns out to be in $L^2(\nu_n)$. 
Consequently, we have that the process 
\begin{equation*}
\mathrm{Exp}(\mathcal{M})^n_{s,t}
=\exp \bigg( \rho \mathcal{X}^n_t(\varphi) - \rho \mathcal{X}^n_s (\varphi)
- \int_s^t \exp(-\rho\mathcal{X}^n_r(\varphi)) 
(\partial_r + L_n)\exp(\rho \mathcal{X}^n_r(\varphi)) dr \bigg),
\end{equation*} 
is a martingale, see \cite{ethier1986markov}. 
Then, applying Taylor's theorem to the exponential functions in the integrand, we can expand $\mathrm{Exp}(\mathcal{M})^n_{s,t}$ in terms of $\rho$ as
\begin{equation}
\label{eq:exponential_martingale_expansion}
\begin{aligned}
 \mathrm{Exp}(\mathcal{M})^n_{s,t} 
=\exp \bigg( \rho 
 \mathcal{M}^n_{s,t}(\varphi) 
-\frac{\rho^2}{2!} 
\langle \mathcal{M}^n(\varphi)\rangle_{s,t} 
- \sum_{i=1}^3\frac{\rho^{2+i}}{(2+i)!} \int_s^t \mathcal{R}_i (r) dr
\bigg) ,
\end{aligned}    
\end{equation}
where 
\begin{equation*}
\begin{aligned}
\mathcal{R}_1 (t)
&= L_n (\mathcal{X}^n_t(\varphi))^3 
-3\mathcal{X}^n_t(\varphi) L_n (\mathcal{X}^n_t(\varphi))^2
+ 3(\mathcal{X}^n_t(\varphi))^2 L_n \mathcal{X}^n_t(\varphi) \\
&= \frac{\theta(n)}{n^{9/2}}\sum_{j\in\mathbb{Z}} 
(\Xi_j(t) -\Xi_{j+1}(t))^3 (T^-_{v_n^1t}\nabla^{n,1}\varphi^n_j)^3 ,
\end{aligned}
\end{equation*}
\begin{equation*}
\begin{aligned}
\mathcal{R}_2 (t)
&= L_n (\mathcal{X}^n_t(\varphi))^4 
-4\mathcal{X}^n_t(\varphi) L_n (\mathcal{X}^n_t(\varphi))^3
+ 6(\mathcal{X}^n_t(\varphi))^2 L_n (\mathcal{X}^n_t(\varphi))^2
\\&\quad
- 4(\mathcal{X}^n_t(\varphi))^3 L_n \mathcal{X}^n_t(\varphi) 
dr \\
&= \frac{\theta(n)}{n^6} \sum_{j\in\mathbb{Z}} (\Xi_j-\Xi_{j+1})^4 (T^-_{v_n^1t}\nabla^{1,n} \varphi^n_j)^4 ,
\end{aligned}
\end{equation*}
and $\mathcal{R}_3$ is a reminder term. 
Here we used the short-hand notation $\mathcal{M}^n_{s,t}(\varphi) = \mathcal{M}^n_t(\varphi) - \mathcal{M}^n_s(\varphi)$ and $\langle \mathcal{M}^n (\varphi)\rangle_{s,t}
= \langle \mathcal{M}^n (\varphi)\rangle_{t}
-\langle \mathcal{M}^n (\varphi)\rangle_{s}$. 
Moreover, we see that the bound 
\begin{equation*}
\| \mathcal{R}_1\|_{L^2(\nu_n)}
\lesssim \frac{1}{n^{3/2}} 
 \| \Xi_j \|_{L^{6}(\nu_n)}^3 
\bigg( \frac{1}{n} \sum_{j\in\mathbb{Z}} (\nabla^{1,n}\varphi^n_j)^3 \bigg) ,
\end{equation*}
and 
\begin{equation*}
\| \mathcal{R}_2\|_{L^1(\nu_n)}
\lesssim \frac{1}{n^5}  
 \| \Xi_j \|_{L^{4}(\nu_n)}^4 
\bigg( \frac{1}{n} \sum_{j\in\mathbb{Z}} (\nabla^{1,n}\varphi^n_j)^4 \bigg) ,
\end{equation*}
which will be used below. 
Noting that $E_{\nu_n}[\mathrm{Exp}(\mathcal{M})^n_{s,t}]=1$, take the expectation of the identity \eqref{eq:exponential_martingale_expansion} and expand the exponential function in $\rho$. 
This enables us to obtain the desired bound.
In the identity \eqref{eq:exponential_martingale_expansion}, we differentiate with respect to $\rho$ four times and then take $\rho=0$.
Moreover, note that the function $\rho \mapsto e^{F(\rho)}$ satisfies 
\begin{equation*}
\frac{d^4}{d\rho^4} e^{F(\rho)} = \big( 
(F^\prime(\rho))^4
+ 6 (F^\prime(\rho))^2 F^{\prime\prime}(\rho)
+ 3 (F^{\prime\prime}(\rho))^2 
+ 4 F^\prime(\rho) F^{(3)}(\rho)
+ F^{(4)}(\rho) \big)e^{F(\rho)} . 
\end{equation*}
Now, we take the expectation of the identity \eqref{eq:exponential_martingale_expansion} to obtain
\begin{equation*}
1 = E_{\nu_n} [\exp(F(\rho))], 
\end{equation*}
where $F(\rho)$ is the quantity inside the parentheses of the identity \eqref{eq:exponential_martingale_expansion}. 
Then, we differentiate the last identity four times and set $\rho=0$: 
\begin{equation*}
\begin{aligned}
0 
&= \frac{d^4}{d\rho^4} E_{\nu_n}[\exp(F(\rho))] \bigg|_{\rho=0}\\
&= E_{\nu_\rho}[(F^\prime(\rho))^4
+ 6 (F^\prime(\rho))^2 F^{\prime\prime}(\rho)
+ 3 (F^{\prime\prime}(\rho))^2 
+ 4 F^\prime(\rho) F^{(3)}(\rho)
+ F^{(4)}(\rho)]\\
&= E_{\nu_n}[\mathcal{M}^n_{s,t}(\varphi)^4]
+ 6E_{\nu_n}[\mathcal{M}^n_{s,t}(\varphi)^2 
\langle \mathcal{M}^n (\varphi) \rangle_{s,t}]
+ 3 E_{\nu_n}[\langle \mathcal{M}^n (\varphi) \rangle_{s,t}^2] \\
&\quad+ 4 E_{\nu_n}\Big[ \mathcal{M}^n_{s,t}(\varphi) \int_s^t \mathcal{R}_1(r) dr \Big]
+ E_{\nu_n} \Big[\int_s^t \mathcal{R}_2(r) dr \Big] . 
\end{aligned}
\end{equation*}
Then, we have the following estimate. 
\begin{equation*}
\begin{aligned}
E_{\nu_n}[\mathcal{M}^n_{s,t}(\varphi)^4]
&= \bigg| 6E_{\nu_n}[\mathcal{M}^n_{s,t}(\varphi)^2 
\langle \mathcal{M}^n (\varphi) \rangle_{s,t}]
+ 3 E_{\nu_n}[\langle \mathcal{M}^n (\varphi) \rangle_{s,t}^2] \\
&\quad+ 4 \int_s^t E_{\nu_n} [\mathcal{M}^n_{s,t}(\varphi)  \mathcal{R}_1(r)] dr
+ \int_s^t E_{\nu_n} [\mathcal{R}_2(r)] dr \bigg| \\
&\le 3A E_{\nu_n}[\mathcal{M}^n_{s,t}(\varphi)^4]
+ (3A^{-1}+3) E_{\nu_n}[\langle \mathcal{M}^n (\varphi) \rangle_{s,t}^2]\\
&\quad+ 2 (t-s) E_{\nu_n}[\mathcal{M}^n_{s,t}(\varphi)^2] 
+ (t-s) \big( 2\|\mathcal{R}_1 \|_{L^2(\nu_n)}^2  
+ \|\mathcal{R}_2\|_{L^1(\nu_n)} \big) ,
\end{aligned}
\end{equation*}
for any $A>0$ where we used Young's inequality. 
Set $A=1/6$ and move $E[\mathcal{M}^n_{s,t}(\varphi)^4]$ on the right-hand side to the left-hand side. 
Moreover, recall that the quadratic variation is given by \eqref{eq:qv_scalar}, which yields the bound
\begin{equation*}
E_{\nu_n} [\langle \mathcal{M}^n(\varphi)\rangle_{s,t}^k]
\lesssim (t-s)^k, 
\end{equation*}
for each $k=1,2$. 
Hence we obtain the desired \Add{bound} according to the moment estimates for $\mathcal{R}_1$ and $\mathcal{R}_2$. 
\end{proof}

With this fourth moment bound at hand, we show the tightness of the martingale part. 
This is established as follows. 
\begin{equation*}
\begin{aligned}
\mathbb{P}_n \bigg( \sup_{\substack{|s-t|\le \delta\\0\le s,t\le T} } 
\big| \mathcal{M}^n_t(\varphi)- \mathcal{M}^n_s (\varphi) \big|
> \varepsilon \bigg)
&\le \varepsilon^{-4} \mathbb{E}_n \bigg[ \sup_{\substack{|s-t|\le \delta\\0\le s,t\le T}} 
\big| \mathcal{M}^n_t(\varphi)- \mathcal{M}^n_s (\varphi) \big|^4 \bigg] \\
&\lesssim \varepsilon^{-4} \delta^{-1} \mathbb{E}_n \big[ \big( \mathcal{M}^n_\delta(\varphi)\big)^4 \big].  
\end{aligned}
\end{equation*}
Here, in the second inequality, we used Doob's inequality and stationarity. 
Now, recall that the fourth moment bound (Lemma \ref{lem:martingale_fourth_moment_estimate}) yields 
\begin{equation*}
\delta^{-1} \mathbb{E}_n \big[ \big( \mathcal{M}^n_\delta(\varphi)\big)^4 \big]
\lesssim \delta + n^{-3} , 
\end{equation*}
which vanishes as $n\to\infty$ and then $\delta\to0$. 
Hence the condition \eqref{eq:tighness_condition_prohorov} is verified and we complete the proof of the tightness of the martingale part.

\subsubsection{Symmetric Part}
For the symmetric part, note that a direct $L^2$-computation yields 
\begin{equation*}
\begin{aligned}
\mathbb{E}_n \big[\big|\mathcal{S}^n_{t_1}(\varphi)
- \mathcal{S}^n_{t_2} (\varphi)\big|^2 \big]
&\lesssim  |t_1-t_2|
\int_{t_1}^{t_2} \frac{1}{n} \sum_{j\in\mathbb{Z}} 
E_{\nu_n} [(\overline{\eta}_j+ c_3\beta_n \overline{\zeta}_j )^2] 
(\Delta^n T^-_{v_n^1s}\varphi^n_j)^2 ds \\
&\lesssim  | t_1-t_2|^2 \| \partial_x^2 \varphi \|^2_{L^2(\mathbb{R})}. 
\end{aligned}
\end{equation*}
By the Kolmogorov-Chentsov criterion (Proposition \ref{prop:kolmogorov_chentsov}), combined with the continuity of $\mathcal{S}^n_\cdot$, the condition \eqref{eq:tighness_condition_prohorov} is verified, from which we deduce the tightness of the symmetric part.

\subsubsection{Anti-symmetric Part}
Finally, we are in a position to show tightness of the anti-symmetric part $\mathcal{B}^n_t(\varphi)$. 
A key ingredient for the proof is the so-called second-order Boltzmann-Gibbs principle, which enables us to replace a quadratic term by its local average. 
For each real sequence $(g_j)_{j \in \mathbb{Z}}$, we define its local average as follows.
\begin{equation*}
\overrightarrow{g}^\ell_j = \frac{1}{\ell} \sum_{i=0}^{\ell-1} g_{j+i}. 
\end{equation*}

\begin{proposition}[The second-order Boltzmann-Gibbs principle]
\label{prop:2BG}
For any $T>0$ and $\varphi\in\mathcal{S}(\mathbb{R})$, it holds that \begin{equation*}
\begin{aligned}
\mathbb{E}_n \bigg[ \sup_{0 \le t \le T} \bigg| 
\int_0^t \sum_{j\in \mathbb{Z}} 
\big( \overline{\xi}_j(s) \overline{\xi}_{j+1}(s) 
- (\overrightarrow{\xi}^\ell_{j}(s))^2 \big)
\nabla^{n,1} T^-_{v_n^1s}\varphi^n_j ds \bigg|^2 \bigg] 
\lesssim 
\bigg(\frac{T\ell}{n} + \frac{T^2n}{\ell^2} \bigg)  \|\partial_x\varphi\|^2_{L^2(\mathbb{R})}.
\end{aligned}
\end{equation*}
\end{proposition}
The proof is completely analogous to~\cite[Theorem 1]{gonccalves2017second} so that we omit it here. 
With this result at hand, we can show tightness of the anti-symmetric part. 
Indeed, we have by Proposition~\ref{prop:2BG} and stationarity that 
\begin{equation*}
\begin{aligned}
&\mathbb{E}_n \bigg[\bigg| \mathcal{B}^n_{t_2}(\varphi)
- \mathcal{B}^n_{t_1}(\varphi) 
- c_3 \int_{t_1}^{t_2} \sum_{j\in\mathbb{Z}} \big(\overrightarrow{\xi}^\ell_{j}(s)\big)^2  
\nabla^{1,n} T^-_{v_n^1s} \varphi^n_j ds \bigg|^2 \bigg]\\
&\quad \lesssim  \bigg( \frac{(t_2-t_1)\ell}{n} + \frac{(t_2-t_1)^2n}{\ell^2} \bigg)\| \partial_x \varphi \|^2_{L^2(\mathbb{R})} .
\end{aligned}
\end{equation*}
On the other hand, a direct $L^2$-computation gives 
\begin{equation*}
\begin{aligned}
\mathbb{E}_n \bigg[ \bigg| \int_{t_1}^{t_2} 
\sum_{j\in\mathbb{Z}} 
\big(\overrightarrow{\xi}^\ell_{j}(s)\big)^2  
\nabla^{1,n} T^-_{v_n^1s} \varphi^n_j ds
\bigg|^2 \bigg]
\le  \frac{(t_2-t_1)^2n}{\ell} \|\partial_x \varphi\|^2_{L^2(\mathbb{R})}. 
\end{aligned}
\end{equation*}
When $1/n^2 \le t_2-t_1 \le 1$, we take $\ell$ with order $(t_2-t_1)^{1/2}n$ to obtain 
\begin{equation*}
\mathbb{E}_n\big[\big|\mathcal{B}^n_{t_2}(\varphi)
- \mathcal{B}^n_{t_1}(\varphi) \big|^2 \big] 
\lesssim  (t_2-t_1)^{3/2}
\| \partial_x \varphi \|^2_{L^2(\mathbb{R})} .
\end{equation*}
On the other hand, when $t_2-t_1 \le 1/n^2$, a direct estimate brings us 
\begin{equation*}
\mathbb{E}_n\big[\big|\mathcal{B}^n_{t_2}(\varphi)
- \mathcal{B}^n_{t_1}(\varphi) \big|^2 \big] 
\lesssim  (t_2-t_1)^2 n \| \partial_x \varphi \|^2_{L^2(\mathbb{R})} 
\lesssim  (t_2-t_1)^{3/2} \| \partial_x \varphi \|^2_{L^2(\mathbb{R})} .
\end{equation*}
This ends the proof of tightness for the anti-symmetric part by the Kolmogorov-Chentsov criterion (Propsition \ref{prop:kolmogorov_chentsov}) and continuity of the process.

\subsection{Identification of Limit Points}
Recall the martingale decomposition \eqref{eq:KPZ_martingale_decomposition}. 
We have already proved in Lemma \ref{lem:tightness} that the sequences $\{ \mathcal{X}^n_t : t \in [0, T] \}_{ n \in \mathbb{N} } $, $\{ \mathcal{M}^n_t : t \in [0, T] \}_{ n \in \mathbb{N} } $, $\{ \mathcal{S}^n_t : t \in [0, T] \}_{ n \in \mathbb{N} } $ and $\{ \mathcal{B}^n_t : t \in [0, T] \}_{ n \in \mathbb{N} } $ are tight in $D([0,T], \mathcal{S}^\prime(\mathbb{R}))$.
Let $\mathscr{Q}^n$ be the distribution of 
\begin{equation*}
\{
(\mathcal{X}^n_t, \mathcal{M}^n_t, 
\mathcal{S}^n_t, \mathcal{B}^n_t):
t\in [0,T]
\}.
\end{equation*}
We proved that there exists a subsequence $n$, which is denoted by the same letter with an abuse of notation, such that $\{\mathscr{Q}^n\}_{n}$ converges to a limit point $\mathscr{Q}$. 
We let $\mathcal{X}, \mathcal{M}, \mathcal{S}$ and $\mathcal{B}$ be the respective limits in distribution of each component. 
Since the tightness is shown in the uniform topology in $D([0,T],\mathcal{S}^\prime(\mathbb{R}))$, then these limiting processes almost surely have continuous trajectories.
In what follows, we identify these limit points as  stationary energy solutions of the stochastic Burgers equation in the sense of Definition \ref{def:energysol}. Since this solution is unique-in-law, convergence then follows.

\subsubsection{Martingale Part}
Now, we shift to the martingale part. 
Recall the expression \eqref{eq:qv_scalar} when $\theta(n)=n^2$, $\mathfrak{u}_n=c_3\beta_n$ and $v_n=v_n^1$: 
\begin{equation*}
\begin{aligned}
\langle \mathcal{M}^n (\varphi) \rangle_t 
&= \frac{1}{2n} \int_0^t \sum_{j \in \mathbb{Z}} 
(\eta_j(s) -\eta_{j+1}(s))^2 \big( \nabla^{n,1} T^-_{v_n^1s} \varphi^n_j \big)^2 ds \\
&\quad+ \frac{(c_3\beta_n)^2 }{2n} \int_0^t \sum_{j \in \mathbb{Z}} 
(\zeta_j(s) -\zeta_{j+1}(s))^2 \big( \nabla^{n,1} T^-_{v_n^1s} \varphi^n_j \big)^2 ds \\
&\quad+ \frac{c_3 \beta_n}{n}
\int_0^t \sum_{j\in\mathbb{Z}} 
(\eta_j(s) - \eta_{j+1}(s)) (\zeta_j(s) -\zeta_{j+1}(s))
(\nabla^{n,1} T^-_{v_n^1s}\varphi^n_j)^2 ds.
\end{aligned}
\end{equation*}
By Markov's inequality and the Schwarz's inequality,  
\begin{equation*}
\lim_{n\to\infty}
\mathbb{P}_n \bigg(\sup_{0\le t\le T} \Big| \langle \mathcal{M}^n(\varphi)\rangle_t -t\| \partial_x \varphi \|^2_{L^2(\mathbb{R})}  \Big| > \varepsilon \bigg)
= 0 
\end{equation*}
for any $\varepsilon>0$.  
In particular, the sequence of processes $\{\langle \mathcal{M}^n (\varphi ) \rangle_t:t\in [0,T] \}_n $ converges in distribution on $D([0,T],\mathbb{R})$ to a deterministic path $\{ t\| \partial_x \varphi \|^2_{L^2(\mathbb{R})} :t\in [0,T] \}$ as $n\to\infty$.

Then we can show that any limit point $\mathcal{M}$ is a continuous martingale whose quadratic variation is $ t \| \partial_x \varphi\|^2_{L^2(\mathbb{R})}$ in the following way. 
First, note that the limit point $\mathcal{M}$ is a martingale since is is obtained as a limit of martingales with respect to the uniform topology. 
Moreover, note that by the triangle inequality and Doob's inequality, we have the bound 
\begin{equation*}
\begin{aligned}
\mathbb{E}_n \bigg[\sup_{0\le s \le t} \big| \mathcal{M}^n_s(\varphi)
- \mathcal{M}^n_{s-}(\varphi) \big| \bigg]
\le 2 \mathbb{E}_n \bigg[\sup_{0\le s \le t} | \mathcal{M}^n_s(\varphi) |^2 \bigg]^{1/2} 
\le 2 \mathbb{E}_n \big[\langle \mathcal{M}^n(\varphi)\rangle_t \big]^{1/2},
\end{aligned}
\end{equation*}
and the utmost right-hand side is bounded by a constant which is independent of $n$. 
Therefore, by Corollary VI.6.30 of \cite{jacod2013limit}, we obtain the convergence $(\mathcal{M}^n(\varphi), \langle\mathcal{M}^n(\varphi)\rangle)\to 
(\mathcal{M}(\varphi), \langle\mathcal{M}(\varphi)\rangle)$ in distribution.
Combining with the convergence of the quadratic variation, we conclude that $\langle \mathcal{M}(\varphi)\rangle_t= t\|\partial_x \varphi\|^2_{L^2(\mathbb{R})}$.

\subsubsection{Symmetric Part}
For the symmetric part, we can easily show that 
\begin{equation*}
\begin{aligned}
\mathcal{S}_t(\varphi)
= \frac{1}{2} \int_0^t \mathcal{X}_s(\partial_x^2\varphi) ds .
\end{aligned}
\end{equation*}

\subsubsection{Anti-symmetric Part}
Now, it suffices to characterize the limit point for the anti-symmetric part to complete the proof of Theorem~\ref{thm:SBE_derivation}. 
Define a modified fluctuation field 
\begin{equation*}
\tilde{\mathcal{X}}^n_t(\varphi)
= \frac{1}{\sqrt{n}} \sum_{j\in\mathbb{Z}} 
\overline{\xi}_j(t) T^-_{v_n^1t} \varphi^n_j .
\end{equation*}
Note that for fixed $j\in\mathbb{Z}$, we have $\tilde{\mathcal{X}}^n_t (\iota_\varepsilon(\frac{j-v_n^1t}{n},\cdot))
= \sqrt{n}\overrightarrow{\xi}^{\varepsilon n}_j$,
recalling the definition of $\iota_\varepsilon(x;\cdot)$ given below \eqref{eq:def_quadratic_function_approximation}. 
Hence, we have that 
\begin{equation*}
\begin{aligned}
\int_0^t \sum_{j\in\mathbb{Z}}
\big( \overrightarrow{\xi}^{\varepsilon n}_j(s)\big)^2 
\nabla^{1,n} T^-_{v_n^1s} \varphi^n_j ds
= \frac{1}{n}\int_0^t \sum_{j\in\mathbb{Z}}
\tilde{\mathcal{X}}^n_s \big( \iota_{\varepsilon}({\textstyle \frac{j-v_n^1s}{n}};\cdot) \big)^2
\nabla^{1,n} T^-_{v_n^1s} \varphi^n_j ds.
 \end{aligned}
\end{equation*}
In addition, set $\Xi_j = \eta_j + c_3 \beta_n \zeta_j$. 
Recall that $\mathcal{X}^n$ is the fluctuation fields of the conserved quantity $\Xi_j$ so that  
\begin{equation*}
\begin{aligned}
\int_0^t \sum_{j\in\mathbb{Z}}
\big(\overrightarrow{\Xi}^{\varepsilon n}_j(s)\big)^2 
\nabla^{1,n} T^-_{v_n^1s} \varphi^n_j ds
= \frac{1}{n}\int_0^t \sum_{j\in\mathbb{Z}}
\mathcal{X}^n_s \big( \iota_{\varepsilon}({\textstyle \frac{j-v_n^1s}{n}};\cdot) \big)^2
\nabla^{1,n} T^-_{v_n^1s} \varphi^n_j ds .
\end{aligned}
\end{equation*}
The difference between these quantities is estimated by the following lemma. 
\begin{lemma}
\label{lem:oscillator_quadratic_field_replacement}
We have that 
\begin{equation*}
\begin{aligned}
&\mathbb{E}_n \bigg[\sup_{0\le t\le T} \bigg|
\int_0^t \sum_{j\in\mathbb{Z}} 
\big(  \big(\overrightarrow{\xi}^{\ell}_j(s)\big)^2 
- \big(\overrightarrow{\Xi}^{\ell}_j(s)\big)^2 \big)
\nabla^{1,n}\varphi^n_j ds
\bigg|^2 \bigg]
\lesssim T^2 \bigg( \frac{n\beta_n^2}{\ell} + \frac{n}{\ell^2} \bigg) 
\|\partial_x\varphi\|^2_{L^2(\mathbb{R})} .
\end{aligned}
\end{equation*}
\end{lemma}
\begin{proof}
By Schwarz's inequality and Young's inequality, for any $A>0$ we have the bound
\begin{equation*}
\begin{aligned}
& \mathbb{E}_n \bigg[\sup_{0\le t\le T} \bigg|
\int_0^t \sum_{j\in\mathbb{Z}} 
\big(  \big(\overrightarrow{\xi}^{\ell}_j(s)\big)^2 
- \big(\overrightarrow{\Xi}^{\ell}_j(s)\big)^2 \big)
\nabla^{1,n}\varphi^n_j ds
\bigg|^2 \bigg]\\
&\quad\le T^2 \sum_{j\in\mathbb{Z}} 
E_{\nu_n} \big[ \big\{  \big(\overrightarrow{\xi}^{\ell}_j\big)^2 
- \big(\overrightarrow{\Xi}^{\ell}_j\big)^2 -\Sigma_\ell \big\}^2 \big] (\nabla^{1,n}\varphi^n_j)^2  \\
&\quad\le \frac{T^2}{A} \sum_{j\in\mathbb{Z}} 
E_{\nu_n} \big[\big(\overrightarrow{\xi}^{\ell}_j - \overrightarrow{\Xi}^{\ell}_j \big)^2 \big] (\nabla^{1,n}\varphi^n_j)^2
+ AT^2 
\sum_{j\in\mathbb{Z}} 
E_{\nu_n} \big[ \big( \overrightarrow{\xi}^{\ell}_j + \overrightarrow{\Xi}^{\ell}_j\big)^2 \big] (\nabla^{1,n}\varphi^n_j)^2 \\
&\qquad+ 2(\Sigma_\ell)^2 \sum_{j\in\mathbb{Z}} (\nabla^{1,n}\varphi^n_j)^2 \\
&\quad\lesssim T^2 \bigg( \frac{n\beta_n^4}{A\ell} + \frac{An}{\ell} + \frac{n}{\ell^2}\bigg)
\| \partial_x\varphi\|^2_{L^2(\mathbb{R})} ,
\end{aligned}
\end{equation*}
where $\Sigma_\ell= E_{\nu_n}[(\overrightarrow{\xi}^\ell_j)^2 - (\overrightarrow{\Xi}^\ell_j)^2]$ which satisfies the bound $\Sigma_\ell^2\lesssim \ell^{-2} $.
Here in the first term of the utmost right-hand side of the last display, we used the Taylor expansion $\xi_j-\Xi_j=\beta_n^2\varepsilon_j$ for some reminder term $\varepsilon_j$, which is not necessarily centered.  
Now, choose $A=\beta_n^{-2}$ in the above estimate to complete the proof. 
\end{proof}

With this estimate at hand, now we identify the limit points of the anti-symmetric part. 
Note that the sequence $\{ \mathcal{X}^n\}_n$ is tight, which subsequently converges to some $\mathcal{X}$ which clearly satisfies the condition \textbf{(S)}. 
Therefore, we obtain the limit 
\begin{equation*}
\mathcal{A}^\varepsilon_{s,t} (\varphi)
= \lim_{n\to\infty} \frac{1}{n}
\int_s^t \sum_{j\in\mathbb{Z}}
\mathcal{X}^n_r \big( \iota_{\varepsilon}({\textstyle \frac{j-v_n^1r}{n}};\cdot) \big)^2
\nabla^{n,1} T^-_{v_n^1r} \varphi^n_j dr, 
\end{equation*}
where $\mathcal{A}^\varepsilon_{s,t}(\varphi)$ in the left-hand side of the last display coincides with the process we defined in \eqref{eq:def_quadratic_function_approximation} with $u=\mathcal{X}$.
Here note that the convergence does not follow immediately since the function $\iota_\varepsilon$ does not belong to $\mathcal{S}(\mathbb{R})$. 
This function, however, can be approximated by elements in $\mathcal{S}(\mathbb{R})$ so that the convergence is justified.
(See \cite[Section~5.3]{gonccalves2014nonlinear} for details.) 
Now, by Proposition~\ref{prop:2BG} and Lemma \ref{lem:oscillator_quadratic_field_replacement}, we have that 
\begin{equation*}
\begin{aligned}
&\mathbb{E}_n \bigg[\bigg|
\mathcal{B}^n_{t}(\varphi) 
- \mathcal{B}^n_{s}(\varphi)
- \int_s^t \sum_{j\in\mathbb{Z}} 
\big( \overrightarrow{\Xi}^{\varepsilon n}_j(r)\big)^2 
\nabla^{1,n} T^-_{v_n^1r} \varphi^n_j dr\bigg|^2 \bigg]\\
&\quad \lesssim \bigg(\frac{(t-s)\varepsilon n}{n} 
+ \frac{(t-s)^2n}{(\varepsilon n)^2} 
+ \frac{n\beta_n^2}{\varepsilon n} \bigg)
\| \partial_x \varphi \|^2_{L^2(\mathbb{R})} .
\end{aligned}
\end{equation*}
Then, let $n\to+\infty$ to obtain 
\begin{equation}
\label{eq:energy_condition_estimate}
\mathbb{E} \big[\big| 
\mathcal{B}_{t}(\varphi) 
- \mathcal{B}_{s} (\varphi)
- c_3 \mathcal{A}^\varepsilon_{t,s}(\varphi)
\big|^2 \big] 
\lesssim  \varepsilon(t-s)
\| \partial_x \varphi \|^2_{L^2(\mathbb{R})} .
\end{equation}
The triangle inequality deduces the energy condition \textbf{(EC)}.
As a result, by Proposition~\ref{prop:nonlinear} we get the existence of the limit 
\begin{equation*}
\mathcal{A}_{t} (\varphi)
= \lim_{\varepsilon\to0}
\mathcal{A}^\varepsilon_{0,t}
 (\varphi). 
\end{equation*}
Moreover, the estimate \eqref{eq:energy_condition_estimate} with $s=0$ shows that $\mathcal{B}= c_3\mathcal{A}$. 

Finally, we note that all the above estimates hold also for the reversed process $\{ \mathcal{X}^n_{T-t} : t \in [0,T]\}$ by repeating the argument for the dynamics generated by the adjoint operator $L^*_n$, and thus the third condition of Definition \ref{def:energysol} is satisfied. 
From this it follows that the limiting process $\mathcal{X} $ is the energy solution of the SB equation \eqref{eq:sbe_main_theorem} and thus we complete the proof of Theorem~\ref{thm:SBE_derivation}.

\section{Proof of Theorem \ref{thm:3/2L_derivation}: The 3/2-L\'{e}vy Regime}
\label{sec:3/2Levy}
\subsection{The Quadratic Martingale}
Hereafter we consider the case~(ii) in \eqref{eq:two_regimes}. 
The fluctuation field \eqref{eq:xi_fluctuation_def}, which we are concerned with, after multiplied by $-\lambda^{-1}$, noting that we can multiply any constant since it does not affect the result, is given by 
\begin{equation*}
\mathcal{X}^n_t(\varphi) = \frac{1}{\sqrt{n}}\sum_{j\in \mathbb{Z}}
\big(\overline{\zeta}_j(t) - \lambda \overline{\eta}_j(t)\big) 
\varphi^n_j , 
\end{equation*}
for each $\varphi\in\mathcal{S}(\mathbb{R})$ where recall that we defined $\varphi^n_j=\varphi(j/n)$.  
Here we used the same notation for the field as before by an abuse of notation. 
In what follows, we are concerned with the correlation function for simplicity, but note that we can also study fluctuation fields themselves\Add{, see \cite{bernardin2018interpolation} for this setting}. 
Recall the correlation function~\eqref{eq:correlation}. 
Fix any compactly supported function $g:\mathbb{R}\to\mathbb{R}$.
Let $\mathscr{S}^n_t(\cdot)\in D([0,T],\mathcal{S}^\prime(\mathbb{R}))$ be a field which is defined by
\begin{equation*}
\begin{aligned}
\mathscr{S}^n_t(f)
= \frac{1}{n}\sum_{j,j^\prime\in\mathbb{Z}}
S_{t}(j^\prime-j) 
g\bigg(\frac{j}{n}\bigg) f\bigg(\frac{j^\prime}{n}\bigg)
= \frac{1}{2}\mathbb{E}_n\big[\mathcal{X}^n_0(g)
\mathcal{X}^n_t(f) \big] ,
\end{aligned}
\end{equation*}
for each $f\in C^\infty_c(\mathbb{R})$. 
From \eqref{eq:action_antisymmetric_X} taking $\mathfrak{u}_n=\mathfrak{u}^2_n=-1/\lambda$ and $v_n=v^2_n=0$, the action of the generator $L_n$ on the field $\mathcal{X}^n_\cdot$ is given by 
\begin{equation}
\label{eq:generator_action_levy_original}
\begin{aligned}
L_n \mathcal{X}^n_s(\varphi)  
= \frac{\theta(n)}{2n^2} 
\mathcal{X}^n_t(\Delta^n \varphi^n) 
+ \frac{\theta(n)\alpha_n}{n^{3/2}} 
\sum_{j \in \mathbb{Z}}
\overline{\xi}_j(s) \overline{\xi}_{j+1}(s)
\nabla^{1,n} \varphi^n_j  
+ E^n_t. 
\end{aligned}
\end{equation}
Here the error term $E^n_t$ satisfies the bound 
\begin{equation}
\label{eq:levy_error_bound_generator}
\mathbb{E}_n\bigg[ \sup_{0\le t\le T} \bigg| \int_0^t E^n_s ds\bigg|^2 \bigg]
\lesssim  \frac{\theta(n)^2\alpha_n^2}{n^4},
\end{equation}
using the fact that $2+2\lambda\mathfrak{u}_n=0$.
In the sequel, we introduce a quadratic field $\mathcal{Q}^n_\cdot$ which belongs to $D([0,T],\mathcal{S}^\prime(\mathbb{R}^2))$ and it is defined on smooth symmetric functions $h:\mathbb{R}^2 \to \mathbb{R}$ by 
\begin{equation}
\label{eq:definition_quadratic_field_Q}
\mathcal{Q}^n_t(h) = \frac{1}{n} \sum_{\substack{j,j^\prime\in\mathbb{Z},\\ j\neq j^\prime}}
\overline{\xi}_j(t) \overline{\xi}_{j^\prime}(t) 
h \bigg( \frac{j}{n}, \frac{j^\prime}{n} \bigg).
\end{equation}
Moreover, set $\mathscr{Q}^n_t(h) 
= (1/2) \mathbb{E}_n \big[ \mathcal{X}^n_0(g) \mathcal{Q}^n_t(h) \big] $. 
From \eqref{eq:generator_action_levy_original} we have that 
\begin{equation}
\label{eq:martingale_decomposition_original}
\begin{aligned}
\frac{d}{dt} \mathscr{S}^n_t(f)
&= \frac{1}{2} \mathbb{E}_n \big[ \mathcal{X}^n_0(g) L_n \mathcal{X}^n_t(f) \big]  \\
&= \frac{\theta(n)}{2n^2} \mathscr{S}^n_t(\Delta^nf)
+ \frac{2\theta(n)\alpha_n}{n^{3/2}}
\mathscr{Q}^n_t(\nabla^{n,1} f\otimes \delta) 
+ \tilde{E}^n_t.
\end{aligned}
\end{equation}
Above, 
\begin{equation*}
(\nabla^{n,1} f \otimes \delta)
\bigg(\frac{j}{n},\frac{j^\prime}{n} \bigg)
= \frac{n}{2} \nabla^{1,n}f \bigg(\frac{j}{n}\bigg) 
\mathbf{1}_{j^\prime=j+1}
+ \frac{n}{2} \nabla^{1,n}f \bigg( \frac{j-1}{n}\bigg) 
\mathbf{1}_{j^\prime=j-1} ,
\end{equation*}
and $\tilde{E}^n_t$ is an error term which satisfies the bound \eqref{eq:levy_error_bound_generator}.  
In the sequel, let $\Delta^n$ be the two-dimensional discrete Laplacian defined by
\begin{equation*}
\Delta^n h \bigg(\frac{j}{n},\frac{j^\prime}{n}\bigg) 
= n^2 \bigg[ 
h\bigg(\frac{j+1}{n},\frac{j^\prime}{n}\bigg) 
+ h\bigg(\frac{j-1}{n},\frac{j^\prime}{n}\bigg) 
+ h\bigg(\frac{j}{n},\frac{j^\prime+1}{n}\bigg) 
+ h\bigg(\frac{j}{n},\frac{j^\prime-1}{n}\bigg) 
- 4 h\bigg(\frac{j}{n},\frac{j^\prime}{n}\bigg) 
\bigg]. 
\end{equation*}
Above we have used the same symbol for the one-dimensional Laplacian with an abuse of notation.
Moreover, we define discrete derivative operators by 
\begin{equation*}
\mathscr{A}_n h \bigg(\frac{j}{n},\frac{j^\prime}{n}\bigg) 
= n \bigg[ 
h \bigg(\frac{j}{n},\frac{j^\prime-1}{n}\bigg)
+ h \bigg(\frac{j-1}{n},\frac{j^\prime}{n}\bigg)
- h \bigg(\frac{j}{n},\frac{j^\prime+1}{n}\bigg)
- h \bigg(\frac{j+1}{n},\frac{j^\prime}{n}\bigg)\bigg] ,
\end{equation*}
and
\begin{equation*}
\mathscr{D}_n h \bigg(\frac{j}{n}\bigg) 
=\frac{n}{2} \bigg[
h\bigg(\frac{j+1}{n},\frac{j}{n}\bigg)
+ h\bigg(\frac{j}{n},\frac{j+1}{n}\bigg)
- h\bigg(\frac{j-1}{n},\frac{j}{n}\bigg)
- h\bigg(\frac{j}{n},\frac{j-1}{n}\bigg)
\bigg] .
\end{equation*}
Here note that $\mathscr{A}_nh(j/n,j/n)=-2\mathscr{D}_nh(j/n)$. 
Additionally, define 
\begin{equation*}
\begin{aligned}
\tilde{\mathscr{D}}_nh 
\bigg(\frac{j}{n},\frac{j^\prime}{n}\bigg)
&= n^2 \bigg[ \frac{1}{2}\tilde{\mathscr{E}}_nh 
\bigg( \frac{j}{n}\bigg) - \frac{1+2\alpha_n}{2}
\tilde{\mathscr{F}}_nh \bigg( \frac{j}{n} \bigg)
\bigg]\mathbf{1}_{j^\prime=j+1} \\
&\quad+ n^2 \bigg[ \frac{1}{2}\tilde{\mathscr{E}}_nh 
\bigg( \frac{j^\prime}{n}\bigg) - \frac{1+2\alpha_n}{2}
\tilde{\mathscr{F}}_nh \bigg( \frac{j^\prime}{n} \bigg)
\bigg]\mathbf{1}_{j^\prime=j-1} ,
\end{aligned}
\end{equation*}
where 
\begin{equation*}
\tilde{\mathscr{E}}_nh \bigg(\frac{j}{n}\bigg)
= \frac{1}{2}\bigg[ 
h\bigg( \frac{j}{n},\frac{j+1}{n}\bigg)
+ h\bigg( \frac{j+1}{n},\frac{j}{n}\bigg)
- 2h\bigg( \frac{j}{n},\frac{j}{n}\bigg) \bigg],
\end{equation*}
and 
\begin{equation*}
\tilde{\mathscr{F}}_nh \bigg(\frac{j}{n}\bigg)
= h\bigg( \frac{j+1}{n},\frac{j+1}{n}\bigg)
- h\bigg( \frac{j}{n},\frac{j}{n}\bigg) .
\end{equation*}
Define 
\begin{equation*}
{\mathscr{E}}_nh\bigg(\frac{j}{n},\frac{j^\prime}{n}\bigg)
= h\bigg(\frac{j+1}{n},\frac{j^\prime}{n}\bigg)-h\bigg(\frac{j}{n},\frac{j^\prime}{n}\bigg) ,
\end{equation*} 
and note that 
\begin{equation*}
\mathscr{E}_nh\bigg(\frac{j}{n},\frac{j}{n}\bigg)
=\tilde{\mathscr{E}}_nh \bigg(\frac{j}{n}\bigg).
\end{equation*}

For a function $h:\mathbb Z^2\to\mathbb R$ we use the notation $\|h\|_{2,n}$ for the discrete $L^2$-norm of $h$, i.e., 
\begin{equation}
\|h\|_{2,n}\coloneqq\sqrt{\frac{1}{n^2}\sum_{j,j'\in\mathbb Z}h\bigg(\frac jn,\frac {j'}{n}\bigg)^2}.\end{equation}

Then, the action of the generator on $\mathcal{Q}^n$ which is defined in \eqref{eq:definition_quadratic_field_Q} is calculated as follows. 

\begin{lemma}
\label{lem:action_Q}
We have that 
\begin{equation*}
\begin{aligned}
L_n \mathcal{Q}^n_t(h)
= \frac{\theta(n)}{2n^2} \mathcal{Q}^n_t(\Delta^nh
- 2n\alpha_n \mathscr{A}_nh)
+ \frac{2\theta(n)\alpha_n}{n^{3/2}} \mathcal{X}^n_t(\mathscr{D}_nh)
+ \frac{2\theta(n)}{n^2} \mathcal{Q}^n_t(\tilde{\mathscr{D}_n} h) 
+ E^n_t(h), 
\end{aligned}
\end{equation*}
where the error term satisfies the bound 
\begin{equation}
\label{eq:error_bound_LQ}
\begin{aligned}
&\mathbb{E}_n \bigg[\sup_{0\le t \le T} \bigg| \int_0^t E^n_s(h)ds \bigg|^2 \bigg] \\
&\quad\lesssim \theta(n)^2 \alpha_n^2  c_{k_*}^2\beta_n^{2(k_*-2)} 
\bigg(
\frac{T^2}{n^2}
\|\mathscr{A}_nh\|^2_{2,n} 
+ \frac{T^2}{n^2}
\|\mathscr{E}_nh\|^2_{2,n}
+ \frac{T}{\theta(n)}
\|h\|^2_{2,n} 
\bigg) .
\end{aligned}
\end{equation}
Here \Add{recall that} $k_*\in \{3,4,\ldots\}$ denotes the smallest finite number such that $c_{k_*}=V^{(k_*)}(0)\neq 0$.
Moreover, the above error bound is improved when $k_*=3$ as follows.
For any $\ell\in\mathbb{N}$ we have that 
\begin{equation}
\label{eq:error_bound_LQ_k3case}
\begin{aligned}
&\mathbb{E}_n \bigg[\sup_{0\le t \le T} \bigg| \int_0^t E^n_s(h)ds  \bigg|^2 \bigg] \\
&\quad\lesssim \theta(n)^2 \alpha_n^2  c_{3}^2\beta_n^{2} 
\bigg[ \frac{T^2}{n^2}
\| \mathscr{A}_nh \|^2_{2,n} 
+ \bigg(\frac{T^2\ell}{\theta(n)n^2} + \frac{T}{n^2\ell}\bigg)
\| \mathscr{E}_nh \|^2_{2,n}
+ \frac{T}{\theta(n)n^2} 
\sum_{j\in\mathbb{Z}} h\bigg( \frac{j}{n},\frac{j+1}{n}\bigg)^2 \bigg] .
\end{aligned}
\end{equation}
\end{lemma}
\begin{proof}
Here we derive the principal part of the action of the generator on the quadratic field $\mathcal{Q}^n_\cdot$. 
The error bounds \eqref{eq:error_bound_LQ_k3case} and \eqref{eq:error_bound_LQ} will be given in Lemma~\ref{lem:reminder_estimate_k3case} and \ref{lem:reminder_estimate}, respectively.
Let us begin with a computation of the symmetric part.
Note that
\begin{equation*}
\begin{aligned}
2S(\overline{\xi}_j \overline{\xi}_{j^\prime})
&= \overline{\xi}_{j^\prime} \Delta \xi_j
+ \overline{\xi}_{j} \Delta\xi_{j^\prime} 
- (\xi_{j+1}-\xi_{j})^2
\mathbf{1}_{j^\prime=j+1}
- (\xi_{j}-\xi_{j-1})^2 
\mathbf{1}_{j^\prime=j-1}\\
&\quad+ [(\xi_{j+1}-\xi_{j})^2 + (\xi_j-\xi_{j-1})^2]  
\mathbf{1}_{j^\prime=j} ,
\end{aligned}
\end{equation*}
where $\Delta$ denotes the discrete Laplacian defined by $\Delta g_j= g_{j+1} + g_{j-1}-2g_j$ for each $(g_j)_{j}$. 
Moreover, we use the short-hand notation $h_{j,j^\prime}=h(j/n,j^\prime/n)$ in the proof.
Then we compute  
\begin{equation*}
\begin{aligned}
S\mathcal{Q}^n_t(h)
&= \frac{1}{2n} \sum_{j\neq j^\prime} 
(\overline{\xi}_{j^\prime} \Delta \xi_j
+ \overline{\xi}_{j} \Delta\xi_{j^\prime}) h_{j,j^\prime} 
- \frac{1}{2n} \sum_{j\in\mathbb{Z}}
(\xi_{j+1}-\xi_{j})^2 (h_{j,j+1}+h_{j+1,j})\\
&= \frac{1}{2n^3} \sum_{j,j^\prime} 
\overline{\xi}_{j} \overline{\xi}_{j^\prime} 
(\Delta^n h)_{j,j^\prime} 
- \frac{1}{n} \sum_{j\in\mathbb{Z}} \overline{\xi}_j \Delta \xi_j h_{j,j}
- \frac{1}{2n} \sum_{j\in\mathbb{Z}} 
(\xi_{j+1}-\xi_{j})^2 (h_{j,j+1}+h_{j+1,j})\\
&= \frac{1}{2n^2} \mathcal{Q}^n_t (\Delta^n h)
+ \frac{1}{2n^3} \sum_{j\in\mathbb{Z}} 
(\overline{\xi}_j)^2 (\Delta^n h)_{j,j}
- \frac{1}{n} \sum_{j\in\mathbb{Z}} \overline{\xi}_j \Delta \xi_j h_{j,j}\\
&\quad- \frac{1}{2n} \sum_{j\in\mathbb{Z}} 
(\xi_{j+1}-\xi_{j})^2 (h_{j,j+1}+h_{j+1,j})\\
&= \frac{1}{2n^2} \mathcal{Q}^n_t(\Delta^nh)
+ \frac{1}{n} \sum_{j\in\mathbb{Z}} \overline{\xi}_j \overline{\xi}_{j+1}
(h_{j,j+1}+h_{j+1,j} -h_{j,j}-h_{j+1,j+1}). 
\end{aligned}
\end{equation*}
On the other hand, we have for the anti-symmetric part that  
\begin{equation*}
A(\overline{\xi}_j \overline{\xi}_{j^\prime})
= (\xi_{j-1} - \xi_{j+1}) \xi^\prime_j \overline{\xi}_{j^\prime}
+ (\xi_{j^\prime-1} - \xi_{j^\prime+1}) \xi^\prime_{j^\prime} 
\overline{\xi}_{j}.  
\end{equation*}
Here, note that $\xi^\prime_j 
= V_{\beta_n}^{\prime\prime}(\eta_j) 
= 1+ \frac{c_{k_*}}{(k_*-2)!}(\beta_n\xi_j)^{k_*-2} 
+ O(\beta_n^{k_*-1})$. 
Accordingly, we have that 
\begin{equation*}
A\mathcal{Q}^n_t(h) 
= I_n + \frac{c_{k_*}}{(k_*-2)!}\beta_n^{k_*-2} J_n
+ O(\beta_n^{k_*-1}), 
\end{equation*} 
where 
\begin{equation*}
I_n = \frac{1}{n} \sum_{j\neq j^\prime}
[(\xi_{j-1}-\xi_{j+1})\overline{\xi}_{j^\prime}
+ (\xi_{j^\prime-1}-\xi_{j^\prime+1}) \overline{\xi}_j]
h_{j,j^\prime},
\end{equation*}
and 
\begin{equation*}
J_n = \frac{1}{n} \sum_{j\neq j^\prime}
[(\xi_{j-1}-\xi_{j+1})\Add{g_*}(\xi_j) \overline{\xi}_{j^\prime}
+ (\xi_{j^\prime-1}-\xi_{j^\prime+1})
\Add{g_*}(\xi_{j^\prime}) \overline{\xi}_j] h_{j,j^\prime} .
\end{equation*}
\Add{
Above, we defined the function $g_*$ by 
\begin{equation}
\label{eq:g_specific_definition}
g_*(x)=x^{k_*-2}.
\end{equation}
}
Then, we have  
\begin{equation}
\label{eq:I}
\begin{aligned}
I_n
&= \frac{1}{n} \sum_{j,j^\prime}
[(\xi_{j-1}-\xi_{j+1}) \overline{\xi}_{j^\prime}
+ (\xi_{j^\prime-1} - \xi_{j^\prime+1}) \overline{\xi}_{j}]h_{j,j^\prime}
- \frac{2}{n} \sum_{j\in\mathbb{Z}} 
(\xi_{j-1}-\xi_{j+1})\overline{\xi}_j h_{j,j} \\
&= \frac{1}{n^2}\sum_{j,j^\prime} \overline{\xi}_{j} \overline{\xi}_{j^\prime} (\mathscr{A}_n h)_{j,j^\prime}
- \frac{2}{n^2} \sum_{j\in\mathbb{Z}} \overline{\xi}_j \overline{\xi}_{j+1} 
(\tilde{\mathscr{F}}_n h)_j \\
&= -\frac{1}{n}\mathcal{Q}^n_t(\mathscr{A}_nh) 
- \frac{1}{n^2} \sum_{j\in\mathbb{Z}} (\overline{\xi}_{j})^2 (\mathscr{A}_nh)_{j,j} 
- \frac{2}{n^2} \sum_{j\in\mathbb{Z}} \overline{\xi}_j \overline{\xi}_{j+1} 
(\tilde{\mathscr{F}}_n h)_j .
\end{aligned}
\end{equation}
Additionally, note that 
\begin{equation*}
\begin{aligned}
(\overline{\xi}_j)^2
&= \bigg( \eta_j -\lambda + \frac{c_{k_*}}{(k_*-1)!} \beta_n^{k_*-2} \eta_j^{k_*-1} + O(\beta_n^{k_*-1}) \bigg)^2 \\
&= \lambda^2 + \eta_j^2 - 2\lambda \eta_j 
+ \frac{2c_{k_*}}{(k_*-1)!} \beta_n^{k_*-2} \eta_j^{k_*-1} (\eta_j-\lambda) + O(\beta_n^{k_*-1}) \\
&=\lambda^2 + 2(\zeta_j - \lambda \eta_j) 
+ \frac{2c_{k_*}}{k_*!} \beta_n^{k_*-2}
\eta_j^{k_*-1} \big[(k_*-1)\eta_j-k_*\lambda \big] 
+ O(\beta_n^{k_*-1}) .
\end{aligned}
\end{equation*}
In particular, the second term in the utmost right-hand side of \eqref{eq:I} gives $2\theta(n)\alpha_nn^{-3/2}\mathcal{X}^n_t(\mathscr{D}_nh)$ with an error term which comes from the third term of the last display, which is substituted back into the second term in the last line of \eqref{eq:I}.
This error term is treated with the Schwarz inequality in the following way.
\begin{equation*}
\begin{aligned}
& \mathbb{E}_n \bigg[ \sup_{0\le t \le T} \bigg| 
\frac{1}{n^2} \int_0^t \sum_{j\in \mathbb{Z}} 
F(\eta_j(s)) (\mathscr{A}_n h)_{j,j}(s) ds\bigg|^2 \bigg]\\
&\quad\le \frac{T}{n^4} \int_0^T \mathbb{E}_n 
\bigg[ \bigg( \sum_{j\in \mathbb{Z}} F(\eta_j(s)) (\mathscr{A}_nh)_{j,j}(s) \bigg)^2 \bigg] ds 
\lesssim \frac{T^2}{n^2} \| \mathscr{A}_nh\|^2_{L^2(\mathbb{R}^2)} ,
\end{aligned}
\end{equation*}
where $F(\eta)=\eta^{k_*-1}(\eta-\lambda)$ and in the second estimate we used the fact that for any symmetric function $h$ it holds that $\sum_{j} (\mathscr{A}_nh)_{j,j}=0$.  
This bound gives an error term which is proportional to $\| \mathscr{A}_nh\|^2_{2,n}$ in \eqref{eq:error_bound_LQ}.  
Finally, the estimates \eqref{eq:error_bound_LQ_k3case} and \eqref{eq:error_bound_LQ} for $J_n$ follow from \Add{Lemma~\ref{lem:reminder_estimate} with $g=g_*$ and Lemma~\ref{lem:reminder_estimate_k3case} below}. 
\end{proof}

In the sequel, we give a quantitative estimate for the error term $J_n$ which appears in the proof of Lemma~\ref{lem:action_Q}. 
For that purpose, we make use of the following Kipnis-Varadhan estimate~\cite{kipnis1986central, chang2001equilibrium}, for each smooth mean-zero function $F:[0,T]\times \Omega \to\mathbb R$
\begin{equation*}
\begin{aligned}
\mathbb{E}_n \bigg[\sup_{0\le t \le T} \bigg|
\int_0^t F(s,\eta(s)) ds \bigg|^2 \bigg]
\lesssim \int_0^T \| F(s,\cdot)\|^2_{-1,n} ds,
\end{aligned}
\end{equation*}
where $\| \cdot \|_{-1,n}$-norm is defined through the variational formula
\begin{equation*}
\|F\|^2_{-1,n}
=\sup_{f} \bigg\{ 2\langle F,f\rangle_{L^2(\nu_n)}
- \langle f, -\theta(n) Sf\rangle_{L^2(\nu_n)} \bigg\},
\end{equation*}
where the supremum is taken over all $L^2(\nu_n)$-local functions.

\begin{lemma}
\label{lem:reminder_estimate}
Let 
\begin{equation}\label{eq:for_Jn}
J_n^g 
= \frac{1}{n} \sum_{j\neq j^\prime}
[(\xi_{j-1}-\xi_{j+1}) g(\xi_j) \overline{\xi}_{j^\prime} 
+ (\xi_{j^\prime-1}-\xi_{j^\prime+1})
g(\xi_{j^\prime}) \overline{\xi}_j ]
h\bigg( \frac{j}{n},\frac{j^\prime}{n}\bigg),
\end{equation}
where $g:\mathbb{R}\to\mathbb{R}$ is any continuous function such that $E_{\nu_n}[g(\eta_j)^2]<+\infty$. 
Then, 
\begin{equation}
\label{eq:reminder_estimate}
\mathbb{E}_n \bigg[ \sup_{0\le t \le T} \bigg| \int_0^t J_n^g(s) ds \bigg|^2 \bigg]
\lesssim 
\frac{T^2}{n^2}
\| \mathscr{E}_nh\|^2_{2,n}
+ \frac{T}{\theta(n)} 
\|h\|^2_{2,n} . 
\end{equation}
\end{lemma}
\begin{proof}
In this proof, we write the short-hand notation $h_{j,j^\prime}=h(j/n,j^\prime/n)$ .
In the sequel, we may assume that $g(\xi_j)$ is centered, namely, that $E_{\nu_n}[g(\xi_j)]=0$.
Let $\Psi_j = \xi_j g(\xi_{j+1}) - \xi_{j+1} g(\xi_j)$. 
Since $h$ is symmetric, we have 
\begin{equation*}
\begin{aligned}
J_n^g
&= \frac{2}{n} \sum_{j\neq j^\prime}
(\xi_{j-1}-\xi_{j+1}) g(\xi_j) \overline{\xi}_{j^\prime} h_{j,j^\prime} \\
&= \frac{2}{n} \sum_{j\neq j^\prime}
\big( \xi_{j-1}g(\xi_j) - \xi_{j} g(\xi_{j+1}) \big) \overline{\xi}_{j^\prime} h_{j,j^\prime}
+ \frac{2}{n} \sum_{j\neq j^\prime}
\big( \xi_j g(\xi_{j+1}) - \xi_{j+1} g(\xi_j) \big) 
\overline{\xi}_{j^\prime} h_{j,j^\prime} \\
&= \frac{2}{n} \sum_{j,j^\prime}
\big( \xi_{j-1}g(\xi_j) - \xi_{j} g(\xi_{j+1}) \big) \overline{\xi}_{j^\prime} h_{j,j^\prime}
+ \frac{2}{n} \sum_{j^\prime\neq j,j+1}
\big( \xi_j g(\xi_{j+1}) - \xi_{j+1} g(\xi_j) \big) 
\overline{\xi}_{j^\prime} h_{j,j^\prime} \\
&\quad- \frac{2}{n} \sum_{j\in\mathbb{Z}} 
\big( \xi_{j-1}g(\xi_j)-\xi_jg(\xi_{j+1}) \big) \overline{\xi}_j h_{j,j} 
+ \frac{2}{n} \sum_{j\in\mathbb{Z}}
\big( \xi_j g(\xi_{j+1}) - \xi_{j+1} g(\xi_j) \big) \overline{\xi}_{j+1} h_{j,j+1} \\
&= \frac{2}{n^2} \sum_{j,j^\prime} 
\overline{\xi}_j g(\xi_{j+1}) \overline{\xi}_{j^\prime} 
(\mathscr{E}_n h)_{j,j^\prime} 
+ \frac{2}{n} \sum_{j^\prime\neq j,j+1}
\Psi_j \overline{\xi}_{j^\prime} h_{j,j^\prime} \\
&\quad- \frac{2}{n} \sum_{j\in\mathbb{Z}} 
\big( \xi_{j-1}g(\xi_j)-\xi_jg(\xi_{j+1}) \big) \overline{\xi}_j h_{j,j} 
+ \frac{2}{n} \sum_{j\in\mathbb{Z}}
\big( \xi_j g(\xi_{j+1}) - \xi_{j+1} g(\xi_j) \big) \overline{\xi}_{j+1} h_{j,j+1} \\
&\eqqcolon H_1 + H_2 + H_3, 
\end{aligned}
\end{equation*}
where we set 
\begin{equation*}
H_1
= \frac{2}{n^2} \sum_{j,j^\prime} 
\overline{\xi}_j g(\xi_{j+1}) \overline{\xi}_{j^\prime} 
(\mathscr{E}_n h)_{j,j^\prime} ,\quad 
H_2= 
\frac{2}{n} \sum_{j^\prime\neq j,j+1}
\Psi_j \overline{\xi}_{j^\prime} h_{j,j^\prime} ,
\end{equation*}
and 
\begin{equation*}
H_3 =
- \frac{2}{n} \sum_{j\in\mathbb{Z}} 
\big( \xi_{j-1}g(\xi_j)-\xi_jg(\xi_{j+1}) \big) \overline{\xi}_j h_{j,j} 
+ \frac{2}{n} \sum_{j\in\mathbb{Z}}
\big( \xi_j g(\xi_{j+1}) - \xi_{j+1} g(\xi_j) \big) \overline{\xi}_{j+1} h_{j,j+1} .
\end{equation*}

First, $H_1$ is estimated by a direct $L^2$-computation as follows. 
\begin{equation*}
\begin{aligned}
&\mathbb{E}_n 
\bigg[ \sup_{0\le t \le T} \bigg|
\int_0^t \frac{2}{n^2}
\sum_{j\neq j^\prime} 
\overline{\xi}_j g(\xi_{j+1}) \overline{\xi}_{j^\prime}
(\mathscr{E}_nh)_{j,j^\prime} ds\bigg|^2 
\bigg] \\
&\quad\lesssim \frac{T}{n^4} 
\int_0^T 
\mathbb{E}_n \bigg[ \bigg( 
\sum_{j\neq j^\prime} 
\bar\xi_j g(\xi_{j+1}) \overline{\xi}_{j^\prime}
(\mathscr{E}_nh)_{j,j^\prime} \bigg)^2 \bigg] ds
\lesssim \frac{T^2}{n^2} 
\| \mathscr{E}_nh \|^2_{2,n} . 
\end{aligned}
\end{equation*}

To give an estimate of $H_2$, we apply the Kipnis-Varadhan inequality for this term.
Note that $\Psi_j(\eta)$ satisfies $\Psi_j(\eta^{j,j+1})=-\Psi(\eta)$. 
We take an arbitrary local function $f:\Omega\to \mathbb{R}$. 
Then by Young's inequality, we have for any $f\in L^2(\nu_n)$ that
\begin{equation*}
\begin{aligned}
&\bigg\langle \frac{4}{n}\sum_{j\neq j^\prime}
\Psi_j(\eta) \overline{\xi}_{j^\prime}
h_{j,j^\prime}, f(\eta)
\bigg\rangle_{L^2(\nu_n)}\\
&\quad= \frac{2}{n} E_{\nu_n} \bigg[\sum_{j\in\mathbb{Z}} 
\Psi_j(\eta) \nabla_{j,j+1}f(\eta)
\sum_{j^\prime\in\mathbb{Z}} 
\overline{\xi}_{j^\prime} 
h_{j,j^\prime} 
\bigg] \\
&\quad\le \frac{A}{n} \sum_{j\in\mathbb{Z}}
E_{\nu_n} [(\nabla_{j,j+1}f(\eta))^2] 
+ \frac{1}{An} \sum_{j\in\mathbb{Z}}
E_{\nu_n} \bigg[ \Psi_j(\eta)^2
\bigg(\sum_{j^\prime\in\mathbb{Z}} \overline{\xi}_{j^\prime}
h_{j,j^\prime}\bigg)^2 \bigg],
\end{aligned}
\end{equation*}
for any $A>0$. 
In particular, when $A=n\theta(n)/4$ the first term in the utmost right-hand side is nothing but the Dirichlet form of the generator $\theta(n)S$. 
In addition, the other remainder term is estimated by \begin{equation*}
\begin{aligned}
\frac{4}{n^2\theta(n)} \sum_{j,j^\prime}
E_{\nu_n} [\Psi_j(\eta)^2 (\overline{\xi}_j)^2 h^2_{j,j^\prime} ]
\lesssim \frac{1}{n^2\theta(n)}
\sum_{j,j^\prime} h^2_{j,j^\prime}.   
\end{aligned}
\end{equation*}
We obtain the desired bound by the Kipnis-Varadhan inequality, i.e., 
\begin{equation*}
\begin{aligned}
\mathbb{E}_n \bigg[\sup_{0\le t \le T}
\bigg| \int_0^t \frac{2}{n} \sum_{j^\prime \neq j,j+1}
\Psi_j(s) \overline{\xi}_{j^\prime}(s) h_{j,j^\prime} ds \bigg|^2 \bigg] 
\lesssim \frac{T^2}{\theta(n)} \| h\|^2_{2,n} . 
\end{aligned}
\end{equation*}

Finally, for $H_3$, we compute 
\begin{equation*}
\begin{aligned}
H_3
&= - \frac{2}{n} \sum_{j\in\mathbb{Z}} 
\overline{\xi}_{j}\overline{\xi}_{j+1} g(\xi_{j+1})
(h_{j+1,j+1}-h_{j,j+1})
+ \frac{2}{n} \sum_{j\in\mathbb{Z}}
\big((\overline{\xi}_j)^2 g(\xi_{j+1}) h_{j,j}
- (\overline{\xi}_{j+1})^2 g(\xi_j) h_{j,j+1}\big) \\
&= - \frac{2}{n^2} \sum_{j\in\mathbb{Z}} 
\overline{\xi}_{j-1}\overline{\xi}_j g(\xi_j)
\big( (\tilde{\mathscr{F}}_nh)_{j} -(\Tilde{\mathscr{E}}_n h)_j \big)
- \frac{2}{n^2} \sum_{j\in\mathbb{Z}}
(\overline{\xi}_j)^2 g(\xi_{j+1}) 
(\Tilde{\mathscr{E}}_n h)_j \\
&\quad+ \frac{2}{n} \sum_{j\in\mathbb{Z}}
\big( (\overline{\xi}_j)^2 g(\xi_{j+1}) -
(\overline{\xi}_{j+1})^2 g(\xi_j) \big) h_{j,j+1}.
\end{aligned}
\end{equation*}
This is bounded as follows: 
\begin{equation*}
\begin{aligned}
& \mathbb{E}_n \bigg[ \sup_{0\le t \le T} \bigg|\int_0^t  H_3(s) ds \bigg|^2 \bigg]\\
&\quad\lesssim
\frac{T^2}{n^4} \bigg(\sum_{j\in\mathbb{Z}}  (\Tilde{\mathscr{F}}_n h)_j^2
+ \sum_{j\in\mathbb{Z}} (\Tilde{\mathscr{E}}_n h)_j^2 \bigg) 
+ \frac{T}{\theta(n)n^2} 
\bigg( \sum_{j\in\mathbb{Z}} (\Tilde{\mathscr{E}}_n h)_j^2
+ \sum_{j\in\mathbb{Z}} (h_{j,j})^2 \bigg) 
\\
&\quad\le \frac{T^2}{n^4} \bigg(\sum_{j,j^\prime\in\mathbb{Z}}  (\mathscr{F}_n h)_{j,j^\prime}^2
+ \sum_{j,j^\prime\in\mathbb{Z}} (\mathscr{E}_n h)_{j,j^\prime}^2 \bigg) 
+ \frac{T}{\theta(n)n^2} 
\bigg( \sum_{j,j^\prime\in\mathbb{Z}} (\mathscr{E}_n h)_{j,j^\prime}^2
+ \sum_{j,j^\prime\in\mathbb{Z}} h_{j,j^\prime}^2 \bigg) 
\\
&\quad\lesssim \frac{T^2}{n^2} \bigg( \| \mathscr{F}_nh \|^2_{L^2(\mathbb{R}^2)} + \| \mathscr{E}_n h\|^2_{L^2(\mathbb{R}^2)}\bigg) 
+ \frac{T}{\theta(n)} \bigg( \| \mathscr{E}_n h\|^2_{L^2(\mathbb{R}^2)} + \| h \|^2_{L^2(\mathbb{R}^2)}  \bigg) .
\end{aligned}
\end{equation*}
These bounds can be absorbed in those of $H_1$ and $H_2$, noting that  $(\mathscr{F}_nh)_{j,j^\prime} = (\mathscr{E}_nh)_{j,j^\prime+1} - (\mathscr{E}_nh)_{j,j^\prime}$. 
Hence we complete the proof. 
\end{proof}

On the other hand, for the case when $k_*=3$, we have the following improved estimate, which can be obtained by the fact that \Add{$J_n^{g_*}$ is written in a gradient form when $k_*=3$ where $g_*$ is the function defined in \eqref{eq:g_specific_definition}. 
}

\begin{lemma}[The case $k_*=3$.]
\label{lem:reminder_estimate_k3case}
\Erase{Recall \eqref{eq:for_Jn} and take $g(x)=x^{k_*-2}=x$, so that}
\Add{Let}
\begin{equation*}
\Add{\tilde{J}_n} = \frac{1}{n} \sum_{j\neq j^\prime}
[(\bar\xi_{j-1}\bar\xi_j-\bar\xi_j\bar\xi_{j+1}) \overline{\xi}_{j^\prime}
+ (\bar\xi_{j^\prime-1}\bar\xi_{j^\prime}-\bar\xi_{j^\prime}\bar\xi_{j^\prime+1})\overline{\xi}_j] 
h\bigg( \frac{j}{n},\frac{j^\prime}{n}\bigg).
\end{equation*}
Then, for any $\ell\in \mathbb{N}$, we have that 
\begin{equation*}
\begin{aligned}
\mathbb{E}_n \bigg[\sup_{0\le t\le T}\bigg|\int_0^t \Add{\tilde J_n}(s)ds \bigg|^2\bigg] 
\lesssim \bigg( \frac{T^2\ell^2}{n^2\theta(n)} + \frac{T}{n^2\ell} \bigg) 
\| \mathscr{E}_nh\|^2_{2,n} 
+ \frac{T}{\theta(n)n^2} 
\sum_{j\in\mathbb{Z}} 
h\bigg( \frac{j}{n},\frac{j+1}{n}\bigg)^2.
\end{aligned}
\end{equation*}
\end{lemma}
\begin{proof}
Here we simply write $h_{j,j^\prime}=h(j/n,j^\prime/n)$, as in the proof of Lemma~\ref{lem:action_Q}.
Since $h$ is symmetric we have that
\begin{equation*}
\begin{aligned}
\Add{\tilde J_n} &= \frac{2}{n}\sum_{j\neq j^\prime}
(\overline{\xi}_{j-1}\overline{\xi}_{j} 
- \overline{\xi}_{j}\overline{\xi}_{j+1})
\overline{\xi}_{j^\prime}h_{j,j^\prime} \\
&= \frac{2}{n}\sum_{j,j^\prime}
(\overline{\xi}_{j-1}\overline{\xi}_{j} 
- \overline{\xi}_{j}\overline{\xi}_{j+1})
\overline{\xi}_{j^\prime}h_{j,j^\prime} 
-\frac{2}{n}\sum_{j} (\overline{\xi}_j)^2 (\overline{\xi}_{j-1}-\overline{\xi}_{j+1}) h_{j,j}\\
&= \frac{2}{n^2}\sum_{j,j^\prime}
\overline{\xi}_{j} \overline{\xi}_{j+1} \overline{\xi}_{j^\prime}
(\mathscr{E}_nh)_{j,j^\prime} 
- \frac{1}{n}\sum_j (\overline{\xi}_j)^2 (\overline{\xi}_{j-1}-\overline{\xi}_{j+1}) h_{j,j}\\
&= \frac{2}{n^2}\sum_{j^\prime\neq j,j+1}
\overline{\xi}_{j} \overline{\xi}_{j+1} \overline{\xi}_{j^\prime}
(\mathscr{E}_nh)_{j,j^\prime}
- \frac{2}{n}\sum_j (\overline{\xi}_j)^2 (\overline{\xi}_{j-1}-\overline{\xi}_{j+1}) h_{j,j}\\
&\quad+ \frac{2}{n}\sum_{j\in\mathbb{Z}} 
(\overline{\xi}_j)^2 \overline{\xi}_{j+1} (h_{j+1,j}-h_{j,j})
+ \frac{2}{n}\sum_{j\in\mathbb{Z}}
(\overline{\xi}_{j-1})^2 \overline{\xi}_j 
(h_{j+1,j+1}-h_{j,j+1}) \\
&= \frac{2}{n^2}\sum_{j^\prime\neq j,j+1}
\overline{\xi}_{j} \overline{\xi}_{j+1} \overline{\xi}_{j^\prime}
(\mathscr{E}_nh)_{j,j^\prime}
+ \frac{2}{n} \sum_{j\in\mathbb{Z}}
\big((\overline{\xi}_j)^2 \overline{\xi}_{j+1}
- (\overline{\xi}_{j+1})^2 \overline{\xi}_j \big) h_{j,j+1}
\eqqcolon K_1 + K_2. 
\end{aligned}
\end{equation*}

To estimate $K_1$ we do the following. 
Fix $\ell\in\mathbb{N}$.  
At first instance we split the sum in two cases. First either $j^\prime\in R_{j+1}^\ell:=\{j+1,\cdots, j+1+\ell\} $ or  $j^\prime\in L_{j}^\ell:=\{j-\ell,\cdots, j-1\}$ or $j^\prime\notin {R_{j+1}^\ell\cup L_{j}^\ell}$.  All the cases are treated analogously. In the first case  we make the replacement $\xi_j$ by $\overleftarrow{\xi}_j^\ell$, in the second we replace $\xi_{j+1}$ by $\overrightarrow{\xi}_{j+1}^\ell$ and in the third case we do one of them (does not matter which). Note that by the one-block estimate, see~\cite{gonccalves2014nonlinear} for example, we have that  
\begin{equation}
\label{eq:multiscale_estimate_replacement}
\begin{aligned}
&\mathbb{E}_n 
\bigg[ \sup_{0\le t \le T} \bigg|
\int_0^t \frac{1}{n^2}
\sum_{j\in\mathbb Z}\sum_{j^\prime\notin{L_j^\ell}} 
\Big(\overline{\xi_j}(s)-\overleftarrow{\xi}^\ell_j(s)\Big) \bar\xi_{j+1}(s) \overline{\xi}_{j^\prime}(s)
(\mathscr{E}_nh)_{j,j^\prime} ds\bigg|^2 
\bigg] \leq \frac{T\ell^2}{n^2\theta(n)}\|\mathcal E_nh\|^2_{2,n}.
\end{aligned}
\end{equation}
Indeed, to obtain the last estimate, we proceed as follows. Set
\begin{equation*}
\sum_{j\in\mathbb Z} \sum_{j^\prime\notin{L_j^\ell}}
(\overline{\xi_j} - \overleftarrow{\xi}^\ell_j)
\overline{\xi}_{j+1} \overline{\xi}_{j^\prime}
(\mathcal{E}_nh)_{j,j^\prime}
= \sum_{j\in\mathbb Z} (\xi_j - \overleftarrow{\xi}^\ell_j)
\overline{\xi}_{j+1} \Phi_j,
\end{equation*}
where $\Phi_j = \sum_{j^\prime\notin{L_j^\ell}} \overline{\xi}_{j^\prime} (\mathcal{E}_nh)_{j,j^\prime}$. 
We have that 
\begin{equation*}
\begin{aligned}
 \sum_{j\in\mathbb{Z}} (\xi_j - \overleftarrow{\xi}^\ell_j)
\overline{\xi}_{j+1} \Phi_j 
= \sum_{j\in\mathbb{Z}} \sum_{i=0}^{\ell-2}
\overline{\xi}_{j+1} (\xi_{j-i}-\xi_{j-i-1}) \psi_i \Phi_j 
= \sum_{k\in\mathbb{Z}} F_k (\xi_k -\xi_{k-1}),
\end{aligned}
\end{equation*}
where $F_k=\sum_{i=0,\ldots,\ell-2} \overline{\xi}_{k+i+1} \Phi_{k+i} \psi_{i}$ and $\psi_i = (\ell-i-1)/\ell$. 
Note that the functional $F_k(\eta)$ is invariant under the action $\eta 
\mapsto \eta^{k,k-1}$. 
Therefore, for any $f\in L^2(\nu_n)$, it holds 
\begin{equation*}
\begin{aligned}
\bigg\langle 2  \sum_{j\in\mathbb{Z}} (\xi_j - \overleftarrow{\xi}^\ell_j)
\overline{\xi}_{j+1} \Phi_j , f(\eta) \bigg\rangle_{L^2(\nu_n)} 
&= 2 \sum_{j\in\mathbb{Z}} E_{\nu_n} [F_j (\xi_j-\xi_{j-1})f] \\
&= 2 \sum_{j\in\mathbb{Z}}
E_{\nu_n}[F_{j+1}\overline{\xi}_j (\nabla_{j,j+1}f(\eta))] .
\end{aligned}
\end{equation*}
Then by Young's inequality for any $A>0$ the above display is bounded by
\begin{equation*}
\begin{aligned}
A \sum_{j\in\mathbb{Z}} E_{\nu_n}
[(\nabla_{j,j+1}f)^2]
+ \frac{1}{A} \sum_{j\in\mathbb{Z}} 
E_{\nu_n} [F_{j+1}^2
(\overline{\xi}_j)^2] .
\end{aligned}
\end{equation*}
In particular, when we take $A=\theta(n)/4$, the first term is nothing but the Dirichlet form. 
For the reminder term, note 
\begin{equation*}
E_{\nu_n}[F_j^2]
= E_{\nu_n}
\bigg[ \bigg( \sum_{i=0}^{\ell-2}
\overline{\xi}_{j+i+1}
\psi_{i}
\sum_{j^\prime>j}
\overline{\xi}_{j^\prime+i}
(\mathcal{E}_nh)_{j^\prime\notin{L_j^\ell}} 
\bigg)^2 \bigg] 
\lesssim  \ell^2 \sum_{j,j^\prime} 
(\mathcal{E}_nh)_{j,j^\prime}^2 ,
\end{equation*}
which follows from a crude $L^2$-estimate and the fact that $\psi_i$ is bounded by a constant. 
Hence the desired estimate follows from the Kipnis-Varadhan inequality.

Now observe that from the Schwarz inequality we have that
\begin{equation}
\label{eq:multiscale_estimate_averageitself}
\begin{aligned}
&\mathbb{E}_n 
\bigg[ \sup_{0\le t \le T} \bigg|
\int_0^t \frac{1}{n^2}
\sum_{j,j^\prime\notin{L_j^\ell}} 
\overleftarrow{\xi}^\ell_j(s) 
\overline{\xi}_{j+1}(s) 
\overline{\xi}_{j^\prime}(s)
(\mathscr{E}_nh)_{j,j^\prime} ds\bigg|^2 
\bigg] 
\le\frac{T^2}{n^2\ell}
\|\mathcal E_nh\|^2_{2,n}.
\end{aligned}\end{equation}
The proof is a consequence of the fact that  
\begin{equation*}
\begin{aligned}
&E_{\nu_n} \bigg[ \bigg( \frac{1}{n^2}
\sum_{j,j^\prime\notin{L_j^\ell}} 
\overleftarrow{\xi}^\ell_j 
\overline{\xi}_{j+1}
\overline{\xi}_{j^\prime}
(\mathscr{E}_nh)_{j,j^\prime} \bigg)^2 
\bigg]
=\frac{1}{n^4}\sum_{j,j^\prime\notin{L_j^\ell}} 
E_{\nu_n}[(\overleftarrow{\xi}^\ell_j)^2 (\bar\xi_{j+1})^2 (\overline{\xi}_{j^\prime})^2]
(\mathscr{E}_nh)_{j,j^\prime}^2,
\end{aligned}
\end{equation*}
which follows from $\mathbb{E}_n[\overline{\xi}_j \overline{\xi}_{j^\prime}]=0$ provided $j\neq j^\prime$, and the fact that 
$\mathbb{E}_n [(\overleftarrow{\xi}^\ell_j)^2 ] \le \ell^{-1}$.
Putting the two last estimates \eqref{eq:multiscale_estimate_replacement} and \eqref{eq:multiscale_estimate_averageitself} together we get that
\begin{equation*}
\begin{aligned}
\mathbb{E}_n 
\bigg[ \sup_{0\le t \le T} \bigg|
\int_0^t \frac{1}{n^2}
\sum_{j,j^\prime\neq j,j+1} 
\overline{\xi_j}(s) 
\overline{\xi}_{j+1}(s) 
\overline{\xi}_{j^\prime}(s)
(\mathscr{E}_nh)_{j,j^\prime} ds\bigg|^2  
\lesssim \bigg( \frac{T^2\ell^2}{n^2\theta(n)} + \frac{T}{n^2\ell} \bigg)  \| \mathscr{E}_nh \|^2_{2,n}. 
\end{aligned}
\end{equation*}
Finally, the remainder term $K_2$ can be estimated by the Kipnis-Varadhan inequality, which yields
\begin{equation*}
\mathbb{E}_n \bigg[ \sup_{0\le t\le T} \bigg| 
\int_0^t \frac{1}{n} \sum_{j\in\mathbb{Z}}
\big((\overline{\xi}_j)^2 \overline{\xi}_{j+1}
- (\overline{\xi}_{j+1})^2 \overline{\xi}_j \big)(s) 
h_{j,j+1} (t) dt \bigg|^2
\bigg]
\lesssim \frac{T}{\theta(n)n^2} \sum_{j\in\mathbb{Z}} h_{j,j+1}^2.  
\end{equation*}
This is verified exactly in the same way as the estimate of $H_1$ in the proof of Lemma~\ref{lem:reminder_estimate}.

To finish the proof of Lemma~\ref{lem:reminder_estimate_k3case}, we optimize over $\ell$, with $\ell=\sqrt n$ and we are done. 
\end{proof}

\subsection{Proof of Theorem~\ref{thm:3/2L_derivation}}
As a consequence of Lemma \ref{lem:action_Q}, we have that 
\begin{equation}
\label{eq:martingale_decomposition_Q}
\begin{aligned}
\frac{d}{dt} \mathscr{Q}^n_t(h) 
&=  \frac{1}{2}\mathbb{E}_n
\big[ \mathcal{X}^n_0(g) L_n \mathcal{Q}^n_t (h) \big] \\
&= \frac{\theta(n)}{2n^2} \mathscr{Q}^n_t(\Delta^nh - 2n\alpha_n \mathscr{A}_nh)
+ \frac{2\theta(n)\alpha_n}{n^{3/2}} \mathscr{S}^n_t(\mathscr{D}_nh) 
+ \frac{2\theta(n)}{n^2} \mathscr{Q}^n_t(\tilde{\mathscr{D}_n} h)  
+ E^n_t(h) ,
\end{aligned}
\end{equation}
where $E_n$ satisfies the bound \eqref{eq:error_bound_LQ}. 
Here recall the martingale decomposition \eqref{eq:martingale_decomposition_original}.
To relate \eqref{eq:martingale_decomposition_original} with \eqref{eq:martingale_decomposition_Q}, we solve the following Poisson equation.\\
\begin{equation}
\label{eq:Poisson}
\frac{1}{2} \Delta^n h
\bigg(\frac{j}{n},\frac{j^\prime}{n} \bigg)
- n\alpha_n \mathscr{A}_n h
\bigg(\frac{j}{n},\frac{j^\prime}{n} \bigg)
= 2n^{1/2} \alpha_n 
(\nabla^{n,1} \varphi \otimes \delta)
\bigg(\frac{j}{n},\frac{j^\prime}{n} \bigg) . 
\end{equation}
Hereafter we denote the solution of \eqref{eq:Poisson} by $h_n$.
Then, we have the following quantitative estimate for $h_n$ and its derivatives.

\begin{lemma}
\label{lem:poisson_sol_bound}
Let $h_n$ be the solution of the discrete Poisson equation \eqref{eq:Poisson}. 
Then, 
\begin{equation*}
\| h_n \|^2_{2,n} \lesssim n^{-1/2},\quad 
\| \mathscr{E}_n h \|^2_{2,n} \lesssim 1 , \quad 
\| \mathscr{A}_n h \|^2_{2,n} \lesssim n^{-1/2} .
\end{equation*}
\end{lemma}

As a consequence, the utmost right-hand side of the bound~\eqref{eq:error_bound_LQ} vanishes provided 
\begin{equation*}
\lim_{n\to\infty}
\theta(n)^2 \alpha_n^2
\beta_n^{2(k_*-2)} 
\bigg( 
\frac{1}{n^2} 
\| \mathscr{A}_nh_n\|^2_{2,n}
+ \frac{1}{n^2} 
\| \mathscr{E}_nh_n\|^2_{2,n}
+ \frac{1}{\theta(n)} 
\| h_n\|^2_{2,n}
\bigg)
=0.
\end{equation*}
According to Lemma~\ref{lem:poisson_sol_bound}, this is satisfied when $\beta_n=O(n^{-b})$ with $b\ge 1/(2k_*-4)$. 
Moreover, note that the error term in the decomposition~\eqref{eq:martingale_decomposition_original} vanishes under the condition of Theorem~\ref{thm:3/2L_derivation}.
Thus we hereafter write the error term simply by $o_n(1)$, which is always negligible when $n\to\infty$.
On the other hand, when $k_*=3$, we see that the critical exponent for $\beta_n$ improves by Lemma \ref{lem:reminder_estimate_k3case} and the following estimate.

\begin{lemma}
\label{lem:est_fourier_1d}
Let $h_n$ be the solution of the discrete Poisson equation \eqref{eq:Poisson}. 
Then 
\begin{equation*}
\frac{1}{n}\sum_{j\in\mathbb{Z}} 
h_n\bigg( \frac{j}{n},\frac{j+1}{n}\bigg)^2
\lesssim 1.
\end{equation*}
\end{lemma}

The proof of Lemmas \ref{lem:poisson_sol_bound} and \ref{lem:est_fourier_1d} is given in Section \ref{subsec:poisson_bound}. 
Now we combine the two identities \eqref{eq:martingale_decomposition_original} and \eqref{eq:martingale_decomposition_Q} to obtain
\begin{equation*}
\begin{aligned}
\frac{d}{dt} \mathscr{S}^n_t(f) 
&= -\frac{d}{dt} \mathscr{Q}^n_t(h)
+ \mathscr{S}^n_t \bigg( \frac{\theta(n)}{2n^2}\Delta^nh_n 
+\frac{2\theta(n)\alpha_n}{n^{3/2}}\mathscr{D}_nh_n
\bigg) 
+ \frac{2\theta(n)}{n^2}\mathscr{Q}^n_t(\tilde{\mathscr{D}_n}h_n) 
+ o_n(1) .
\end{aligned}
\end{equation*}
By integrating the last display in time we obtain 
\begin{equation}
\label{eq:martingale_decomposition_poisson}
\begin{aligned}
\mathscr{S}^n_T(f)-\mathscr{S}^n_0(f)
&= \int_0^T \mathscr{S}^n_t \bigg( \frac{\theta(n)}{2n^2}\Delta^nh_n 
+\frac{2\theta(n)\alpha_n}{n^{3/2}}\mathscr{D}_nh_n
\bigg) dt \\
&\quad+ \mathscr{Q}^n_0(h_n) -\mathscr{Q}^n_T(h_n)
+ \frac{2\theta(n)}{n^2}\int_0^T \mathscr{Q}^n_t(\tilde{\mathscr{D}}_nh_n)dt
+ o_n(1) .
\end{aligned}
\end{equation}
Since $\| h_n\|^2_{L^2(\mathbb{R}^2)}\lesssim n^{-1/2}$ by Lemma~\ref{lem:poisson_sol_bound}, the second and third terms in the utmost right-hand side of \eqref{eq:martingale_decomposition_poisson} vanish as $n\to+\infty$. 
Moreover, the contribution of the first term in the utmost right-hand side of \eqref{eq:martingale_decomposition_poisson} is captured by the following result~\cite[Lemma~3.2]{bernardin2018weakly}.

\begin{proposition}
\label{prop:levy_operator}
If $a=\min\{ 3/2+3\kappa/2 , 2 \}$ and $\kappa \in (0,\infty)$, then 
\begin{equation*}
\lim_{n\to\infty} 
\frac{1}{n} \sum_{j\in\mathbb{Z}} \bigg|
\big( n^{a-2} \Delta^n f 
- 2\gamma n^{a-\kappa-3/2} \mathscr{D}_n h_n\big) \bigg(\frac{j}{n}\bigg)
- \mathbb{L}_{\gamma, \kappa} f\bigg(\frac{j}{n}\bigg) \bigg|^2 = 0,
\end{equation*}
where $\mathbb{L}_{\gamma,\kappa}$ is the operator defined by \eqref{eq:levy_operator}. 
\end{proposition}

We apply Proposition~\ref{prop:levy_operator} by setting $\theta(n) = n^a$ and $\alpha_n = \gamma n^{-\kappa}$ to see that the first term on the utmost right-hand side of \eqref{eq:martingale_decomposition_poisson} gives rise to the L\'{e}vy operator $\mathbb{L}_{\gamma,\kappa}$. 
Finally, the fourth term in the utmost right-hand side of \eqref{eq:martingale_decomposition_poisson} is controlled with the following result.

\begin{lemma}
Let $h_n$ be the solution of the Poisson equation given in \eqref{eq:Poisson}, $\theta(n)=n^a$ and $\alpha_n= \gamma n^{-\kappa}$ with $a=\min\{ \frac{3}{2}(1+\kappa), 2\}$ and $\kappa\in (0,1)$.
For any $T>0$, we have that
\begin{equation*}
\lim_{n\to\infty}\mathbb{E}_{n} 
\bigg[\bigg( \int_0^T \mathcal{Q}^n_s(n^{a-2} \tilde{\mathscr{D}}_n h_n) ds \bigg)^2 \bigg]
= 0.
\end{equation*}
\end{lemma}
\begin{proof}
The proof is similar to previous works in a case of the harmonic potential and it is exactly the same as ~\cite[Lemma~3.3]{bernardin2018weakly}. 
The only difference is the definition of the quadratic field $\mathcal{Q}^n$ where the variables $\xi_j$'s were replaced by $\eta_j$ in the harmonic case. 
The proof given in~\cite{bernardin2018weakly}, however, is valid also under this difference so we omit a detailed description of the proof here. 
\end{proof}

Finally, we need to show tightness of the fluctuation fields, but this is given in ~\cite[Section~5.2]{bernardin20163} so that we omit the details here.
Then the proof of Theorem~\ref{thm:3/2L_derivation} is completed.

\subsection{Quantitative Estimate for the Solution of the Poisson Equation}
\label{subsec:poisson_bound}
Here we give a proof of the estimates in Lemmas \ref{lem:poisson_sol_bound} and \ref{lem:est_fourier_1d} with Fourier analysis.
For $f:\frac{1}{n}\mathbb{Z}^d \to \mathbb{R}$, let $\widehat{f}_n:\mathbb{R}^d \to \mathbb{R}$ be the Fourier transform of $f$ defined by 
\begin{equation*}
\widehat{f}_n(k)
= \frac{1}{n^d} \sum_{j\in\mathbb{Z}^d} 
f\bigg(\frac{j}{n} \bigg)
e^{\frac{2\pi\mathsf{i} k\cdot j}{n}}. 
\end{equation*}
Here $\mathsf{i}=\sqrt{-1}$ denotes the imaginary unit. 
Then, the inverse Fourier transform is given as  
\begin{equation*}
f \bigg(\frac{j}{n} \bigg)
= \int_{[-\frac{n}{2}, \frac{n}{2}]^d} 
\widehat{f}_n (k) 
 e^{- \frac{2\pi\mathsf{i} j\cdot k}{n}}  dk .
\end{equation*}

Now, we give the proof of Lemma~\ref{lem:poisson_sol_bound}.

\begin{proof}[Proof of Lemma \ref{lem:poisson_sol_bound}]
The first assertion is shown in \cite[Appendix~D]{bernardin20163} so that we focus on the other estimates. 
Applying Fourier transform to the Poisson equation \eqref{eq:Poisson}, we have that 
\begin{equation}
\label{eq:hn_fourier_transformation}
\widehat{h}_n (k,k^\prime)
= \frac{1}{\sqrt{n}}
\frac{\mathsf{i} \Omega(\frac{k}{n}, \frac{k^\prime}{n}) }{-\alpha_n^{-1} \Lambda(\frac{k}{n},\frac{k^\prime}{n}) 
- \mathsf{i} \Omega(\frac{k}{n},\frac{k^\prime}{n})} 
\widehat{\varphi}_n (k+k^\prime), 
\end{equation}
where 
\begin{equation*}
\Lambda \bigg( \frac{k}{n}, \frac{k^\prime}{n} \bigg)
= 4 \bigg[ \sin^2\bigg(\frac{\pi k}{n}\bigg)
+ \sin^2\bigg(\frac{\pi k^\prime}{n}\bigg)
\bigg] ,\quad
\mathsf{i} \Omega \bigg( \frac{k}{n}, \frac{k^\prime}{n} \bigg)
= 2\mathsf{i} \bigg[ \sin\bigg(\frac{2\pi k}{n}\bigg)
+ \sin\bigg(\frac{2\pi k^\prime}{n}\bigg)
\bigg] .
\end{equation*}
The computation is exactly the same as \cite[Appendix C]{bernardin2018weakly}, therefore we omit details. 
Note that 
\begin{equation*}
\widehat{(\mathscr{E}_nh)}_n (k,k^\prime)
= n \big( e^{-\frac{2\pi\mathsf{i}k}{n}}-1 \big) 
\widehat{h}_n (k,k^\prime) ,\quad 
\widehat{(\mathscr{A}_nh)}_n (k,k^\prime)
=\mathsf{i}n 
\Omega \bigg( \frac{k}{n}, \frac{k^\prime}{n} \bigg)
\widehat{h}_n (k,k^\prime) .
\end{equation*}
Hence, according to the Parseval-Plancherel identity, we have that 
\begin{equation*}
\begin{aligned}
\| \mathscr{E}_n h \|^2_{2,n}
&= \iint_{[-\frac{n}{2},\frac{n}{2}]^2}
| \widehat{(\mathscr{E}_nh)}_n (k,k^\prime) |^2
dk dk^\prime \\
&= n \iint_{[-\frac{n}{2},\frac{n}{2}]^2}
\big| e^{-\frac{2\pi\mathsf{i}k}{n}} - 1 \big|^2
\frac{ \Omega(\frac{k}{n}, \frac{k^\prime}{n})^2 |\widehat{\varphi}_n (k+k^\prime)|^2 }{\alpha_n^{-2} \Lambda(\frac{k}{n},\frac{k^\prime}{n})^2 
+ \Omega(\frac{k}{n},\frac{k^\prime}{n})^2} 
dk dk^\prime\\
&= 4n \iint_{[-\frac{n}{2},\frac{n}{2}]^2}
\frac{ \sin^2 (\frac{\pi(\xi-k^\prime)}{n}) 
\Omega(\frac{\xi-k^\prime}{n}, \frac{k^\prime}{n})^2 |\widehat{\varphi}_n (\xi)|^2 }{\alpha_n^{-2} \Lambda(\frac{\xi-k^\prime}{n},\frac{k^\prime}{n})^2 
+ \Omega(\frac{\xi-k^\prime}{n},\frac{k^\prime}{n})^2} 
d\xi dk^\prime .
\end{aligned}
\end{equation*}
Here in the last estimate we used 
$|e^{-\frac{2\pi\mathsf{i}k}{n}}-1|^2
= 4\sin^2 (\pi k/n)$, and then applied a change of variables $\xi=k+k^\prime$.
Similarly, we have that 
\begin{equation*}
\| \mathscr{A}_nh\|^2_{2,n}
= n \iint_{[-\frac{n}{2},\frac{n}{2}]^2}
\frac{ \Omega(\frac{\xi-k^\prime}{n}, \frac{k^\prime}{n})^4 |\widehat{\varphi}_n (\xi)|^2 }{\alpha_n^{-2} \Lambda(\frac{\xi-k^\prime}{n},\frac{k^\prime}{n})^2 
+ \Omega(\frac{\xi-k^\prime}{n},\frac{k^\prime}{n})^2} 
d\xi dk^\prime .
\end{equation*}
Moreover, note that 
\begin{equation*}
\Omega\bigg( \frac{\xi-k^\prime}{n} , \frac{k^\prime}{n} \bigg)
\le 4 | 1-e^{\frac{2\pi\mathsf{i}\xi}{n}} |^2
= 16 \sin^2 \bigg( \frac{\pi \xi}{n}\bigg) .
\end{equation*}
Then, we have that
\begin{equation*}
\| \mathscr{E}_nh \|^2_{2,n}
\lesssim n^3
\int_{-1/2}^{1/2} 
\sin^2 (\pi y)
| \widehat{\varphi}_n (ny) |^2
\tilde{V}_n (y) dy, 
\end{equation*}
and 
\begin{equation*}
\| \mathscr{A}_nh \|^2_{2,n}
\lesssim n^3
\int_{-1/2}^{1/2} 
\sin^4 (\pi y)
| \widehat{\varphi}_n (ny) |^2
\tilde{W}_n (y) dy, 
\end{equation*}
where 
\begin{equation}
\label{eq:kernel_estimate}
\begin{aligned}
\tilde{V}_n(y)
= \int_{-1/2}^{1/2}
\frac{\sin^2(\pi(y-x))}{\alpha_n^{-2} \Lambda(y-x,x)^2 
+ \Omega(y-x,x)^2 } dx 
\le \int_{-1/2}^{1/2}
\frac{\sin^2(\pi(y-x))}{ \Lambda(y-x,x)^2 
+ \Omega(y-x,x)^2 } dx ,
\end{aligned}
\end{equation}
and 
\begin{equation*}
\begin{aligned}
\tilde{W}_n(y)
= \int_{-1/2}^{1/2}
\frac{dx}{\alpha_n^{-2} \Lambda(y-x,x)^2 
+ \Omega(y-x,x)^2 } 
\le \int_{-1/2}^{1/2}
\frac{dx}{ \Lambda(y-x,x)^2 
+ \Omega(y-x,x)^2 } .
\end{aligned}
\end{equation*}
Similarly to the estimate of $h_n$ which is given in~\cite[Appendix C]{bernardin2018weakly}, it suffices to get a bound for $\tilde{V}_n$ by a polynomial function.  
Indeed, suppose $\tilde{V}_n(y)=O(|y|^{q})$ with some $q>0$.
Moreover, note that for $\varphi\in \mathcal{S}(\mathbb{R})$, we have 
\begin{equation*}
| \widehat{\varphi}_n(yn) |
\lesssim \frac{1}{1 + (n|y|)^p}  ,
\end{equation*}
for any $p\ge 1$.
(See \cite[Lemma B.1]{bernardin20163}.)  
Thus, with these estimates at hand we have 
\begin{equation*}
\begin{aligned}
\| \mathscr{E}_nh \|^2_{2,n}
&\lesssim n^3 \int_{-1/2}^{1/2}
\frac{|y|^{2+q}}{(1+(n|y|)^p)^2} dy
\lesssim n^3 \int_{-1/2}^{1/2}
\frac{|y|^{2+q}}{1+(n|y|)^{2p}} dy \\
&\lesssim n^{-q} \int_\mathbb{R}
\frac{|y|^{2+q}}{1 + |y|^{2p}} dy
= O(n^{-q}),
\end{aligned}
\end{equation*}
by taking $p\ge 1$ sufficiently large. 
In addition, Lemma~\ref{lem:residue_bound} below enables us to take $q=0$, by which we complete the proof of the bound for $\mathscr{E}_nh$.  
On the other hand, we have $\tilde{W}_n(y)=|y|^{-3/2}$ by \cite[Lemma F.5]{bernardin20163}, which yields 
\begin{equation*}
\begin{aligned}
\| \mathscr{A}_nh \|^2_{2,n}
\lesssim n^3 \int_{-1/2}^{1/2}
\frac{|y|^{4-3/2}}{1+(n|y|)^{2p}} dy
\lesssim n^{-1/2} \int_\mathbb{R}
\frac{|y|^{q}}{1 + |y|^{2p}} dy
= O(n^{-1/2}).
\end{aligned}
\end{equation*}
Hence we complete the proof of Lemma \ref{lem:poisson_sol_bound}. 
\end{proof}

\begin{lemma}
\label{lem:residue_bound}
There exists some $C>0$ such that $\tilde{V}_n(y) \le C$ for any $y\in \mathbb{R}$.
\end{lemma}
\begin{proof}
Recall the estimate for $\tilde{V}_n$ given in~\eqref{eq:kernel_estimate}. 
This is further bounded as follows. 
\begin{equation*}
\tilde{V}_n(y)
\lesssim \int_{-1/2}^{1/2} 
\frac{\sin(\pi(y-x))}{\sin(\pi(y-x)) + 2 \cos(\pi(y-x))} dx
= V_n(y) .
\end{equation*}
In what follows we give an estimate for $V_n$.
Set $z=e^{2\pi\mathsf{i}x}$ and $w=e^{2\pi\mathsf{i}y}$ for $x,y\in[-1/2,1/2] $. 
Note that 
\begin{equation*}
\begin{aligned}
\sin (\pi(y-x))
= \frac{1}{2\mathsf{i}}(wz^{-1}-w^{-1}z), \quad
\cos (\pi(y-x))
= \frac{1}{2}(wz^{-1}+w^{-1}z). 
\end{aligned}
\end{equation*}
Then, we can write 
\begin{equation*}
V_n(y)
= \frac{1}{2\pi\mathsf{i}}\oint_{\mathscr{C}} g_w(z) dz ,
\end{equation*}
where $\mathscr{C}$ denotes the unit circle which is positively oriented and the meromorphic function $g_w$ is defined by
\begin{equation*}
g_w(z)
= \frac{1}{z}
\frac{wz^{-1}-w^{-1}z}{(wz^{-1}-w^{-1}z) + 2\mathsf{i}(wz^{-1}+w^{-1}z)} 
= \frac{1}{z}
\frac{z^2-w^2}{(1-2\mathsf{i})z^2- (1+2\mathsf{i})w^2 }. 
\end{equation*}
Note that the meromorphic function $g_w$ has three poles inside the unit circle $\mathscr{C}$: $z=0, \pm a_0w$ with $a_0=(1+2\mathsf{i})(1-2\mathsf{i})$. 
Then, by the residue theorem, we have 
\begin{equation*}
V_n(y) = 
\mathrm{Res}(g_w,0)
+ \mathrm{Res}(g_w,a_0w)
+ \mathrm{Res}(g_w,-a_0w), 
\end{equation*}
where $\mathrm{Res}(g_w,\cdot)$ denotes the value of residue of $g_w$ at each pole in $\mathbb{C}$. 
Noting  
\begin{equation*}
\mathrm{Res}(g_w,0) = \frac{1}{1+2\mathsf{i}} ,
\end{equation*}
and
\begin{equation*}
\mathrm{Res}(g_w,\pm a_0w) 
= \frac{1}{aw}
\frac{(a_0^2-1)w^2}{(1-2\mathsf{i}) (a_0w)(2a_0w)}
= \frac{a_0^2-1}{2a_0^2(1-2\mathsf{i})},
\end{equation*}
we see that the summation of the residues inside the circle $\mathscr{C}$ is independent of $w$. 
This means $\tilde{W}_n(y)=O(1)$ and we complete the proof. 
\end{proof}

Finally, we give a proof of Lemma \ref{lem:est_fourier_1d}.

\begin{proof}[Proof of Lemma \ref{lem:est_fourier_1d}]
In the sequel, we simply write $h_n =h$ and $h_{j,j^\prime}=h(j/n,j^\prime/n)$ with an abuse of notation. 
Set $u_{j,j^\prime}=h_{j+j^\prime,j-j^\prime+1}$. 
Then, the Fourier transformation of $u$ is computed as 
\begin{equation*}
\begin{aligned}
\widehat{u}(\xi,\xi^\prime)   
&= \frac{1}{n^2}\sum_{j,j^\prime\in\mathbb{Z}}
h_{j+j^\prime,j-j^\prime+1} 
e^{\frac{2\pi \mathsf{i}}{n}(j\xi + j^\prime \xi^\prime)}
= \frac{1}{n^2}\sum_{k,k^\prime\in\mathbb{Z}}
h_{k,k^\prime+1} 
e^{\frac{\pi \mathsf{i}}{n}( k(\xi+\xi^\prime) + k^\prime(\xi-\xi^\prime))} \\
&= e^{-\frac{\pi \mathsf{i}}{n}(\xi-\xi^\prime) } 
\widehat{h}\bigg( \frac{\xi+\xi^\prime}{2},\frac{\xi-\xi^\prime}{2}\bigg) . 
\end{aligned}
\end{equation*}
Noting $\int_{[-n/2,n/2]} e^{2\pi \mathsf{i}j\xi/n}d\xi= n\delta_{j,0}$, we apply the integration over $\xi^\prime$-variable in the first identity of the last display. 
Then we have that 
\begin{equation*}
\int_{[\frac{n}{2},\frac{n}{2}]}
\widehat{u}(\xi,\xi^\prime)d\xi^\prime
= \frac{1}{n} \sum_{j,j^\prime\in\mathbb{Z}} 
h_{j+j^\prime,j-j^\prime+1} e^{\frac{2\pi\mathsf{i}}{n}j\xi} \delta_{j^\prime,0}
= \frac{1}{n} \sum_{j\in\mathbb{Z}} 
\Erase{h_{j,j^\prime+1}}\Add{h_{j,j+1}} e^{\frac{2\pi\mathsf{i}}{n}j\xi}
= \widehat{h\mathbf{1}_{j^\prime=j+1}}(\xi). 
\end{equation*}
Hence by the Parseval-Plancherel identity, we have that 
\Add{
\begin{equation*}
\frac{1}{n} \sum_{j\in\mathbb{Z}}
h_{j,j+1}^2
= \int_{[\frac{n}{2},\frac{n}{2}]}
\big|\widehat{h\mathbf{1}_{j^\prime=j+1}}(\xi)\big|^2 d\xi
= \int_{[\frac{n}{2},\frac{n}{2}]}
\bigg( \int_{[\frac{n}{2},\frac{n}{2}]}
\widehat{u}(\xi,\xi^\prime)d\xi^\prime\bigg)^2
d\xi.
\end{equation*}
}
Moreover, recall that the Fourier transformation of $h_n$ is given in~\eqref{eq:hn_fourier_transformation}, which yields,  
\begin{equation*}
\widehat{u}(\xi,\xi^\prime)
= e^{-\frac{\pi \mathsf{i}}{n}(\xi-\xi^\prime) }
\widehat{h}\bigg(\frac{\xi+\xi^\prime}{2},\frac{\xi-\xi^\prime}{2} \bigg)
= \frac{1}{\sqrt{n}}
\frac{\mathsf{i} e^{-\frac{\pi \mathsf{i}}{n}(\xi-\xi^\prime) } \Omega(\frac{\xi+\xi^\prime}{2n}, \frac{\xi-\xi^\prime}{2n}) }{-\alpha_n^{-1} \Lambda(\frac{\xi+\xi^\prime}{2n},\frac{\xi-\xi^\prime}{2n}) 
- \mathsf{i} \Omega(\frac{\xi+\xi^\prime}{2n},\frac{\xi-\xi^\prime}{2n})} 
\widehat{\varphi}_n (\xi). 
\end{equation*}
Recalling the definition of $\Omega$ and $\Lambda$, we have a bound $|\widehat{u}(\xi,\xi^\prime)|^2
\le n^{-1} | \widehat{\varphi_n}(\xi)|^2$. 
Therefore,  
\Add{
\begin{equation*}
\frac{1}{n} \sum_{j\in\mathbb{Z}}
h_{j,j+1}^2 
\le \int_{[-\frac{n}{2},\frac{n}{2}]}
|\widehat{\varphi}_n(\xi)|^2 
d\xi
\lesssim  \int_{\mathbb{R}} \frac{d\xi}{1+|\xi|^2}, 
\end{equation*}
}
which is bounded by a universal constant. 
This completes the proof. 
Here we used the bound $| \widehat{\varphi_n}(\xi)|\lesssim (1+|\xi|^p)^{-1}$ for any $\varphi\in \mathcal{S}(\mathbb{R})$ and $p \ge 1$. 
\end{proof}

\appendix

\section{Heuristics from Mode Coupling Theory}

The NLFH theory developed in \cite{spohn2014nonlinear} for a model of chains of oscillators, gives a whole scenario of possible limits for fluctuation fields of the conserved quantities of those   systems. In those models, one can observe non-trivial nonlinear couplings which might happen  in different time scales. 
In~\cite{spohn2015nonlinear} the authors applied the NLFH theory of \cite{spohn2014nonlinear} to the BS model we investigate in this article but for  a fixed potential $V(\cdot)$ and  in this section, we recap those predictions. 
The arguments can be  briefly described as follows.
Following the setting of \cite{spohn2015nonlinear}, we rewrite the invariant measure of the system as 
\begin{equation*}
\mu_{\tau,b}(\eta_j)
= \frac{1}{Z_{\tau,b}} \exp \big( - b(V_{\beta}(\eta_j)+ \tau \eta_j ) \big)d\eta_j .
\end{equation*}
Note that our parameter $\lambda$ corresponds to $-b\tau$ in the measure $\mu_{\tau,b}$. 
Now, we choose $\tau=\tau(v,e), b=b(v,e)$ in such a way that 
\begin{equation*}
E_{\mu_{\tau,b}}[\eta_j]= v, \quad 
E_{\mu_{\tau,b}}[\zeta_j]= e.
\end{equation*}
Then, we write $\overline{\mu}_{v,e}=\mu_{\tau(v,e), b(v,e)}$. 
A simple computation shows that 
\begin{equation*}
L_n \eta_j = j^v_{j-1,j} - j^v_{j,j+1}, \quad 
L_n \zeta_j = j^e_{j-1,j} - j^e_{j,j+1},
\end{equation*}
where 
\begin{equation*}
\begin{aligned}
& j^v_{j,j+1} = \frac{1}{2} \theta(n)(\eta_j - \eta_{j+1}) 
- \theta(n) \alpha_n (\xi_j + \xi_{j+1}), \\
& j^e_{j,j+1} = \frac{1}{2} \theta(n)(\zeta_j - \zeta_{j+1})  - \theta(n)\alpha_n \xi_j \xi_{j+1} .
\end{aligned}
\end{equation*}
For simplicity of notation let us drop down the time scale $\theta(n) $ and the asymmetry $\alpha_n$. 
Note that $E_{\overline{\mu}_{v,e}}[j^v_{j,j+1}]=2\tau(v,e) $ and  $E_{\overline{\mu}_{v,e}}[j^e_{j,j+1}]=-\tau(v,e)^2$.

The column flux matrix is given by
\begin{equation*}
j(v,e)= \begin{pmatrix}
E_{\overline{\mu}_{v,e}}[j^v_{j,j+1}] \\
E_{\overline{\mu}_{v,e}}[j^e_{j,j+1}] 
\end{pmatrix}
= 
\begin{pmatrix}
2\tau \\
- \tau^2 
\end{pmatrix} ,
\end{equation*}
and the Jacobian matrix is equal to 
\begin{equation*}
J = \begin{pmatrix}
\partial_v E_{\overline{\mu}_{v,e}}[j^v_{j,j+1}]
& \partial_e E_{\overline{\mu}_{v,e}}[j^v_{j,j+1}] \\
\partial_v E_{\overline{\mu}_{v,e}}[j^e_{j,j+1}]
& \partial_e E_{\overline{\mu}_{v,e}}[j^e_{j,j+1}] 
\end{pmatrix}=
2\begin{pmatrix}
\partial_v \tau
& \partial_e \tau \\
-\tau \partial_v \tau 
& -\tau \partial_e \tau 
\end{pmatrix} .
\end{equation*}
It is simple to check that the matrix $J$ has two eigenvalues $v_+=0$ and $v_-=2(\partial_v \tau -\tau \partial_e\tau)<0$, see \cite{spohn2015nonlinear}. 
The corresponding eigenvectors are proportional to  $u_+=(\partial_e \tau, -\partial_v\tau)$ and $u_-=(1,-\tau)$. 
To obtain the normal modes we need to find the matrix $R$ that diagonalizes $J$, namely, $RJR^{-1}={\rm diag}(v_+, v_-)$. Note that $R^{-1}$ is the matrix whose columns are the eigenvectors of $J$:
\begin{equation}
R^{-1}=\begin{pmatrix}\partial_e \tau&1 \\
-\partial_v\tau& -\tau
\end{pmatrix},
\end{equation}
and its inverse is equal to  
\begin{equation}
R=\frac{1}{\partial_ v\tau - \tau \partial_e \tau  }\begin{pmatrix}-\tau &-1 \\
\partial_ v\tau&\partial_e \tau 
\end{pmatrix}.
\end{equation}
According to NLFH, the quantities that we should look at are equal to $R (\bar \eta_j, \bar \zeta_j),$
which gives 
\begin{equation}\label{normal_modes}
\begin{aligned}
\phi_+&=\frac{1}{\partial_ v\tau-\tau \partial_e \tau}
(-\tau \bar{\eta}_j-\bar\zeta_j ),\\
\phi_-&=\frac{1}{\partial_ v\tau-\tau \partial_e \tau}\big((\partial_v \tau )\bar{\eta}_j+(\partial_e\tau)\bar{\zeta}_j\big).
\end{aligned}
\end{equation}
Here recall the Taylor expansion \eqref{eq:taylor_expansion_upto_beta_squared} and take the expectation. 
Then, we have 
\begin{equation*}
-\tau = v + c_3 \beta_n e + O(\beta_n^2). 
\end{equation*}
Hence \eqref{normal_modes} explains the quantities that we obtained by the martingale decomposition argument \Add{in Section \ref{sec:cancellation}}, up to a scale factor and an error term with order $O(\beta_n^2)$.

In analogy to anharmonic chains of oscillators the first mode is called the sound mode and the second mode is called the heat mode. 
Remarkably for the BS model one can compute the precise values of the quantities above in terms of the cumulants of the variables involved, but explicit expressions are difficult to obtain. For details we refer the reader to \cite[Appendix C]{spohn2015nonlinear}. 

Now that we have the fluctuations fields for the normal modes, which are fixed, let us see the predictions on the form of the fluctuations for each one of these quantities. 
To do that we look now at the corresponding Hessians of the entries of the Jacobian matrix $J$:

\begin{equation*}
 H^1= \begin{pmatrix}
\partial^2_v E_{\overline{\mu}_{v,e}}[j^v_{j,j+1}]
&\partial_v \partial_e E_{\overline{\mu}_{v,e}}[j^v_{j,j+1}] \\
\partial_v \partial_e E_{\overline{\mu}_{v,e}}[j^v_{j,j+1}]
& \partial^2_e E_{\overline{\mu}_{v,e}}[j^v_{j,j+1}] 
\end{pmatrix}
=2\begin{pmatrix}
\partial^2_v \tau
& \partial_v \partial_e \tau \\
\partial_v \partial_e \tau
& \partial^2_e \tau 
\end{pmatrix}
\end{equation*}

\begin{equation*}
 H^2= \begin{pmatrix}
\partial^2_v E_{\overline{\mu}_{v,e}}[j^e_{j,j+1}]
&\partial_v \partial_e E_{\overline{\mu}_{v,e}}[j^e_{j,j+1}] \\
\partial_v \partial_e E_{\overline{\mu}_{v,e}}[j^e_{j,j+1}]
& \partial^2_e E_{\overline{\mu}_{v,e}}[j^e_{j,j+1}] 
\end{pmatrix}=\tau H^1-2\begin{pmatrix}
(\partial_v \tau)^2
&( \partial_v \tau)(\partial_e \tau) \\
( \partial_v \tau)(\partial_e \tau)
& (\partial_e \tau)^2 
\end{pmatrix}
\end{equation*}

\noindent 
The coupling constants, which are  determined by the above matrices, are given on $i\in\{1,2\}$ by
$$G^i=\frac{1}{2} \sum_{j=1}^2 R_{i,j}[(R^{-1})^\top  H^j R^{-1}],
$$
where $R_{i,j}$ is the entry of the matrix $R$.

From the particular expressions of the instantaneous currents,  one can check that the only non-zero entry of the mode coupling matrix $G^2$ is $G^2_{11}$, i.e. $G^2_{22}=0$ and $G^2_{12}=G^2_{21}=0$ and, in fact, $G^2_{11}<0$. For the matrix $G^1$ all the entries are, a priori, non-zero. 
The classification for the limiting processes that we should get for a microscopic system with two conservation laws is  summarized in the tables below, which are taken from  \cite{spohn2015nonlinear}. Since $G^2_{22}=G^2_{12}=G^2_{21}=0$ and $G^2_{11}<0$, the only possible pairs that we can obtain for the BS model are: 
\begin{itemize}
\item[(A)] 
if $G^1_{11}\neq 0$ then we get the pair \textcolor{blue}{(KPZ\,,\,$5/3$-L\'evy)}.
\item[(B)] 
if $G^1_{11}= 0$ then we get the pair \textcolor{red}{(Diff\,,\,$3/2$-L\'evy)}.
\end{itemize}

\begin{table}[htpb]
\begin{center}
\caption{Classification (I). }
\label{tab:mcm_1}
\begin{tabular}{ |c | c | c | c| }
\hline
$G^1_{11}=G^2_{22}=1$ & $G^1_{22}$  & $G^2_{11}$ & (Mode 1, Mode 2) \\
\hline
& 0,1 & 0,1& (KPZ,KPZ)  \\
\hline
\end{tabular}
\end{center}
\end{table}

\begin{table}[htpb]
\begin{center}
\caption{Classification (II). }
\label{tab:mcm_2}
\begin{tabular}{ |c | c | c | c| }
\hline
$G^1_{11}=1$, $G^2_{22}=0$ & $G^1_{22}$  & $G^2_{11}$ & (Mode 1, Mode 2) \\
\hline 
& 0,1 & 1& \textcolor{blue}{(KPZ\,,\,$5/3$-L\'evy)}  \\
\hline
& 1 & 0& (mod. KPZ\,,\,Diff)  \\
\hline
& 0 & 0& (KPZ\,,\,Diff)  \\
\hline
\end{tabular}
\end{center}
\end{table}

\begin{table}[htpb]
\caption{Classification (III). }
\label{tab:mcm_3}
\begin{center}
	\begin{tabular}{ |c | c | c | c| }
		\hline
		$G^1_{11}=0$, $G^2_{22}=0$ & $G^1_{22}$  & $G^2_{11}$ & (Mode 1, Mode 2) \\
		\hline 
		& 1 & 1& (Gold-L\'evy\,,\, Gold-L\'evy)  \\
		\hline
		& 1 & 0& ($3/2$-L\'evy \,,\,Diff)  \\
		\hline
		& 0 & 1& \textcolor{red}{(Diff\,,\,$3/2$-L\'evy)}\\
		\hline
		& 0 & 0& 	(Diff\,,\,Diff)  \\
		\hline
	\end{tabular}
\end{center}
\end{table}

In \cite{spohn2015nonlinear} several numerical simulations \Erase{haven}\Add{have} been done for specific potentials. 
In the case of the FPU potential $V(\eta)= \eta^2/2+\alpha\eta^3/3+\eta^4/4$ the two modes should behave as case (A) i.e.  \textcolor{blue}{(KPZ\,,\,$5/3$-L\'evy)} and the same should occur for the Toda lattice potential $V(\eta)=e^{-\eta}-1+\eta$. For the harmonic potential $V(\eta)=\eta^2/2$, one should observe the behavior of case (B), i.e. \textcolor{red}{(Diff\,,\,$3/2$-L\'evy)} and this has been rigorously proved in \cite{bernardin20163, bernardin2018weakly}.

In the previous tables there is no pair of universality classes of the form (KPZ\,,\,$3/2$-L\'evy), or the exchanged case ($3/2$-L\'evy\,,\,KPZ).
The classification which is based on coupling matrices is, however, conducted for a fixed potential and in the strong asymmetric regime, corresponding in our case to $\alpha_n=O(1)$. 
In particular, our current oscillator model belongs to different universality classes if the driving potential varies.

We finish this digression with a simple case corresponding to $\beta_n=0$, for which we can make precise computations. 
Note that if $\beta_n\to0$, then $\nu_n \to \mathcal{N}(\tau/b, 1/b)$ and $E_{\nu_n}[\xi_j]=\tau/b$ so that 
\begin{equation*}
v = \tau/b ,\quad 
e = (1/2)(1/b + (\tau/b)^2 ).
\end{equation*}
Hence, $b=(2e-v^2)^{-1}$ and $\tau = v(2e-v^2)^{-1}$. 
Then, 
\begin{equation*}
j(v,e) = \begin{pmatrix}
2\theta(n) \alpha_n v \\
\theta(n) \alpha_n v^2  
\end{pmatrix}
, \quad 
J= \begin{pmatrix}
2\theta(n) \alpha_n & 0 \\
2\theta(n)\alpha_n v & 0
\end{pmatrix}.
\end{equation*}
The 
eigenvalues of the matrix $J$ correspond to $v_+=2\theta(n)\alpha_n$ and $v_-=0$,
whose eigenvectors are proportional to $u_+=(1,v_+)$ and $v_-=(0,1)^\top$, respectively.  
Hence, we can diagonalize the matrix $J$ as 
\begin{equation*}
R J R^{-1} = \mathrm{diag}[2\theta(n)\alpha_n, 0], \quad 
R = \begin{pmatrix}
1 & 0 \\
-v_+ & 1 
\end{pmatrix},\quad R^{-1}=
\begin{pmatrix}
1 & 0 \\
v_+ & 1
\end{pmatrix}.
\end{equation*}
Therefore, the linear combination of the  fields that we should look at  are associated to the quantities $(\phi_+, \phi_-) =R(\overline{\eta}_j, \overline{\zeta}_j)
= (\overline{\eta}_j ,-v\overline{\eta}_j  + \overline{\zeta}_j)$ and these are exactly the quantities what we obtained in Section \ref{sec:cancellation} with $\beta_n=0$.

Now that we have the fluctuation fields fixed, let us see the predictions on the form of the fluctuations for each one of these quantities. 
To do that we look now at the corresponding Hessians of the entries of the Jacobian matrix $J$:
\begin{equation}\label{eq:hessians}
 H^1=\begin{pmatrix}
0& 0\\
0& 0
\end{pmatrix} \quad \textrm{and }\quad 
 H^2=\begin{pmatrix}
2\theta_n\alpha_n& 0\\
0& 
0 
\end{pmatrix} 
\end{equation}
A simple computation shows that 
$(R^{-1})^\top  H^1 R^{-1}$ corresponds to the null matrix while 
$(R^{-1})^\top H^2 R^{-1}= H^2$.
From this  we get 
	$$G^1=\begin{pmatrix}
	0 & 0\\0& 0
	\end{pmatrix} \quad \textrm{and}\quad G^2=\begin{pmatrix}
		\theta_n\alpha_n & 0\\0 & 0
	\end{pmatrix}.
$$
Comparing to the results of the table above we expect to observe, in the strong asymmetric regime corresponding to $\alpha_n=O(1)$ the behavior 
\textcolor{red}{(Diff\,,\,$3/2$-L\'evy)}, i.e., diffusive and fractional behavior.

\section*{Acknowledgments}
The authors thank Makiko Sasada, Gunter M. Sch\"{u}tz and Hayate Suda for their interest in this work and for giving fruitful comments. K.H. was supported by JSPS KAKENHI Grant Number JP22J12607.
P.G. thanks  FCT/Portugal for financial support through the
projects UIDB/04459/2020 and UIDP/04459/2020.  This project has received funding from the European Research Council (ERC) under  the European Union's Horizon 2020 research and innovative program (grant agreement   n. 715734).
Data sharing not applicable to this article as no datasets were generated or analysed during the current study.

\bibliographystyle{abbrv}
\bibliography{ref}

\end{document}